\documentclass{amsart}
\usepackage{amssymb}
\usepackage{amsmath}
\usepackage{amsthm}
\usepackage{mathrsfs}
\usepackage{verbatim}
\usepackage{xcolor}
\usepackage{ulem}
\usepackage{mathtools}
\usepackage{physics}
\usepackage{enumerate}

\newtheorem{theorem}{Theorem}[section]

\newtheorem{lemma}[theorem]{Lemma}
\newtheorem{proposition}[theorem]{Proposition}
\newtheorem{corollary}[theorem]{Corollary}
\newtheorem*{corollary*}{Corollary}
\theoremstyle{definition}
\newtheorem{definition}[theorem]{Definition}

\theoremstyle{remark}
\newtheorem{remark}[theorem]{Remark}
\numberwithin{equation}{section}
\title[]{A Dirichlet form approach to the potential theory of $m$-subharmonic functions}

\newcommand{\N}{\mathbb{N}}
\newcommand{\R}{\mathbb{R}}
\newcommand{\C}{\mathbb{C}}
\newcommand{\Q}{\mathbb{Q}}
\newcommand{\mc}[1]{\mathcal{#1}}
\newcommand{\supp}[1]{\operatorname{supp}[#1]}

\author[]{Bowoo Kang, Jaechul Lee, and Zhenghao Li}
\date{\today}
\begin{document}
\begin{abstract}
    On a bounded open subset of $\C^n$, we study the Dirichlet form induced by an admissible pair $(\theta,\mu)$, where $\mu$ is a positive Radon measure and
    \[
    \theta=dd^c u_1\wedge\cdots\wedge dd^c u_{m-1}\wedge(dd^c \lvert z\rvert^2)^{n-m}
    \]
    for locally bounded $m$-subharmonic functions $u_1,\dots,u_{m-1}$. We prove that locally bounded $m$-subharmonic functions belong locally to the domains of these Dirichlet forms and are quasi-continuous with respect to the induced capacities. We establish a monotonicity property for Dirichlet forms, which motivates the introduction of $\mu$-quasi-monotone convergence of $m$-subharmonic functions. The natural embeddings into these domains are continuous with respect to this convergence. These results yield new proofs of the commutativity and continuity of complex Hessian operators.

    Using this framework, we characterize quasi-continuity with respect to the induced capacity for an arbitrary $m$-subharmonic function in terms of its pole set having zero measure $\mu$, or equivalently, zero induced capacity. On bounded strictly $m$-pseudoconvex domains, we prove that the Bedford--Taylor capacity for $m$-subharmonic functions and the capacities induced by the corresponding Dirichlet forms have the same null sets. We further characterize $m$-polar sets in terms of Dirichlet forms and stochastic processes, extending the corresponding characterization of pluripolar sets due to Fukushima and Okada to the $m$-subharmonic setting.
\end{abstract}

\maketitle
\tableofcontents

\section{Introduction}
Let $\Omega \subset \C^n$ be a bounded open subset and $\mu$ be a positive Radon measure on $\Omega$ satisfying $\supp\mu=\Omega$. A Dirichlet form on the real space $L^2(\Omega; \mu)$ is a non-negative definite symmetric bilinear form on a dense linear subspace of $L^2(\Omega; \mu)$ that is Markovian and closed. Fukushima and Okada \cite{FO87} showed that, by considering the Dirichlet form associated with a  closed positive $(n-1, n-1)$-current on $\Omega$, one can study plurisubharmonic functions and complex Monge--Amp\`ere operators together with a probabilistic interpretation of pluripolar sets. 

The modern theory of Dirichlet forms was initiated by Beurling and Deny \cite{Beurling59,Deny70} and developed substantially by Fukushima, who established a correspondence between Dirichlet forms and symmetric Markov processes \cite{Fu71}. Silverstein \cite{Sil74} subsequently provided another proof of this correspondence using a more direct approach. Since this correspondence reveals a deep connection between analytic and probabilistic potential theory, Dirichlet forms have become an important tool in various areas, including geometric analysis and stochastic analysis, see, e.g., \cite{Ro93} and the references therein. Owing to the flexibility of the Dirichlet form method, Vigny \cite{Vi07} adapted it, without invoking probabilistic methods, to develop the complex Sobolev space theory of Dinh and Sibony. Introduced in \cite{DS06}, this theory has found various applications in complex Monge--Amp\`ere equations and complex dynamics, see, e.g., \cite{Vi15, Vu20, DKN22, DMV23}. More recently, Nguyen \cite{Nguyen26} extended the construction of a functional capacity in \cite{Vi07} to the local setting. It was also shown in \cite{Vi07, Nguyen26} that the complex Sobolev space can be embedded into a Dirichlet space associated with a closed positive $(n-1, n-1)$-current.

Let $\theta$ be a closed positive $(n-1, n-1)$-current on $\Omega$. The non-negative definite symmetric bilinear form $\mathcal{E}^{\theta}$ on $\mathscr{C}_c^{\infty}(\Omega, \R)$ induced by $\theta$ is defined by
\begin{align}\label{equation:Etheta}
    \mathcal{E}^{\theta}(\varphi, \psi) := \int_{\Omega}d\varphi \wedge d^c\psi \wedge \theta,\quad \varphi, \psi \in \mathscr{C}_c^{\infty}(\Omega, \R).
\end{align}
This symmetric form $\mathcal{E}^\theta$ generalizes the Dirichlet integral from classical potential theory. If $\mathcal{E}^\theta$ is closable, then its smallest closed extension is a Dirichlet form on $L^2(\Omega;\mu)$. Using this Dirichlet form, Okada and Fukushima established various properties of plurisubharmonic functions in a series of works \cite{Okada82,Oka83, FO84,Fu86a, Fu86b, FO87,Fukushima88}. For instance, Okada \cite{Okada82} constructed the completion of $\mathscr{C}_c^\infty(\Omega,\R)$ with respect to the norm induced by $\mathcal{E}^{\theta}$, while Fukushima \cite{Fu86b} established a stochastic minimum principle for complex Monge--Amp\`ere operators. In their work \cite{FO87}, Fukushima and Okada developed a probabilistic counterpart of the Bedford--Taylor theory \cite{BT82}, providing quantitative proofs of several fundamental results, including the commutativity and continuity of complex Monge--Amp\`ere operators, as well as the comparison theorem and minimum principle. They also obtained a stochastic characterization of pluripolar sets.

The notion of $m$-subharmonic functions unifies and generalizes the notions of subharmonic functions and plurisubharmonic functions, while $m$-subharmonicity coincides with subharmonicity when $m=1$ and with plurisubharmonicity when $m=n$. We observe that the Dirichlet form method developed by Fukushima and Okada \cite{FO87} for the study of plurisubharmonic functions is sufficiently flexible to suggest its applicability to $m$-subharmonic functions, complex Hessian operators, and $m$-polar sets.

G\aa rding's work \cite{Garding} on hyperbolic polynomials laid the foundation for the theory of $m$-subharmonic functions. Let $S_m$ denote the $m$-th elementary symmetric polynomial on $\R^n$. G\aa rding proved that $S_m$ is hyperbolic at every point in $\Gamma_m$, where $\Gamma_m$ is the connected component of $\{S_m\neq 0\}$ containing $(1,\dots,1)$. One has
\begin{align*}
    \Gamma_m=\{S_1>0\}\cap\cdots\cap\{S_m>0\}, \quad \overline{\Gamma_m}=\{S_1\geq0\}\cap\cdots\cap\{S_m\geq0\}.
\end{align*}
Let $u\in \mathscr{C}^2(\Omega,\R)$. We denote by $H_u$ the complex Hessian matrix of $u$, and by $\lambda(H_u)\in \R^n$ its eigenvalues. We say that $u$ is $m$-subharmonic on $\Omega$ if $\lambda(H_u)\in \overline{\Gamma_m}$ on $\Omega$, or equivalently, the complex Hessian operator
\begin{align*}
    (dd^cu)^{k} \wedge \omega^{n-k}=C(n,k) S_k\left(\lambda(H_u)\right)dV \geq 0
\end{align*}
for $k=1,\dots, m$. Here, $\omega := dd^c\lvert z\rvert^2$, $V$ denotes the Lebesgue measure on $\C^n$, and $C(n,k)>0$ is a constant depending only on $n$ and $k$. Using G\aa rding's inequality, B\l ocki \cite{Blocki05} extended the notion of $m$-subharmonicity to upper semi-continuous functions. Following the idea of Bedford and Taylor \cite{BT82}, by induction,
\begin{align*}
    dd^cu_k \wedge dd^cu_{k-1} \wedge \cdots \wedge dd^cu_1 \wedge \omega^{n-m} := dd^c(u_kdd^cu_{k-1}\wedge \cdots \wedge dd^cu_1 \wedge \omega^{n-m})
\end{align*}
is a closed positive $(n-m+k,n-m+k)$-current for $1\leq k\leq m$ and $u_1,\dots,u_k\in SH_m(\Omega)\cap L_{loc}^\infty(\Omega)$. It also coincides with the wedge product of continuous forms when $u_1,\dots,u_k \in SH_m(\Omega)\cap \mathscr{C}^2(\Omega,\R)$. Moreover, the Bedford--Taylor capacity $\mathrm{Cap}_{BT, m}$ for $m$-subharmonic functions, $m$-relative extremal functions $u_{E, \Omega, m}$, $m$-negligible sets, and $m$-polar sets are well-defined and enjoy properties analogous to those in the theory of plurisubharmonic functions. Many important results have been established in this direction, see, e.g., \cite{Blocki05,Lu12,Nguyen13,Lu13Liouville,DK14,KN16,DK17,GN18,KN26} and the references therein.

We now introduce the setting of this paper and state our main results, which are divided into three parts. For $\theta\in \mathcal{M}_m(\Omega)$, we consider the symmetric form $\mathcal{E}^\theta$ on $L^2(\Omega;\mu)$ defined by (\ref{equation:Etheta}), where
\begin{align*}
    \mathcal{M}_m(\Omega) := \{dd^cu_1 \wedge \cdots \wedge dd^cu_{m-1} \wedge \omega^{n-m}:  u_1, \dots, u_{m-1}\in SH_m(\Omega)\cap L_{loc}^\infty(\Omega)\}.
\end{align*}
We say that $(\theta,\mu)$ is an admissible pair if $\mathcal{E}^\theta$ is closable. In this case, we denote by $\mathcal{F}^\theta$ the domain of the smallest closed extension of $\mathcal{E}^\theta$. Then, $\mathcal{F}^\theta$ is a Hilbert space with the inner product $\mathcal{E}^\theta_1$,
\begin{align*}
    \mathcal{E}^\theta_1 (u,v):= \mathcal{E}^\theta(u,v)+\int_\Omega uv \ d\mu, \quad  u,v\in \mathcal{F^\theta},
\end{align*}
and $(\mathcal{F}^\theta,\mathcal{E}^\theta)$ is a regular Dirichlet space on $L^2(\Omega;\mu)$.
Thus, the potential theory of Dirichlet forms applies to our setting. More precisely, $(\mathcal{F}^\theta,\mathcal{E}^\theta)$ naturally induces a Choquet capacity $\mathrm{Cap}_\theta$ on $\Omega$, with respect to which the notions of $\mathcal{E}^\theta$-polar sets and $\mathcal{E}^\theta$-quasi-continuity are understood.

Our first objective, carried out in Section $3$, is to study the basic relationship between the Dirichlet space $(\mathcal{F}^\theta,\mathcal{E}^\theta)$ and $m$-subharmonic functions, following the strategy of \cite[Section 2]{FO87} for plurisubharmonic functions.

We first establish a Poincar\'e-type inequality for $m$-subharmonic functions and then use it to prove that $(\theta,V)$ is an admissible pair for every $\theta\in \mathcal{M}_m^+(\Omega)$. Here, the subset $\mathcal{M}_m^+(\Omega)$ of $\mathcal{M}_m(\Omega)$ is defined by
\begin{align*}
    \mathcal{M}_m^+ (\Omega):=&\{dd^cu_1 \wedge \cdots \wedge dd^cu_{m-1} \wedge \omega^{n-m}: u_1, \dots, u_{m-1}\in SH_m^+(\Omega)\},\\
    SH_m^+(\Omega):=&\{v+\delta |z|^2:  v\in SH_m(\Omega)\cap L^\infty(\Omega),   \delta>0 \}.
\end{align*}
We then establish the relationship between $(\mathcal{F}^\theta,\mathcal{E}^\theta)$
and $m$-subharmonic functions. For every admissible pair $(\theta,\mu)$, we prove that locally bounded $m$-subharmonic functions are $\mathcal{E}^\theta$-quasi-continuous and locally belong to the domain $\mathcal{F}^\theta$. Moreover, the map from $SH_m(\Omega)\cap L_{loc}^\infty(\Omega)$ to $\mathcal{F}^\theta$ is continuous with respect to monotone convergence.
\begin{theorem}\label{Theorem: domain of dirichlet form}
    Let $\theta \in \mathcal{M}_m(\Omega)$, and $(\mathcal{F}^{\theta}, \mathcal{E}^{\theta})$ be the Dirichlet space for an admissible pair $(\theta, \mu)$. Let $u \in SH_m(\Omega) \cap L^{\infty}_{loc}(\Omega)$. Then, $u$ is $\mathcal{E}^\theta$-quasi-continuous and $\eta u \in \mathcal{F}^{\theta}$ for every $\eta \in \mathscr{C}_c^{\infty}(\Omega, \R)$. Moreover,  if $(u_k)_{k\in\N}\subset SH_m(\Omega)$ is monotone, i.e., increasing or decreasing, locally uniformly bounded, and converges to $u$ $\mu$-a.e. on $\Omega$, then $\eta u_k$ converges to $\eta u$ in $\mathcal{E}_1^{\theta}$.
\end{theorem}

Our second objective, carried out in Section $4$, is to use the analytic results obtained in the first objective to give new proofs of the fundamental results for complex Hessian operators and to further investigate the relationship between Dirichlet forms and $\mu$-quasi-monotone convergence for $m$-subharmonic functions.

The commutativity of complex Hessian operators and their continuity under monotone convergence are classical results in pluripotential theory. The continuity was proved in \cite{Blocki05,Lu12,Nguyen13} by adapting the argument of Bedford and Taylor \cite{BT82} for complex Monge--Amp\`ere operators, while commutativity follows as a consequence. Adapting the strategy of \cite[Section 3]{FO87} for complex Monge--Amp\`ere operators and using the results obtained in our first objective, we provide new proofs of these results for complex Hessian operators that are independent of the original ones. We remark that the proof of the continuity of the Monge--Amp\`ere operators under monotone convergence in \cite[Theorem 3.2]{FO87} is more quantitative than the original proof of Bedford and Taylor, in that it provides direct estimates for each gradient term. Our proof for complex Hessian operators enjoys the same feature. In the course of establishing these results, we also obtain the following monotonicity property for Dirichlet forms, inspired by the arguments in the proofs of \cite[Theorems 2.3 and 3.2]{FO87}.
\begin{theorem}\label{Theorem: monotonicity, introduction}
    Let $(\mathcal{F}^\theta,\mathcal{E}^\theta)$ and $(\mathcal{F}^{\Theta},\mathcal{E}^{\Theta})$ be the Dirichlet spaces for admissible pairs $(\theta,\mu)$ and $(\Theta,\nu)$, respectively, where
    \begin{align*}
        \theta=dd^c w_1\wedge\dots\wedge dd^c w_{m-1}\wedge\omega^{n-m},\quad \Theta=dd^c h_1\wedge\dots\wedge dd^c h_{m-1}\wedge\omega^{n-m},
    \end{align*}
    with $w_1,\dots,w_{m-1}, h_1,\dots,h_{m-1}\in SH_m(\Omega)\cap L_{loc}^{\infty}(\Omega)$. Let $u_1,u_2,v_1,v_2\in SH_m(\Omega)\\\cap L^{\infty}_{loc}(\Omega)$ and $K\subset\Omega$ be a compact set. Assume that
    \begin{align*}
        &w_k\leq h_k~ \text{on} ~K, \quad w_k=h_k ~\text{on} ~\Omega\setminus K,\quad \text{for}~k=1,\dots,m-1,\\
        &u_k\leq v_k~ \text{on} ~K, \quad u_k=v_k ~\text{on} ~\Omega\setminus K,\quad \text{for}~k=1,2.
    \end{align*}
    Then for every $\xi,\gamma \in \mathscr{C}_c^{\infty}(\Omega,\R)$ such that $\xi = \gamma=1$ on an open neighborhood of $K$, we have
    \begin{align*}
        \mathcal{E}^{\theta}(\xi u_1, \gamma u_2) \geq \mathcal{E}^{\Theta}(\xi v_1, \gamma v_2).
    \end{align*}
\end{theorem}

The notion of quasi-monotone convergence for $\omega_X$-plurisubharmonic functions on a compact K\"ahler manifold $(X,\omega_X)$ was introduced and studied in \cite{GT23}. Motivated by their monotone sandwiching formulation, we introduce a local notion of $\mu$-quasi-monotone convergence for $m$-subharmonic functions with respect to an arbitrary positive Radon measure $\mu$ on $\Omega$. When $\mu=V$, it includes monotone convergence and, in the plurisubharmonic setting, gives a localized version of the notion in \cite{GT23}. This notion provides a tool for further investigation of the relationship between Dirichlet forms and $m$-subharmonic functions. For any admissible pair $(\theta,\mu)$, we prove that the map from $SH_m(\Omega)\cap L_{loc}^\infty(\Omega)$ to $\mathcal{F}^\theta$ is continuous with respect to $\mu$-quasi-monotone convergence, thus further generalizing Theorem \ref{Theorem: domain of dirichlet form}.
\begin{theorem}\label{Theorem: E theta convergence and quasi-monotone convergence}
    Let $\theta\in\mathcal{M}_m(\Omega)$, and $(\mathcal{F}^\theta,\mathcal{E}^\theta)$ be the Dirichlet space for an admissible pair $(\theta,\mu)$. Assume that $(u_k)_{k\in\N}\subset SH_m(\Omega)$ is locally uniformly bounded and $u_k$ converges $\mu$-quasi-monotonically to $u\in SH_m(\Omega)\cap L_{loc}^\infty(\Omega)$. Then $(\eta u_k)_{k\in\N}$ converges to $\eta u$ in $\mc E^\theta_1$ for every $\eta\in \mathscr{C}_c^\infty(\Omega,\R)$.
\end{theorem}
The above theorem is new even in the setting of plurisubharmonic functions. As an application, we obtain the continuity of complex Hessian operators under $V$-quasi-monotone convergence.
\begin{theorem}\label{thm:quasi-monotone}
    Assume that $u_{j,k},u_j\in SH_m(\Omega)$ for $j=0,\dots,m$ and $k\in\N$ are locally uniformly bounded and satisfy the following conditions.
    \begin{itemize}
        \item[(i)] Either $\lim_{k\to\infty}u_{0,k}=u_0$ on $\Omega$, or $u_{0,k}$ converges $V$-quasi-monotonically to $u_0$  on $\Omega$.
        \item[(ii)] $u_{j,k}$ converges $V$-quasi-monotonically to $u_j$  on $\Omega$ for every $j=1,\dots,m$.
    \end{itemize}
    Then
    \begin{align*}
        u_{0,k} dd^c u_{1,k}\wedge\dots\wedge dd^c u_{m,k}\wedge\omega^{n-m}\to u_0 dd^c u_1\wedge\dots\wedge dd^c u_m\wedge\omega^{n-m},\quad k\to\infty,
    \end{align*}
    in the sense of currents.
\end{theorem}
We emphasize that, on the one hand, apart from Remarks \ref{Remark: Etheat convergence and BT capacity} and \ref{Remark:convergence of complex Hessian under BT capacity convergence}, the proofs of the results in our first two objectives do not rely on the Bedford--Taylor capacity $\mathrm{Cap}_{BT,m}$ for $m$-subharmonic functions. On the other hand, the arguments used in the proofs of Theorems \ref{Theorem: E theta convergence and quasi-monotone convergence} and \ref{thm:quasi-monotone} can also be adapted to extend these results to the setting in which the $m$-subharmonic functions converge in $\mathrm{Cap}_{BT,m}$, see the two aforementioned remarks for discussions.

Our third objective, carried out in Section 5, is to establish new characterizations of $m$-polar sets in terms of Dirichlet forms and stochastic processes through capacity comparisons and a probabilistic representation of the Choquet capacity $\mathrm{C}_{\sharp,m}$. 

We first revisit the $\mathcal{E}^\theta$-quasi-continuity of locally bounded $m$-subharmonic functions established in Theorem \ref{Theorem: domain of dirichlet form}, extending this property to unbounded functions and providing equivalent characterizations. This result is new even in the setting of plurisubharmonic functions.
\begin{theorem}
    Let $\theta\in \mathcal{M}_m(\Omega)$, $(\mathcal{F}^{\theta}, \mathcal{E}^{\theta})$ be the Dirichlet space for an admissible pair $(\theta, \mu)$, and $u \in SH_m(\Omega)$. Then the following statements are equivalent.
    \begin{itemize}
        \item [(i)] $u$ is $\mathcal{E}^{\theta}$-q.c.
        \item [(ii)] $\mu\left(\{u = -\infty\}\right) = 0$.
        \item [(iii)] $\mathrm{Cap}_{\theta}(\{u = -\infty\}) = 0$.
    \end{itemize}
\end{theorem}

As a supplementary observation on $V$-quasi-monotone convergence, using an argument similar to that in \cite{GT23}, we show that for any sequence of locally uniformly bounded $m$-subharmonic functions, convergence in $\mathrm{Cap}_{BT, m}$ and $V$-quasi-monotone convergence are equivalent after passing to a subsequence. 
Next, we establish the quantitative comparison between $\mathrm{Cap}_{BT,m}$ and $\mathrm{Cap}_{\theta}$ for $\theta\in\mathcal{M}^+_m(\Omega)$, where $\mathrm{Cap}_{\theta}$ is the capacity induced by the admissible pair $(\theta,V)$. Consequently, on a bounded strictly $m$-pseudoconvex domain $\Omega$, a Borel subset $E$ satisfies $\mathrm{Cap}_{BT, m}(E,\Omega) = 0$ if and only if $\mathrm{Cap}_{\theta}(E)=0$ for every $\theta\in\mathcal{M}_m^+(\Omega)$. This generalizes \cite[Lemma 7.8 and Corollary 7.9]{Nguyen26}. Then, we incorporate a stochastic perspective into the characterization of $m$-polar sets. We show that the Choquet capacity $\mathrm{C}_{\sharp,m}$ introduced in this paper can be computed in terms of the conformal diffusion associated with the Dirichlet form $\mathcal{E}^\theta$, generalizing \cite[Theorem 6.5]{FO87}.
\begin{theorem}\label{Theorem : stochastic part}
    Let $\Omega \subset \C^n$ be a bounded strictly $m$-pseudoconvex domain, and $E \subset \Omega$ be a Borel subset.  For every $\theta\in \mathcal{M}_m^+(\Omega)$, let $M^{\theta}$ be a conformal diffusion associated with the Dirichlet form $\mathcal{E}^{\theta}$. Then
    \begin{align*}
        \mathrm{C}_{\sharp, m}(E) = \sup_{\theta \in \mathcal{M}_m^+(\Omega)}\int_{\Omega}P_z^{\theta}(\sigma_E < \infty)dV(z).
    \end{align*}
\end{theorem}
Combining these results, we obtain the following characterizations of $m$-polar sets, substantially generalizing \cite[Theorem 6.6]{FO87}.
\begin{theorem}\label{Theorem: main result}
    Let $\Omega \subset \C^n$ be a bounded strictly $m$-pseudoconvex domain, and $E \subset \Omega$ be a Borel subset. For every $\theta\in \mathcal{M}_m^+(\Omega)$, let $M^{\theta}$ be a conformal diffusion associated with the Dirichlet form $\mathcal{E}^{\theta}$. Then the following statements are equivalent.
    \begin{itemize}
        \item [(i)] $u_{E, \Omega, m}^* \equiv 0$.
        \item [(ii)] $\mathrm{C}_{\sharp, m}(E) = 0$.
        \item [(iii)] $\mathrm{Cap}_{BT,m}(E, \Omega) = 0$.
        \item [(iv)] $\mathrm{Cap}_\theta(E) = 0$ for every $\theta \in \mathcal{M}_m^+(\Omega)$.
        \item [(v)] $P_z^{\theta}(\sigma_{E} < \infty) = 0$ $V$-a.e. on $\Omega$ for every $\theta \in \mathcal{M}_m^+(\Omega)$.       
    \end{itemize}
\end{theorem}

\subsection*{Organization} 
In Section $2$, we introduce preliminaries on Dirichlet forms, $m$-subharmonic functions, and $\mu$-quasi-monotone convergence. In Section $3$, we study the basic relationship between Dirichlet spaces and $m$-subharmonic functions. In Section $4$, we apply the Dirichlet form method to complex Hessian operators, establish the monotonicity property for Dirichlet forms, and prove continuity results under $\mu$-quasi-monotone convergence. In Section $5$, we characterize $m$-polar sets.
\subsection*{Notations}
In this paper, $\Omega$ denotes a bounded open subset of $\C^n$, $\theta$ denotes a closed positive $(n-1,n-1)$-current on $\Omega$, and $\mu$ denotes a positive Radon measure on $\Omega$ satisfying $\supp \mu=\Omega$, e.g., the Lebesgue measure $V$ on $\C^n$. We denote the space of real-valued smooth functions on $\Omega$ by $\mathscr{C}^{\infty}(\Omega,\R)$, and the space of real-valued smooth functions on $\Omega$ with compact supports by $\mathscr{C}_c^{\infty}(\Omega,\R)$. For every open subset $D$ of $\Omega$, we identify $\mathscr{C}_c^{\infty}(D,\R)$ with the subspace of $\mathscr{C}_c^{\infty}(\Omega,\R)$ consisting of functions supported in $D$. For $1\leq p \leq \infty$, we denote the space of extended real-valued $L^p$ functions with respect to $\mu$ by $L^p(\Omega;\mu)$, and write $L^p(\Omega;V)$ simply as $L^p(\Omega)$. Set $d^c:=\frac{i}{2\pi}(\bar{\partial}-\partial)$, so that $dd^c=\frac{i}{\pi}\partial\overline{\partial}$. Set $\omega:=dd^c |z|^2$.

\subsection*{Acknowledgement}
We would like to thank Ngoc Cuong Nguyen, our advisor and mentor, for suggesting this problem and for his invaluable guidance. We are grateful to Jun-Muk Hwang for mentoring the third author and for providing the first and second authors with the opportunity to visit the IBS Center for Complex Geometry during this project. We would also like to thank Hyunsoo Ahn and Shuang Su for helpful discussions.

Bowoo Kang and Jaechul Lee were supported by the National Research Foundation of Korea (NRF, Grant No. RS-2026-25470686). Bowoo Kang was also supported by the DFG-Project VU 126/1-2. Zhenghao Li was supported by the Institute for Basic Science (IBS-R032-D1-2026-a00).

\section{Preliminaries}

In this section, we introduce the preliminaries concerning Dirichlet forms and the induced capacity $\mathrm{Cap}_\theta$, $m$-subharmonic functions and the capacity $\mathrm{C}_{\sharp, m}$, as well as $\mu$-quasi-monotone convergence of $m$-subharmonic functions.

\subsection{Dirichlet forms and the capacity $\mathrm{Cap}_\theta$}
In this subsection, we briefly recall Dirichlet forms, Dirichlet spaces, and the capacity $\mathrm{Cap}_\theta$. We refer the reader to \cite[Chapters 1 and 2]{FOT11} for further details.

We consider the following symmetric bilinear form on $\mathscr{C}_c^{\infty}(\Omega,\R)$, as in (\ref{equation:Etheta}).
\begin{align*}
    \mathcal{E}^\theta: \mathscr{C}_c^{\infty}(\Omega,\R) \times \mathscr{C}_c^{\infty}(\Omega,\R)&\to \R,\\
    (\varphi, \psi)&\mapsto \int_\Omega d\varphi\wedge d^c \psi\wedge\theta.
\end{align*}
This map is non-negative definite, i.e., $\mathcal{E}^\theta(\varphi,\varphi)\geq0$ for every $\varphi\in \mathscr{C}_c^{\infty}(\Omega,\R)$, and hence it is a semi-inner product on $\mathscr{C}_c^\infty(\Omega,\R)$. 
One computes via integration by parts that
\begin{align}\label{Equation:another formula of Dirichlet form}
    \mathcal{E}^{\theta}(\varphi, \psi) = -\int_{\Omega}\varphi dd^c\psi \wedge \theta=\frac{i}{\pi}\int_\Omega \partial\varphi\wedge\bar{\partial}\psi\wedge\theta.
\end{align}
Since $\supp \mu=\Omega$, every element of the Hilbert space $L^2(\Omega;\mu)$ has at most one continuous representative. Thus, $\mathcal{E}^\theta$ induces a \textit{symmetric form} on the real Hilbert space $L^2(\Omega;\mu)$, namely, a semi-inner product on a dense linear subspace of $L^2(\Omega;\mu)$. We also denote by $\mathcal{E}^\theta$ this symmetric form.
\begin{lemma}\label{Lemma:Etheta Markovian}
    The symmetric form $\mathcal{E}^\theta$ is Markovian.    
\end{lemma}
\begin{proof}
    For every $\varepsilon>0$, take $f_\varepsilon\in \mathscr{C}^\infty(\R,\R)$ that satisfies the following properties.
    \begin{itemize}
        \item[(i)] $-\varepsilon\leq f_\varepsilon\leq 1+\varepsilon$.
        \item[(ii)] $f_\varepsilon(t)=t$ for $0\leq t\leq 1$.
        \item[(iii)] $0\leq f_\varepsilon(t_2)-f_\varepsilon(t_1)\leq t_2 -t_1$ whenever $t_1<t_2$.
    \end{itemize}
    Notice that $0\leq f'_\varepsilon\leq 1$. For every $\varphi\in \mathscr{C}_c^\infty(\Omega,\R)$, we have $f_\varepsilon(\varphi)\in \mathscr{C}_c^\infty(\Omega,\R)$ and 
    \[
        \mathcal{E}^\theta(f_\varepsilon(\varphi),f_\varepsilon(\varphi))=\frac{i}{\pi}\int_\Omega (f'_\varepsilon(\varphi))^2\partial \varphi\wedge \overline{\partial}\varphi \wedge \theta\leq \frac{i}{\pi}\int_\Omega \partial \varphi\wedge \overline{\partial}\varphi \wedge \theta=\mathcal{E}^\theta(\varphi,\varphi).
    \]
    Therefore, $\mathcal{E}^\theta$ is Markovian.
\end{proof}

The notions of \textit{closability} and \textit{closedness} are fundamental to symmetric forms. It is well known that a symmetric form is \textit{closable} if and only if it admits a \textit{closed} extension. A Markovian closed symmetric form is called a \textit{Dirichlet form}. Following \cite{FO87}, we say that $(\theta,\mu)$ is an \textit{admissible pair} if $\mathcal{E}^\theta$ is closable, that is, for $(\varphi_k)_{k\in\N}\subset \mathscr{C}_c^\infty(\Omega,\R)$,
\begin{align*}
    \lim_{k,\ell\to\infty}\mathcal{E}^\theta(\varphi_k-\varphi_{\ell},\varphi_k-\varphi_{\ell})=0,~\lim_{k\to\infty}\int_\Omega \varphi_k^2 ~d\mu=0~\implies~ \lim_{k\to\infty}\mathcal{E}^\theta(\varphi_k,\varphi_k)=0.
\end{align*}
The following example of an admissible pair was provided by \cite[Lemma 1]{FO84}, see also \cite[Lemma 7.2]{Nguyen26} for an alternative proof.
\begin{align}\label{Example:admissible pair}
    (\theta,\theta\wedge \omega) \ \text{is an admissible pair whenever} \ \supp{\theta\wedge \omega}=\Omega.
\end{align}
Let $D$ be an open subset of $\Omega$. One simple observation is
\begin{align}\label{admissible pair and open subset}
    (\theta,\mu) ~\text{is an admissible pair}~\implies~ (\theta|_D,\mu|_D) ~\text{is an admissible pair}.
\end{align}

Assume that $(\theta,\mu)$ is an admissible pair. \cite[Theorem 3.1.1]{FOT11} together with Lemma \ref{Lemma:Etheta Markovian} implies that the smallest closed extension of $\mathcal{E}^\theta$ remains Markovian and is therefore a Dirichlet form. We continue to denote this extension by $\mathcal{E}^\theta$, denote its domain by $\mathcal{F}^\theta$, and refer to $(\mathcal{F}^\theta, \mathcal{E}^\theta)$ as the Dirichlet space for the admissible pair $(\theta,\mu)$. Notice that $(\mathcal{F^\theta}, \mathcal{E}^\theta_1)$ is a Hilbert space, where
$$\mathcal{E}^\theta_1 (u,v):= \mathcal{E}^\theta(u,v)+\int_\Omega uv \ d\mu, \quad u,v\in \mathcal{F^\theta}.$$
By \cite[Theorem 3.1.2]{FOT11}, $\mathcal{E}^\theta$ is a \textit{regular} Dirichlet form with the \textit{local property}, and $\mathscr{C}_c^\infty(\Omega,\R)$ is a \textit{special standard core} of $\mathcal{E}^\theta$. In particular, $\mathscr{C}_c^\infty(\Omega,\R)$ is dense in $(\mathcal{F^\theta}, \mathcal{E}^\theta_1)$. We may then consider the capacity $\mathrm{Cap}_\theta$ on $\Omega$, which is induced by the Dirichlet form $\mathcal{E}^\theta$. 
\begin{definition}\cite[Section 2.1]{FOT11}\label{definition:CapTheta}
    For every open subset $D \subset \Omega$, we define
    \begin{align*}
        \mathcal{L}_{D} := \{u \in \mathcal{F}^\theta: u \geq 1 \ \mu\text{-a.e. on} \ D\},
    \end{align*}
    and set
    \begin{align*}
        \mathrm{Cap}_\theta(D) := \begin{cases}
            \inf_{u \in \mathcal{L}_D}\mathcal{E}^\theta_1(u, u),&\mathcal{L}_D \neq\emptyset,\\
            \infty, &\mathcal{L}_D = \emptyset.
        \end{cases}
    \end{align*}
    For every subset $E\subset\Omega$, we define 
    \begin{align*}
        \mathrm{Cap}_\theta(E) := \inf\left\{\mathrm{Cap}_\theta(D): E \subset D\subset\Omega, D~\text{open}\right\}.
    \end{align*}
\end{definition}

It follows from the definition that $\mu(E)\leq\mathrm{Cap}_\theta(E)$ for every Borel subset $E\subset\Omega$.
\begin{theorem}\cite[Theorem 2.1.1]{FOT11}
    The map $\mathrm{Cap}_\theta: 2^\Omega\to [0,\infty]$ is a Choquet capacity on $\Omega$, namely, it has the following properties.
\begin{itemize}
    \item [(i)] $\mathrm{Cap}_\theta(\emptyset)=0$.
    \item [(ii)] If $E_1 \subset E_2$, then $\mathrm{Cap}_\theta(E_1) \leq \mathrm{Cap}_\theta(E_2)$.
    \item [(iii)] If $E_k \uparrow E$, then $\mathrm{Cap}_\theta(E_k)\uparrow\mathrm{Cap}_\theta(E)$.
    \item [(iv)] If $G_k$ are compact and $G_k \downarrow G$, then $\mathrm{Cap}_\theta(G_k)\downarrow\mathrm{Cap}_\theta(G)$.
\end{itemize}
\end{theorem}
By the Choquet capacitability theorem \cite{Choquet}, for every Borel subset $E\subset\Omega$,
\begin{align*}
    \mathrm{Cap}_\theta(E) = \sup \left\{\mathrm{Cap}_\theta(G): G \subset E, ~G~\text{compact}\right\}.
\end{align*}
It follows from \cite[Lemma 2.1.2]{FOT11} that $\mathrm{Cap}_\theta$ is countably subadditive, i.e.,
\begin{align*}
    \mathrm{Cap}_\theta\left(\bigcup_{k=1}^\infty E_k\right)\leq \sum_{k=1}^\infty \mathrm{Cap}_\theta(E_k).
\end{align*}

Let $E\subset\Omega$. We say that $E$ is \textit{$\mc E^\theta$-polar} if $\mathrm{Cap}_\theta (E)=0$.  A statement concerning $z \in E$ is said to hold \textit{$\mathcal{E}^\theta$-quasi-everywhere on $E$} ($\mathcal{E}^\theta$-q.e. on $E$ for short) if there exists an $\mc E^\theta$-polar set $N \subset E$ such that the statement holds for every $z \in E \setminus N$.
Let $u$ be an extended real-valued function defined $\mathcal{E}^\theta$-q.e. on $\Omega$. We say that $u$ is \textit{$\mathcal{E}^\theta$-quasi-continuous} ($\mathcal{E}^\theta$-q.c. for short) if, for every $\varepsilon > 0$, there exists an open set $D\subset\Omega$ such that $\mathrm{Cap}_\theta(D) < \varepsilon$ and $u$ is finite-valued and continuous on $\Omega \setminus D$.

\subsection{$m$-subharmonic functions and the capacity $\mathrm{C}_{\sharp,m}$}
In this subsection, we recall several notions concerning $m$-subharmonic functions and introduce the capacity $\mathrm{C}_{\sharp,m}$. We refer the reader to \cite{Blocki05, Lu12, Nguyen13, DK14} for a comprehensive introduction to the theory of $m$-subharmonic functions.

Fix a standard mollifier $\chi\in\mathscr{C}_c^\infty(\C^n,\R)$ and denote $\chi_\delta(z):=\frac{1}{\delta^{2n}}\chi(\frac{z}{\delta})$ for $\delta>0$. For $u\in SH_m(\Omega)$, we denote by $u^{[\delta]}:=u * \chi_{\delta}$ its standard regularization. Then $u^{[\delta]}\in SH_m(\Omega_\delta)\cap \mathscr{C}^\infty(\Omega_\delta,\R)$ decreases to $u$ as $\delta\downarrow 0$.

Let $1 \leq k \leq m$. The complex Hessian operator
\begin{align*}
    dd^cu_k \wedge dd^cu_{k-1} \wedge \cdots \wedge dd^cu_1 \wedge \omega^{n-m}
\end{align*}
is symmetric in $u_1, \dots, u_k \in SH_m(\Omega) \cap L^{\infty}_{loc}(\Omega)$ and continuous under monotone convergence. We emphasize that, prior to Section $5$, neither of these two properties is assumed a priori, both are established in Section $4$ using new proofs.

\begin{definition}\label{def:m-relative extremal function}
    We say that a subset $E$ of $\C^n$ is \textit{$m$-polar} if for every $z\in E$ there exists an open neighborhood $D$ of $z$ and a function $u\in SH_m(D)$ such that $D\cap E\subset\{u=-\infty\}$.
\end{definition}

\begin{definition}
    We say that a subset $E$ of $\Omega$ is \textit{$m$-negligible} if there exists a family $(u_\alpha)_{\alpha\in I}\subset SH_m(\Omega)$ such that $(u_\alpha)_{\alpha\in I}$ is locally uniformly bounded from above and $E\subset \{u<u^*\}$. Here $u:=\sup_{\alpha\in I}u_\alpha$ and $u^*$ denotes the upper semi-continuous regularization of $u$.
\end{definition}

\begin{definition}
    For every subset $E\subset\Omega$, the \textit{$m$-relative extremal function} of $E$ with respect to $\Omega$ is defined by 
    \[u_{E,\Omega,m}:=\sup\{u\in SH_m(\Omega):u\leq0, u|_E\leq -1\}.\]
\end{definition}

Notice that $u_{E,\Omega,m}$ remains unchanged if, in its definition, $SH_m(\Omega)$ is replaced by the class of functions that are either $m$-subharmonic or identically $-\infty$ on each connected component of $\Omega$. We now introduce the capacity $\mathrm{C}_{\sharp,m}$, a natural generalization of the capacity $\mathrm{C}_{\sharp}$ associated with plurisubharmonic functions. The latter was studied in \cite{Cegrell80,BT82,Fu86a,Fu86b, FO87}.

\begin{definition}
    For every subset $E\subset\Omega$, we define
    \begin{align*}
        \mathrm{C}_{\sharp, m}(E) := -\int_{\Omega}u_{E, \Omega, m}^*(z)dV(z). 
    \end{align*}
\end{definition}

It is easy to check that $\mathrm{C}_{\sharp,m}(\emptyset)=0$ and $\mathrm{C}_{\sharp,m}(E_1) \leq \mathrm{C}_{\sharp,m}(E_2)$ if $E_1 \subset E_2$. We now prove that $\mathrm{C}_{\sharp,m}$ is a Choquet capacity by adapting the arguments of \cite[Lemma 6.1 and Theorem 6.4]{FO87}. We emphasize that the capacity $\mathrm{C}_{\sharp,m}$ will not be used until Section $5$.
\begin{proposition}\label{proposition: Csharpm is Choquet}
    $\mathrm{C}_{\sharp,m}$ satisfies the following properties.
    \begin{itemize}
        \item[(i)] For every subset $E\subset\Omega$,
        \begin{align*}
            \mathrm{C}_{\sharp,m}(E)=\inf \{\mathrm{C}_{\sharp,m}(D): E\subset D\subset\Omega, D ~\text{open}\}.
        \end{align*}
        \item[(ii)] $\mathrm{C}_{\sharp,m}$ is countably subadditive.
        \item[(iii)] If $E$ is $m$-polar, then $\mathrm{C}_{\sharp,m}(E)=0$.
        \item[(iv)] If $E_k \uparrow E$, then $\mathrm{C}_{\sharp,m}(E_k)\uparrow\mathrm{C}_{\sharp,m}(E)$.
    \end{itemize}
    Consequently, by \cite[Proposition 4.3]{GZ17}, $\mathrm{C}_{\sharp, m}$ is a Choquet capacity. In particular, the Choquet capacitability theorem \cite{Choquet} reads that for every Borel subset $E\subset\Omega$
\begin{align*}
    \mathrm{C}_{\sharp,m}(E) = \sup \left\{\mathrm{C}_{\sharp,m}(G): G \subset E, ~G~\text{compact}\right\}.
\end{align*}
\end{proposition}
\begin{proof}
    (i) By Choquet's lemma, there exists a sequence $(u_k)_{k\in\N}\subset SH_m(\Omega)$, such that $(u_k)_{k\in\N}$ is increasing, $u_k\leq 0$, $u_k|_E\leq -1$, and $\left(\sup_{k\in\N}u_k\right)^*=u_{E,\Omega,m}^*$. Hence, $u_k\leq u_{E,\Omega,m}$ and $\lim_{k\to\infty}u_k=u_{E,\Omega,m}~V$-a.e. on $\Omega$. Set
    \[
    D_k:=\left\{z\in\Omega: \left(1+\frac{1}{k}\right)u_k(z)<-1\right\}
    \]
    for every $k\in\N$, which is an open subset of $\Omega$ and contains $E$. We then have
    \[
    \left(1+\frac{1}{k}\right)u_k\leq u_{D_k,\Omega,m}\leq u_{E,\Omega,m},
    \]
    which implies that 
    \[
    \lim_{k\to\infty}u_{D_k,\Omega,m}=u_{E,\Omega,m},\quad V\text{-a.e. on}~ \Omega.
    \]
    Therefore, $u^*_{D_k,\Omega,m}\to u^*_{E,\Omega,m}~V$-a.e. on $\Omega$, and hence $ \mathrm{C}_{\sharp,m}(D_k)\to\mathrm{C}_{\sharp,m}(E)$ follows.
    
    (ii) Let $D_k\subset\Omega$ be an open subset for every $k\in\N$. Set $D:=\cup_{k=1}^\infty D_k$. Then, $u^*_{D_k,\Omega,m}=-1$ on $D_k$,
    $\sum_{k=1}^\infty u^*_{D_k,\Omega,m}\leq -1$ on $D$, and $\sum_{k=1}^\infty u^*_{D_k,\Omega,m}$ is either $m$-subharmonic or identically $-\infty$ on each connected component of $\Omega$. Therefore,
    \[
    \sum_{k=1}^\infty u^*_{D_k,\Omega,m}\leq u_{D,\Omega,m}\leq u^*_{D,\Omega,m},
    \]
    and hence $\mathrm{C}_{\sharp,m}\left(D\right)\leq \sum_{k=1}^\infty \mathrm{C}_{\sharp,m}(D_k)$ follows.
    
    Combining the preceding result with (i), we obtain the desired assertion.
    
    (iii) By \cite[Theorem 1.6.5]{Lu12} and the fact that $\Omega$ is bounded, there exists a negative function $u\in SH_m(\Omega)$ such that $E\subset u^{-1}(-\infty)$. Since $\varepsilon u\leq u_{E,\Omega,m}$ for every $\varepsilon>0$, we deduce that $u_{E,\Omega,m}=0$ outside $u^{-1}(-\infty)$. Therefore, $u^*_{E,\Omega,m}\equiv0$, and hence $\mathrm{C}_{\sharp,m}(E)=0$ follows.

    (iv) Let $v\in SH_m(\Omega)$ be the limit of the decreasing sequence $(u^*_{E_k,\Omega,m})_{k\in\N}$. Set $E':=\{z\in E: v(z)=-1\}$. On the one hand, $v\leq u_{E',\Omega,m}$ implies that
    \[
    \mathrm{C}_{\sharp,m}(E')\leq -\int_\Omega v ~dV = -\lim_{k\to\infty}\int_\Omega u^*_{E _k,\Omega,m}~ dV=\lim_{k\to\infty} \mathrm{C}_{\sharp,m}(E_k).
    \]
    On the other hand,
    \[
    E\setminus E'\subset \bigcup_{k=1}^\infty \left\{z\in E_k: u^*_{E_k,\Omega,m}(z)>-1\right\}\subset \bigcup_{k=1}^\infty \left\{u^*_{E_k,\Omega,m}>u_{E_k,\Omega,m}\right\}.
    \]
    Using the fact that $m$-negligible sets are $m$-polar \cite[Theorem 1.6.7]{Lu12}, together with (ii, iii),
    we conclude that $\mathrm{C}_{\sharp,m}(E\setminus E')=0$.
    
    Therefore, $\mathrm{C}_{\sharp,m}(E)\leq \lim_{k\to\infty} \mathrm{C}_{\sharp,m}(E_k)$. The converse inequality is clear.
\end{proof}

We conclude this subsection by recalling the definition of bounded strictly $m$-pseudoconvex domains and the localization principle for $m$-subharmonic functions.

\begin{definition}\label{Definition: bounded strictly m-psc domain}
    Let $\Omega$ be a bounded domain in $\C^n$ with $\mathscr{C}^2$ boundary. We say that $\Omega$ is a \textit{bounded strictly $m$-pseudoconvex domain} if there exist a $\mathscr{C}^2$ defining function $\rho$ of $\Omega$ and a constant $C>0$ such that $(dd^c\rho)^k\wedge\omega^{n-k}\geq C \omega^n$ on $\overline{\Omega}$ for $1\leq k\leq m$.
\end{definition}

\begin{proposition}[Localization principle]\label{Proposition:localization principle}
    Let $\Omega\subset\C^n$ be a bounded strictly $m$-pseudoconvex domain and $\rho$ be a function satisfying the property in Definition \ref{Definition: bounded strictly m-psc domain}. Then for every compact subset $G$ of $\Omega$ and constants $B_2>B_1>0$, there exist an open set $U$, a compact set $K$, and a constant $A>0$ that satisfy $G\subset U\subset K\subset \Omega$ and the following property.
    
    For every $u\in SH_m(\Omega)$ with $-B_2\leq u\leq -B_1$, the function $\hat{u}:=\max(u, A\rho)$ is $m$-subharmonic on $\Omega$, equal to $u$ on $U$ and equal to $A\rho$ on $\Omega\setminus K$.
\end{proposition}
\begin{proof}
    Take $\delta_1>0$ so that $\rho<-\delta_1$ on $G$. Take $\delta_2\in(0,\delta_1)$ so that $\frac{B_2}{\delta_1}<\frac{B_1}{\delta_2}$. Then $U:=\{\rho<-\delta_1\}$, $K:=\{\rho\leq -\delta_2\}$, and any $A\in(\frac{B_2}{\delta_1},\frac{B_1}{\delta_2})$ meet the requirement.
\end{proof}

\subsection{$\mu$-quasi-monotone convergence of $m$-subharmonic functions}
Guedj and Trusiani \cite[Definition 2.1]{GT23} introduced the notion of quasi-monotone convergence for $\omega_X$-plurisubharmonic functions on a compact K\"ahler manifold $(X,\omega_X)$. In this subsection, we adapt their monotone sandwiching formulation to define a local notion of $\mu$-quasi-monotone convergence for $m$-subharmonic functions on $\Omega$, where $\mu$ is an arbitrary positive Radon measure. Here we don't require $\supp\mu=\Omega$.
\begin{definition}\label{Definition of convergence quasi-monotonically}
    Let  $(u_k)_{k \in \N} \subset SH_m(\Omega)\cap L_{loc}^1(\Omega;\mu)$ and $u \in SH_m(\Omega)\cap L_{loc}^1(\Omega;\mu)$. We say that $u_k$ converges to $u$ $\mu$-quasi-monotonically on $\Omega$ if, for every open subset $D\Subset\Omega$, there exist sequences $(\phi_k)_{k \in \N}, (\psi_k)_{k \in \N} \subset SH_m(D)\cap L_{loc}^1(D;\mu)$ that satisfy the following conditions.
    \begin{itemize}
        \item [(i)] $\phi_k \leq u_k \leq \psi_k$ on $D$ for every $k\in\N$.
        \item [(ii)] $(\phi_k)_{k \in \N}$ is increasing and $(\psi_k)_{k \in \N}$ is decreasing.
        \item [(iii)] Both $\phi_k$ and $\psi_k$ converge to $u$ in $L^1_{loc}(D; \mu)$.
    \end{itemize}
    These conditions imply that $u_k$ converges to $ u$ in
    $L^1_{loc}(\Omega;\mu)$ and $\mu$-a.e. on $\Omega$, and that both $\phi_k$ and $\psi_k$ converge to $u$ $\mu$-a.e. on $D$.
\end{definition}
\begin{remark}
    Let $(u_k)_{k\in\N}\subset SH_m(\Omega)$ and $u\in SH_m(\Omega)$. Assume that either $(u_k)_{k\in\N}$ is decreasing with $\lim_{k\to\infty}u_k=u$ on $\Omega$, or $(u_k)_{k\in\N}$ is increasing with $\lim_{k\to\infty}u_k=u$ $V$-a.e. on $\Omega$. Then $u_k$ converges to $u$ $V$-quasi-monotonically.
\end{remark}

\begin{remark}\label{Remark: mu-quasi-monotone convergence for locally bounded functions}
    Assume that $(u_k)_{k \in \N} \subset SH_m(\Omega)$ converges to $u \in SH_m(\Omega)$ $\mu$-quasi-monotonically on $\Omega$, and $(u_k)_{k\in\N}$ is  locally uniformly bounded on $\Omega$. For every open subset $U\Subset\Omega$, choose an open set $D$ such that $U\Subset D\Subset\Omega$, and let $(\phi_k)_{k\in\N}$ and $(\psi_k)_{k\in\N}$ be the corresponding sequences in Definition~\ref{Definition of convergence quasi-monotonically}. Take a constant $C$ so that $u_k\geq C$ on $U$ for every $k\in\N$, and set $$\varphi_k:=\max(\phi_k,C).$$
    Then $(\varphi_k)_{k\in\N},(\psi_k)_{k\in\N}\subset SH_m(U)$ are uniformly bounded and satisfy the following properties.
        \begin{itemize}
            \item [(i)] $\varphi_k \leq u_k \leq \psi_k$ on $U$ for every $k\in\N$.
            \item [(ii)] $(\varphi_k)_{k \in \N}$ is increasing and $(\psi_k)_{k \in \N}$ is decreasing.
            \item [(iii)] Both $\varphi_k$ and $\psi_k$ converge to $u$ in $L^1_{loc}(U; \mu)$.
            \item[(iv)] $\varphi_k$ converges to $u$ $\mu$-a.e. on $U$ and $\psi_k$ converges to $u$ $\mu$-a.e. on $U$.
        \end{itemize}
    When $\mu=V$, the decreasing sequence $(\psi_k)_{k\in\N}$ converges to $u$ at every point of $U$.
\end{remark}
For $\omega_X$-plurisubharmonic functions, examples of quasi-monotone convergence arise from families of solutions
to Monge--Amp\`ere equations and Monge--Amp\`ere flows,
see \cite[Examples 3.1 and 3.2]{GT23}.
In Section~5, after recalling the Bedford--Taylor capacity
for $m$-subharmonic functions, we will establish its relationship
with $V$-quasi-monotone convergence on an arbitrary bounded
open subset of $\C^n$, as suggested by the discussion in
\cite[Section 3.3.3]{GT23}. Indeed,
$V$-quasi-monotone convergence implies convergence in
Bedford--Taylor capacity, and the converse holds after passing
to a subsequence when the sequence is locally uniformly bounded.

\section{Relationship between $(\mathcal{F}^\theta,\mathcal{E}^\theta)$ and $m$-subharmonic functions}
In this section, we study the Dirichlet space $(\mathcal{F}^{\theta}, \mathcal{E}^{\theta})$ when $\theta\in\mathcal{M}_m(\Omega)$ and $(\theta,\mu)$ is an admissible pair. We show that the results in \cite[Section 2]{FO87} are well extended to $m$-subharmonic functions.

We shall utilize the Schwarz inequality in \cite[Lemma 1]{Okada82}.
\begin{lemma}[Schwarz inequality] For every closed positive $(n-1,n-1)$-current $\theta$ on $\Omega$ and $\xi, \eta, \varphi,\psi\in \mathscr{C}^\infty(\Omega,\R)$, we have  \begin{equation}\label{Equation:Schwarz inequality}
        \left(\int_\Omega \eta\xi d\varphi \wedge d^c\psi \wedge \theta\right)^2 \leq \left(\int_\Omega \eta^2 d\varphi \wedge d^c\varphi \wedge \theta\right)\left(\int_\Omega \xi^2 d\psi \wedge d^c\psi \wedge \theta\right)
    \end{equation}
provided that both $\eta d\varphi$ and $\xi d\psi$ have compact supports.
\end{lemma}

For every $\eta \in \mathscr{C}^{\infty}(\Omega,\R)$ and $\varphi \in \mathscr{C}_c^{\infty}(\Omega,\R)$, (\ref{Equation:Schwarz inequality}) implies that
\begin{equation}\label{Equation:corollary of Schwarz inequality}
    \int_\Omega d(\eta\varphi)\wedge d^c (\eta\varphi)\wedge\theta\leq 2\int_\Omega \eta^2 d\varphi \wedge d^c\varphi \wedge \theta +  2\int_\Omega \varphi^2 d\eta \wedge d^c\eta \wedge \theta.
\end{equation}

We establish the following Poincar\'e-type inequality for $m$-subharmonic functions, generalizing Okada's result \cite[Theorem 1]{Okada82} for $u\in PSH(\Omega)\cap L_{loc}^\infty(\Omega)$ and an arbitrary closed positive $(n-1,n-1)$-current $\theta$.

\begin{lemma}[Poincar\'e-type inequality]\label{Poincare type inequality} Assume that $\theta \in \mathcal{M}_m(\Omega), u \in SH_m(\Omega)\cap L^\infty_{loc}(\Omega)$, and $\varphi \in \mathscr{C}_c^{\infty}(\Omega,\R)$. Let $D\Subset\Omega$ be an open subset containing $\supp\varphi$. Then we have 
    \begin{align*}
        \int_\Omega\varphi^2dd^cu \wedge \theta \leq 8\lVert u\rVert_{L^{\infty}(D)}\mathcal{E}^{\theta}(\varphi, \varphi).
    \end{align*}
\end{lemma}
\begin{proof}
    The proof is adapted from that of \cite[Theorem 1]{Okada82}. We provide it here for the reader's convenience. Assume first that $u \in SH_m(\Omega)\cap \mathscr{C}^{\infty}(\Omega,\R)$. By the Schwarz inequality (\ref{Equation:Schwarz inequality}), we have
    \begin{align*}
        \int_\Omega \varphi^2 dd^c u \wedge \theta = -2\int_D \varphi d\varphi \wedge d^cu \wedge \theta 
         \leq 2\left(\int_D \varphi^2 du \wedge d^cu \wedge \theta\right)^{\frac{1}{2}}\left(\mathcal{E}^{\theta}(\varphi, \varphi)\right)^{\frac{1}{2}}.
    \end{align*}    
 By the Schwarz inequality (\ref{Equation:Schwarz inequality}) again and the above estimate, we have
    \begin{align*}
        &\int_D\varphi^2 du \wedge d^cu \wedge \theta \\
        &= -2\int_D \varphi u d\varphi \wedge d^cu \wedge \theta - \int_D \varphi^2 u dd^cu \wedge \theta \\
        &\leq  2\left(\int_D \varphi^2du \wedge d^cu \wedge \theta\right)^{\frac{1}{2}}\left(\int_D u^2 d\varphi \wedge d^c\varphi \wedge \theta\right)^{\frac{1}{2}}+\lVert u\rVert_{L^{\infty}(D)}\int_D \varphi^2 dd^cu \wedge \theta \\
        &\leq  2\lVert u\rVert_{L^{\infty}(D)}\left(\int_D \varphi^2 du \wedge d^cu \wedge \theta\right)^{\frac{1}{2}}\mathcal{E}^{\theta}(\varphi, \varphi)^{\frac{1}{2}} + \lVert u\rVert_{L^{\infty}(D)}\int_D \varphi^2dd^cu \wedge \theta\\
        &\leq  4\lVert u\rVert_{L^{\infty}(D)}\left(\int_D \varphi^2 du \wedge d^cu \wedge \theta\right)^{\frac{1}{2}}\mathcal{E}^{\theta}(\varphi, \varphi)^{\frac{1}{2}}.
    \end{align*}
     We may assume that $\int_D \varphi^2 du \wedge d^cu \wedge \theta>0$, since otherwise the left-hand side of the desired inequality is zero. Therefore, we have
    \begin{align*}
        \left(\int_D \varphi^2du \wedge d^cu \wedge \theta \right)^{\frac{1}{2}}\leq 4\lVert u\rVert_{L^{\infty}(D)} \mathcal{E}^{\theta}(\varphi, \varphi)^{\frac{1}{2}},
    \end{align*}
    which implies that
    \begin{align*}
        \int_\Omega \varphi^2 dd^cu \wedge \theta \leq 8\lVert u\rVert_{L^{\infty}(D)}\mathcal{E}^{\theta}(\varphi, \varphi). 
    \end{align*}
    
    We next deal with the general setting. 
    Let $D_1\Subset D$ be an open subset containing $\supp\varphi$.
    For $k\in\N$ large enough, the standard regularization $u^{[\frac{1}{k}]}\in SH_m(D)\cap \mathscr{C}^\infty(D,\R)$. Applying the first step on $u^{[\frac{1}{k}]}, D$, and $D_1$, we have
    \begin{align*}
        \int_D \varphi^2 dd^c u^{[\frac{1}{k}]}\wedge\theta\leq 8 \lVert u^{[\frac{1}{k}]}\rVert_{L^{\infty}(D_1)}\int_D d\varphi\wedge d^c\varphi \wedge\theta \leq 8 \lVert u\rVert_{L^{\infty}(D)}\mathcal{E}^{\theta}(\varphi, \varphi).
    \end{align*}
    By the dominated convergence theorem, we have
    \begin{align*}
        \int_D \varphi^2 dd^c u^{[\frac{1}{k}]}\wedge\theta= \int_D u^{[\frac{1}{k}]} dd^c \varphi^2\wedge\theta \to \int_D u dd^c \varphi^2\wedge\theta=\int_\Omega \varphi^2 dd^c u\wedge\theta
    \end{align*}
    as $k\to \infty$. This finishes the proof.
\end{proof}

Using the Poincar\'e-type inequality above, we obtain the following examples of admissible pairs, generalizing \cite[Theorem 2.5]{FO87} to the setting of $m$-subharmonic functions.
\begin{theorem}\label{Theorem:example of admissible pair}
    For any $\theta\in\mathcal{M}_m^+(\Omega)$, $(\theta, V)$ is an admissible pair.
\end{theorem}
\begin{proof}
    Choose $\delta>0$ and $w_1, \cdots, w_{m-1}\in SH_m(\Omega)\cap L^\infty(\Omega)$ so that
    \begin{align*}
        \theta=dd^c(w_1+2\delta |z|^2)\wedge\cdots\wedge dd^c (w_{m-1}+2\delta|z|^2)\wedge\omega^{n-m}.
    \end{align*}
    Denote $v_\ell:=w_\ell+\delta |z|^2$ for $\ell=1,\dots,m-1$. Then we have
    \begin{align*}
        \theta =& ~dd^c(v_1+\delta |z|^2)\wedge\cdots\wedge dd^c (v_{m-1}+\delta|z|^2)\wedge (dd^c|z|^2)^{n-m}\\
        =&~\delta^{m-1}(dd^c\lvert z\rvert^2)^{n-1} + \sum_{\ell = 1}^{m-1}\delta^{m-\ell-1}\sum_{1 \leq j_1 < \cdots < j_\ell \leq m-1}\theta^{j_1,\dots, j_\ell},
    \end{align*}
    where 
    \begin{align*}
        \theta^{j_1,\dots,j_\ell} := dd^cv_{j_1} \wedge \cdots \wedge dd^c v_{j_\ell} \wedge (dd^c\lvert z\rvert^2)^{n-\ell-1}.
    \end{align*}
    Notice that $\supp {\theta^{j_1,\dots,j_\ell}\wedge\omega}=\Omega$ since $\theta^{j_1,\dots,j_\ell}\wedge\omega\geq \delta^\ell (dd^c |z|^2)^n$. Hence,
    \begin{align}\label{Example:corollary of admissible pair}
        ((dd^c |z|^2)^{n-1}, V)\  \text{and} \ (\theta^{j_1,\dots,j_\ell}, \theta^{j_1,\dots,j_\ell}\wedge\omega) \ \text{are admissible pairs}
    \end{align}
    by reason of (\ref{Example:admissible pair}). For every $\varphi,\psi\in \mathscr{C}_c^\infty (\Omega,\R)$, we have
    \begin{align*}
        \mathcal{E}^{\theta}(\varphi, \psi) = \delta^{m-1}\mathcal{E}^{(0)}(\varphi, \psi) + \sum_{\ell = 1}^{m-1}\delta^{m-\ell-1}\sum_{1 \leq j_1 < \cdots < j_\ell \leq m-1}\mathcal{E}^{(j_1,\dots,j_\ell)}(\varphi,\psi),
    \end{align*}
    where
    \begin{align*}
        \mathcal{E}^{(0)}(\varphi, \psi) := \int_{\Omega}d\varphi \wedge d^c\psi \wedge (dd^c| z|^2)^{n-1}, \quad
        \mathcal{E}^{(j_1,\dots, j_\ell)}(\varphi, \psi) := \mathcal{E}^{\theta^{j_1,\dots, j_\ell}}(\varphi, \psi).
    \end{align*}
    
    Now assume that $(\varphi_k)_{k\in\N}\subset \mathscr{C}_c^{\infty}(\Omega,\R)$ constitutes an $\mathcal{E}^{\theta}$-Cauchy sequence and converges to $0$ in $L^2(\Omega)$. Then $(\varphi_k)_{k\in\N}$ is $\mathcal{E}^{(0)}$-Cauchy and $\mathcal{E}^{(j_1,\dots,j_\ell)}$-Cauchy for every multi-index $(j_1,\dots,j_\ell)$, and (\ref{Example:corollary of admissible pair}) implies that $\mathcal{E}^{(0)}(\varphi_k, \varphi_k) \to 0$ as $k \to
    \infty$. For any fixed multi-index $(j_1,\dots,j_\ell)$, it follows from Lemma \ref{Poincare type inequality}
    that
    \begin{align*}
        \int_{\Omega}\varphi_k^2 ~\theta^{j_\ell}\wedge dd^c| z|^2 =\int_\Omega \varphi_k^2 ~dd^c v_{j_\ell}\wedge (dd^c|z|^2)^{n-1}\leq 8\lVert v_{j_\ell}\rVert_{L^\infty(\Omega)}\mathcal{E}^{(0)}(\varphi_k,\varphi_k).
    \end{align*}
    Hence $\varphi_k \to 0$ in $L^2(\Omega; \theta^{j_\ell} \wedge dd^c\lvert z\rvert^2)$ and consequently $\mathcal{E}^{(j_\ell)}(\varphi_k,\varphi_k) \to 0$ as $k \to \infty$ by (\ref{Example:corollary of admissible pair}). Then Lemma \ref{Poincare type inequality} again reads
    \begin{align*}
        \int_{\Omega}\varphi_k^2~\theta^{j_{\ell-1},j_\ell}\wedge dd^c|z|^2 =& \int_{\Omega}\varphi_k^2~dd^c v_{j_{\ell-1}}\wedge \left(dd^c v_{j_\ell}\wedge \left(dd^c|z|^2\right)^{n-2}\right)\\\leq& ~8 \lVert v_{j_{\ell-1}}\rVert_{L^\infty(\Omega)}\mathcal{E}^{(j_\ell)}(\varphi_k, \varphi_k).
    \end{align*}
    Hence $\varphi_k \to 0$ in $L^2(\Omega; \theta^{j_{\ell-1},j_\ell} \wedge dd^c\lvert z\rvert^2)$ and consequently $\mathcal{E}^{(j_{\ell-1},j_\ell)}(\varphi_k,\varphi_k) \to 0$ as $k \to \infty$ by (\ref{Example:corollary of admissible pair}) again. Repeating the argument above, we finally have $\mathcal{E}^{(j_1,\dots,j_\ell)}(\varphi_k,\varphi_k) \to 0$ as $k \to \infty$. Therefore, $\mathcal{E}^\theta(\varphi_k,\varphi_k) \to 0$ as $k\to\infty$. This finishes the proof.
\end{proof}

Let $\theta$ be a closed positive $(n-1,n-1)$-current on $\Omega$, and $(\mathcal{F}^\theta,\mathcal{E}^\theta)$ be the Dirichlet space for an admissible pair $(\theta,\mu)$. We say that an extended real-valued function $v$ on $\Omega$ is in $\mathcal{F}_{loc}^\theta$ or locally in $\mathcal{F}^\theta$, if for every open subset $D\Subset\Omega$, there exists $u\in\mathcal{F}^\theta$ such that $u=v$ on $D$. Motivated by \cite[Theorem 7.6]{Nguyen26}, we characterize $\mathcal{F}^\theta_{loc}$.

\begin{proposition}\label{proposition:characterization for local domain of Dirichlet space}
    Let $(\mathcal{F}^\theta,\mathcal{E}^\theta)$ be the Dirichlet space for an admissible pair $(\theta,\mu)$. Let $v$ be an extended real-valued function on $\Omega$. Then $v\in \mathcal{F}^\theta_{loc}$ if and only if $\eta v\in\mathcal{F}^\theta$ for every $\eta\in \mathscr{C}_c^\infty(\Omega,\R)$.
\end{proposition}

\begin{proof}
    Assume first that $v\in\mathcal{F}_{loc}^\theta$. For every $\eta\in \mathscr{C}_c^\infty(\Omega,\R)$, there exists $u\in\mathcal{F}^\theta$ such that $u=v$ on an open neighborhood of $\supp \eta$. Hence $\eta u=\eta v$ on $\Omega$. Since $\mathscr{C}_c^\infty(\Omega,\R)$ is dense in the Hilbert space $(\mathcal{F}^\theta, \mathcal{E}^\theta_1)$, there exists $(\varphi_k)_{k\in\N}\subset \mathscr{C}_c^\infty(\Omega,\R)$ such that $\varphi_k\to u$ in $\mathcal{E}^\theta_1$.
    Take a constant $C>0$ so that $ d\eta\wedge d^c\eta\leq C dd^c |z|^2$ on $\Omega$.
    By (\ref{Equation:corollary of Schwarz inequality}) and the Poincar\'e-type inequality \cite[Theorem 1]{Okada82}, we have
    \begin{align*}
        &\mathcal{E}^\theta(\eta \varphi_k - \eta \varphi_\ell,\eta \varphi_k - \eta \varphi_\ell)\\
        &\leq 2 \int_\Omega \eta^2 d(\varphi_k - \varphi_\ell)\wedge d^c (\varphi_k - \varphi_\ell)\wedge\theta+2\int_\Omega (\varphi_k - \varphi_\ell)^2 d\eta \wedge d^c\eta\wedge\theta\\
        &\leq~ 2 \lVert\eta\rVert_{L^\infty(\Omega)}^2~\mathcal{E}^\theta(\varphi_k - \varphi_\ell,\varphi_k - \varphi_\ell)+2C\int_\Omega (\varphi_k - \varphi_\ell)^2 dd^c|z|^2\wedge\theta\\
        &\leq~ 2 \lVert\eta\rVert_{L^\infty(\Omega)}^2~\mathcal{E}^\theta(\varphi_k - \varphi_\ell,\varphi_k - \varphi_\ell)+16 C~\sup_\Omega |z|^2~ \mathcal{E}^\theta(\varphi_k - \varphi_\ell,\varphi_k - \varphi_\ell),
    \end{align*}
    which tends to $0$ as $k,\ell\to\infty$. Hence $(\eta \varphi_k)_{k\in\N}$ is $\mathcal{E}^\theta$-Cauchy. Since $\eta \varphi_k \to \eta u$ in $L^2(\Omega;\mu)$, $(\eta\varphi_k)_{k\in\N}\subset\mathcal{F}^\theta$ is $\mathcal{E}^\theta_1$-Cauchy. Therefore, $\eta u\in\mathcal{F}^\theta$ and $\eta \varphi_k \to \eta u$ in $\mathcal{E}^\theta_1$.
    
    Another direction follows directly from the definition of $\mathcal{F}^{\theta}_{loc}$. 
\end{proof}

We state an elementary fact for the convenience of proving Lemma \ref{Lemma:crucial Cauchy sequence}.

\begin{lemma}\label{Lemma:quasi-increasing sequence}
    Assume that $(a_k)_{k\in\N}, (b_{k,\ell})_{k,\ell\in\N}\subset\R$ satisfy the following conditions.
    \begin{itemize}
        \item[(i)] $\lim_{k,\ell\to\infty}b_{k,\ell}=0$.
        \item[(ii)] $(a_k)_{k\in\N}$ is bounded from above (resp. bounded from below).
        \item[(iii)] $a_\ell - a_k + b_{k,\ell}\geq 0$  when $\ell>k$ (resp. $a_\ell - a_k + b_{k,\ell}\leq 0$  when $\ell>k$).
    \end{itemize}
    Then $(a_k)_{k\in\N}$ has a finite limit.
\end{lemma}
\begin{proof}
    Assume first that $(a_k)_{k\in\N}$ is bounded from above and $a_\ell - a_k + b_{k,\ell}\geq 0$  when $\ell>k$.
    For every $\varepsilon>0$, there exists $N\in\N$ such that $|b_{k,\ell}|\leq\varepsilon$ for every $k, \ell\geq N$. Then we have $a_\ell - a_k +\varepsilon\geq 0$ whenever $\ell>k\geq N$. Hence, $(a_k)_{k\in\N}$ is bounded from below, and $\liminf_{\ell\to\infty}a_\ell\geq a_k -\varepsilon$ for every $k\geq N$. By letting $k\to\infty$, we deduce that $\liminf_{\ell\to\infty}a_\ell\geq \limsup_{k\to\infty} a_k -\varepsilon$. Therefore, $(a_k)_{k\in\N}$ has a finite limit. Under another assumption, the proof is similar. 
\end{proof}

The following result is crucial in this paper, generalizing \cite[Lemma 2.1]{FO87} to the setting of $m$-subharmonic functions.
\begin{lemma}\label{Lemma:crucial Cauchy sequence}
    Let $\theta \in \mathcal{M}_m(\Omega)$, $D\subset\Omega$ be an open subset, and $\eta \in \mathscr{C}_c^{\infty}(D,\R)$. Assume that $(u_k)_{k\in\N}\subset SH_m(D) \cap \mathscr{C}^{\infty}(D,\R)$ is decreasing and locally uniformly bounded. Then $(\eta u_k)_{k\in\N}$ is a Cauchy sequence with respect to $\mathcal{E}^{\theta}$, i.e.,  $$\lim_{k,\ell\to\infty}\mathcal{E}^{\theta}(\eta u_k - \eta u_\ell,\eta u_k-\eta u_\ell)=0.$$
\end{lemma}
\begin{proof}
    Since $\eta u_k=0$ outside $\supp \eta$ for every $k\in\N$, we can assume that $(u_k)_{k\in\N}$ is uniformly bounded by shrinking $D$, and $0 \leq u_k \leq M$ on $D$ for some constant $M$. We can further assume that $\Omega=D$, since $\mathcal{E}^{\theta}(\eta u_k - \eta u_\ell,\eta u_k-\eta u_\ell)$ depends only on the restrictions of $\theta$ and  $\eta u_k - \eta u_\ell$ on $D$.
    
    Since $dd^cu_k \wedge \theta$ is a positive current on $D$, we have
    \begin{align}\label{Equation: crucial inequality}
        dd^c(u_k^2)\wedge \theta = 2du_k \wedge d^cu_k \wedge \theta + 2u_kdd^cu_k \wedge \theta \geq 2du_k \wedge d^cu_k \wedge \theta\geq 0
    \end{align}
    in the sense of currents. It follows that
    \begin{align*}
        2\int_{D}\eta^2 du_k \wedge d^cu_k \wedge \theta 
        &\leq \int_D \eta^2 dd^c(u_k^2) \wedge \theta \\ &= \int_D u_k^2dd^c(\eta^2)\wedge \theta
        \\
    &\leq CM^2\int_{\supp \eta}\theta\wedge dd^c|z|^2,
    \end{align*}
    where $C>0$ is a constant depending only on $n$ and $\eta$.
    For every $k,\ell\in\N$, set
    \begin{align*}
        a_k:=&~\mathcal{E}^\theta(\eta u_k,\eta u_k),\\
        b_{k, \ell} :=& \int_D \eta(u_\ell-u_k)u_k dd^c\eta \wedge \theta + 2 \int_D \eta(u_\ell-u_k)d\eta \wedge d^cu_k \wedge \theta.
    \end{align*}
    By the last estimate and (\ref{Equation:corollary of Schwarz inequality}), we have
    \begin{align*}
        a_k&\leq 2\int_D \eta^2 d u_k \wedge d^c u_k \wedge \theta +  2\int_D u_k^2 d\eta \wedge d^c\eta \wedge \theta\\
        &\leq CM^2\int_{\supp \eta}\theta\wedge dd^c|z|^2 +2 M^2 \int_D d\eta\wedge d^c \eta\wedge\theta,
    \end{align*}
    and hence $\sup_{k\in\N}a_k < \infty$ follows.
    By the Schwarz inequality (\ref{Equation:Schwarz inequality}), we have
    \begin{align*}
        &\left(\int_D \eta(u_\ell-u_k)d\eta \wedge d^cu_k \wedge \theta \right)^2\\ &\leq \left(\int_D\eta^2du_k \wedge d^cu_k \wedge \theta\right)\left(\int_D (u_\ell-u_k)^2d\eta \wedge d^c\eta \wedge \theta\right)\\
        &\leq \frac{1}{2}CM^2\int_{\supp \eta}\theta\wedge dd^c|z|^2\int_D (u_\ell-u_k)^2d\eta \wedge d^c\eta \wedge \theta,
    \end{align*}
    which tends to $0$ as $k,\ell\to\infty$ following from the dominated convergence theorem. By the dominated convergence theorem again, we have
    \begin{align*}
        \lim_{k,\ell\to\infty}\int_D \eta(u_\ell-u_k)u_k dd^c\eta \wedge \theta=0.
    \end{align*}
    Therefore, $\lim_{k,\ell\to\infty}b_{k,\ell}=0$.
    
    Since $(u_k)_{k\in\N}$ is decreasing, for every $k < \ell$ we have
        \begin{align*}
            0&\leq \mathcal{E}^{\theta}(\eta u_k - \eta u_\ell, \eta u_k - \eta u_\ell) \\
            &= a_\ell - a_k + 2 \mathcal{E}^{\theta}(\eta u_k, \eta u_k-\eta u_\ell)\\
            &= a_\ell - a_k  + 2b_{k,\ell} +2\int_D\eta^2(u_\ell-u_k) dd^cu_k \wedge \theta\\
            &\leq a_\ell - a_k  + 2b_{k,\ell}.
        \end{align*}
    Since $(a_k)_{k\in\N}$ and $(2b_{k,\ell})_{k,\ell\in\N}$ satisfy the conditions in Lemma \ref{Lemma:quasi-increasing sequence}, $(a_k)_{k\in\N}$ has a finite limit. We then deduce that $\mathcal{E}^{\theta}(\eta u_k - \eta u_\ell, \eta u_k - \eta u_\ell)\to 0$ as $k,\ell\to \infty$. 
\end{proof}

The next result is an immediate consequence of Lemma \ref{Lemma:crucial Cauchy sequence}, and generalizes \cite[Theorem 2.2]{FO87} to the setting of $m$-subharmonic functions.
\begin{theorem}\label{bounded msh in F loc}
    Let $\theta \in \mathcal{M}_m(\Omega)$, and $(\mathcal{F}^{\theta}, \mathcal{E}^{\theta})$ be the Dirichlet space for an admissible pair $(\theta,\mu)$. Assume that $u \in SH_m(\Omega) \cap L^{\infty}_{loc}(\Omega)$ and $\eta\in \mathscr{C}_c^\infty(\Omega,\R)$. Then the following conclusions hold.
    \begin{itemize}
        \item [(i)] $\eta u\in\mathcal{F}^\theta$, and $\eta u^{[\delta]}$ converges to $\eta u$ in $\mathcal{E}_1^{\theta}$ as $\delta \downarrow 0$.
        \item [(ii)] $u$ is $\mathcal{E}^{\theta}$-q.c.
        \item[(iii)] Assume in addition that $v\in SH_m(\Omega) \cap L^{\infty}_{loc}(\Omega)$, $\xi \in \mathscr{C}_c^{\infty}(\Omega,\R)$, and $\xi = 1$ on an open neighborhood of $\supp \eta$. Then
        \begin{align*}
            \mathcal{E}^{\theta}(\eta u, \xi v) = -\int_\Omega \eta u dd^c v \wedge \theta.
        \end{align*}
    \end{itemize}
\end{theorem}

\begin{proof}
    (i) Let $D\Subset \Omega$ be an open subset containing $\supp \eta$. For an arbitrary sequence $(\delta_k)_{k\in\N}\subset\R_{>0}$ decreasing to $0$, $(u^{[\delta_k]})_{k\in\N}\subset SH_m(D)\cap \mathscr{C}^\infty(D,\R)$ is decreasing and uniformly bounded on $D$. By Lemma \ref{Lemma:crucial Cauchy sequence}, $(\eta u^{[\delta_k]})_{k\in\N}$ is $\mathcal{E}^\theta$-Cauchy. Moreover, $\eta u^{[\delta_k]}$ converges to $\eta u$ in $L^2(\Omega;\mu)$, since
    \begin{align*}
        \int_\Omega (\eta u^{[\delta_k]}-\eta u)^2 d\mu\leq \lVert\eta\rVert_{L^\infty(\Omega)}^2\int_{\supp\eta}(u^{[\delta_k]}-u)^2 d\mu,
    \end{align*}
    which tends to $0$ as $k\to\infty$ by the dominated convergence theorem.   
    We deduce that $(\eta u^{[\delta_k]})_{k\in\N}$ is a Cauchy sequence in the Hilbert space $(\mathcal{F}^\theta,\mathcal{E}_1^\theta)$, and then by completeness there exists $v\in \mathcal{F}^\theta$ such that 
    \begin{align*}
        \mathcal{E}^\theta_1(\eta u^{[\delta_k]}-v,\eta u^{[\delta_k]}-v)=\mathcal{E}^\theta(\eta u^{[\delta_k]}-v,\eta u^{[\delta_k]}-v)+\int_\Omega (\eta u^{[\delta_k]}-v)^2d\mu\to 0
    \end{align*}
    as $k\to\infty$.
    Therefore, $\eta u=v\in\mathcal{F}^\theta$, and $\mathcal{E}^\theta_1(\eta u^{[\delta_k]}-\eta u,\eta u^{[\delta_k]}-\eta u )\to 0$ as $k\to\infty$.

    (ii) Fix an $\varepsilon>0$. Let $D$ be an arbitrary relatively compact open subset of $\Omega$. Let $\xi\in \mathscr{C}_c^\infty(\Omega,\R)$ so that $\xi=1$ on $D$. For every $k\in\N$ large enough and $\delta>0$,
    \begin{align*}
        D\cap\{\xi u^{[\frac{1}{k}]} -\xi u>\delta\}=D\cap\{u^{[\frac{1}{k}]}-u>\delta\}
    \end{align*}
    is an open subset of $D$ and hence is an open subset of $\Omega$. Then, we have
    \begin{align*}
        \mathrm{Cap}_\theta(D\cap\{\xi u^{[\frac{1}{k}]} -\xi u>\delta\})&\leq\mathcal{E}^\theta_1(\frac{\xi u^{[\frac{1}{k}]}-\xi u}{\delta},\frac{\xi u^{[\frac{1}{k}]}-\xi u}{\delta})\\
        &=\frac{1}{\delta^2}\mathcal{E}^\theta_1(\xi u^{[\frac{1}{k}]}-\xi u, \xi u^{[\frac{1}{k}]}-\xi u), 
    \end{align*}
    which tends to $0$ as $k\to\infty$.
    Hence, there exists a subsequence $(j_k)_{k\in\N}$ of $\N$ such that
    \begin{align*}
        \mathrm{Cap}_\theta(D\cap\{\xi u^{[\frac{1}{j_k}]} -\xi u>\frac{1}{k}\})\leq \frac{\varepsilon}{2^k}.
    \end{align*}
    Set
    \begin{align*}
        G:=\bigcup_{k\in\N}\left(D\cap\left\{\xi u^{[\frac{1}{j_k}]}-\xi u>\frac{1}{k}\right\}\right).
    \end{align*}
    Then $G$ is an open subset of $D$, $\mathrm{Cap}_\theta(G)\leq \varepsilon$, and $u^{[\frac{1}{j_k}]}-u\leq \frac{1}{k}$ on $D\setminus G$ for every $k\in\N$. Therefore, $u$ is a finite-valued continuous function on $D\setminus G$ since $u^{[\frac{1}{j_k}]}$ uniformly converges to $u$ on it.
    
    Let $(D_k)_{k\in\N}$ be a relatively compact exhaustion of $\Omega$. For every $k\in\N$, there exists an open subset $G_k$ of $D_k$ such that $\mathrm{Cap}_\theta(G_k)<\frac{\varepsilon}{2^k}$ and $u$ is a finite-valued continuous function on $D_k\setminus G_k$. Therefore, $\mathrm{Cap}_\theta(\tilde{G})<\varepsilon$ and $u$ is a finite-valued continuous function on $\Omega\setminus \tilde{G}$, where $\tilde{G}:=\cup_{k\in\N}G_k$ is an open subset of $\Omega$.
    
    (iii) Let $E\Subset \Omega$ be an open subset such that $\supp \eta\subset E$ and $\xi=1$ on $E$. By (i) and the dominated convergence theorem, we have
    \begin{align*}
        \mathcal{E}^{\theta}(\eta u, \xi v) &= \lim_{k\to\infty}\lim_{\ell\to\infty}\mathcal{E}^{\theta}(\eta u^{[\frac{1}{k}]}, \xi v^{[\frac{1}{\ell}]})\\
        &=-\lim_{k\to\infty}\lim_{\ell\to\infty}\int_\Omega \eta u^{[\frac{1}{k}]} dd^c (\xi v^{[\frac{1}{\ell}]})\wedge\theta\\
        &=-\lim_{k\to\infty}\lim_{\ell\to\infty}\int_E \eta u^{[\frac{1}{k}]} dd^c v^{[\frac{1}{\ell}]}\wedge\theta \\
        &=-\lim_{k\to\infty}\int_E \eta u^{[\frac{1}{k}]} dd^c v\wedge\theta\\
        &=-\int_\Omega \eta u dd^c v \wedge \theta.
    \end{align*}
    This finishes the proof.
\end{proof}

Motivated by \cite[the proof of Theorem 2.3]{FO87}, we next establish a weak form of monotonicity for the Dirichlet form $\mathcal{E}^{\theta}$ when $(\theta,\mu)$ is an admissible pair. This property will be used repeatedly throughout this section. The complete form of the monotonicity property will be proved in Theorem \ref{Theorem: Monotonicity for admissible pairs}.
\begin{lemma}\label{Lemma: weak monotonicity for admissible pairs}
    Let $\theta\in\mathcal{M}_m(\Omega)$, and $(\mathcal{F}^{\theta},\mathcal{E}^{\theta})$ be the Dirichlet space for an admissible pair $(\theta,\mu)$. Let $u_1,u_2,v_1,v_2\in SH_m(\Omega)\cap L^{\infty}_{loc}(\Omega)$ and $K\subset\Omega$ be a compact set. Assume that
    \begin{align*}
        u_k\leq v_k~ \text{on} ~K, \quad u_k=v_k ~\text{on} ~\Omega\setminus K
    \end{align*}
    for every $k$. Then for every $\xi,\gamma \in \mathscr{C}_c^{\infty}(\Omega,\R)$ such that $\xi = \gamma=1$ on an open neighborhood of $K$, we have
    \begin{align*}
        \mathcal{E}^{\theta}(\xi u_1, \gamma u_2) \geq \mathcal{E}^{\theta}(\xi v_1, \gamma v_2).
    \end{align*}
\end{lemma}

\begin{proof}
    Let $E\Subset \Omega$ be an open subset such that $K\subset E$ and $\xi =\gamma= 1$ on $E$. For every $\delta>0$ small enough, $dd^c u_2^{[\delta]}\wedge\theta \geq0$ on $E$ and $u_1^{[\delta]}=v_1^{[\delta]}$ on $\Omega\setminus E$, hence
    \begin{align*}
        \mathcal{E}^{\theta}(\xi u_1^{[\delta]}, \gamma u_2^{[\delta]}) &= -\int_{\Omega}\xi u_1^{[\delta]} dd^c(\gamma u_2^{[\delta]}) \wedge \theta\\
        &= -\int_E \xi u_1^{[\delta]}dd^c(\gamma u_2^{[\delta]})\wedge \theta -\int_{\Omega\setminus E}\xi u_1^{[\delta]}dd^c(\gamma u_2^{[\delta]})\wedge \theta\\
        &= -\int_E u_1^{[\delta]}dd^c u_2^{[\delta]}\wedge \theta -\int_{\Omega\setminus E}\xi u_1^{[\delta]}dd^c(\gamma u_2^{[\delta]})\wedge \theta\\
        &\geq -\int_E v_1^{[\delta]}dd^c u_2^{[\delta]}\wedge \theta -\int_{\Omega\setminus E}\xi u_1^{[\delta]}dd^c(\gamma u_2^{[\delta]})\wedge \theta\\
        &= -\int_E v_1^{[\delta]}dd^c u_2^{[\delta]}\wedge \theta -\int_{\Omega\setminus E}\xi v_1^{[\delta]}dd^c(\gamma u_2^{[\delta]})\wedge \theta \\
        &= \mathcal{E}^{\theta}(\xi v_1^{[\delta]}, \gamma u_2^{[\delta]}).
    \end{align*}
    By the same computation as above, we obtain
    \[
        \mathcal{E}^{\theta}(\xi u_1^{[\delta]}, \gamma u_2^{[\delta]})\geq \mathcal{E}^{\theta}(\xi v_1^{[\delta]}, \gamma u_2^{[\delta]})\geq \mathcal{E}^{\theta}(\xi v_1^{[\delta]}, \gamma v_2^{[\delta]}).
    \]
    Letting $\delta\to0$ and applying Theorem \ref{bounded msh in F loc} (i), we obtain the desired assertion.
\end{proof}

The next lemma will simplify the proofs of several subsequent results.
\begin{lemma}[Reduction process]\label{Lemma: reduction process for admissible pairs}
    Let $\Omega$ be a bounded strictly $m$-pseudoconvex domain in $\C^n$, $\theta \in \mathcal{M}_m(\Omega)$, and $(\mathcal{F}^{\theta}, \mathcal{E}^{\theta})$ be the Dirichlet space for an admissible pair $(\theta, \mu)$. Assume that $(u_k)_{k\in \N} \subset SH_m(\Omega)$ and $u \in SH_m(\Omega)$ are uniformly bounded and have a uniform negative upper bound. Then for every $\eta \in \mathscr{C}_c^{\infty}(\Omega, \R)$, there exist uniformly bounded functions $(\hat{u}_k)_{k \in \N} \subset SH_m(\Omega)$ and $\hat{u} \in SH_m(\Omega)$, an open set $U$ and a compact set $K$ that satisfy the following properties.
    \begin{itemize}
        \item [(i)] $\supp\eta \subset U \subset K \subset \Omega$, $\hat{u}_k = u_k$ on $U$ for every $k \in \N$, $\hat{u} = u$ on $U$, $\hat{u}_k = \hat{u}$ on $\Omega \setminus K$ for every $k \in \N$, and $(\hat{u}_k)_{k\in\N}$ is increasing (resp. decreasing) whenever $(u_k)_{k\in\N}$ is increasing (resp. decreasing).
        \item [(ii)] For every $\xi \in \mathscr{C}_c^{\infty}(\Omega, \R)$ with $\xi = 1$ on an open neighborhood of $K$, if $(\xi \hat{u}_k)_{k \in \N}$ is $\mathcal{E}^{\theta}$-Cauchy, then so is $(\eta u_k)_{k \in \N}$. Moreover, if $\xi \hat{u}_k$ converges to $\xi \hat{u}$ in $\mathcal{E}^{\theta}$, then $\eta u_k$ converges to $\eta u$ in $\mathcal{E}^{\theta}$.
    \end{itemize}
\end{lemma}

\begin{proof}
    By Proposition \ref{Proposition:localization principle}, there exist $(\hat{u}_k)_{k\in\N}\subset SH_m(\Omega)$, $\hat{u}\in SH_m(\Omega)$, $U$, and $K$ that satisfy the property (i).

    For every $\xi \in \mathscr{C}_c^{\infty}(\Omega, \R)$ with $\xi = 1$ on an open neighborhood of $K$, let $C > 0$ be a constant so that $2 d\eta\wedge d^c\eta\leq C dd^c |z|^2$. For every $k,\ell\in\N$, by Theorem \ref{bounded msh in F loc} (i) we have
    \begin{align*}
        \mathcal{E}^{\theta}(\eta(u_k-u_{\ell}), \eta (u_k-u_{\ell})) &= \lim_{\delta \rightarrow 0}\mathcal{E}^{\theta}(\eta (u_k^{[\delta]} - u_{\ell}^{[\delta]}), \eta (u_k^{[\delta]} - u_{\ell}^{[\delta]})) \\
        &= \lim_{\delta \rightarrow 0} \mathcal{E}^{\theta}(\eta \xi (\hat{u}_k^{[\delta]} - \hat{u}_{\ell}^{[\delta]}), \eta \xi (\hat{u}_k^{[\delta]} - \hat{u}_{\ell}^{[\delta]})).
    \end{align*}
    By (\ref{Equation:corollary of Schwarz inequality}) and Lemma \ref{Poincare type inequality}, we have
    \begin{align*}
        &\mathcal{E}^{\theta}(\eta\xi(\hat{u}_k^{[\delta]}-\hat{u}_{\ell}^{[\delta]}), \eta\xi(\hat{u}_k^{[\delta]}-\hat{u}_{\ell}^{[\delta]}))\\
        &\leq 2\int_\Omega \xi^2 (\hat{u}_k^{[\delta]}-\hat{u}_{\ell}^{[\delta]})^2 d\eta\wedge d^c\eta\wedge\theta+ 2\int_\Omega \eta^2 d (\xi(\hat{u}_k^{[\delta]}-\hat{u}_{\ell}^{[\delta]}))\wedge d^c(\xi(\hat{u}_k^{[\delta]}-\hat{u}_{\ell}^{[\delta]}))\wedge\theta\\
        &\leq C\int_\Omega \xi^2 (\hat{u}_k^{[\delta]}-\hat{u}_{\ell}^{[\delta]})^2 dd^c|z|^2\wedge\theta+ 2\int_\Omega \eta^2 d (\xi(\hat{u}_k^{[\delta]}-\hat{u}_{\ell}^{[\delta]}))\wedge d^c(\xi(\hat{u}_k^{[\delta]}-\hat{u}_{\ell}^{[\delta]}))\wedge\theta\\
        &\leq (8 C \sup_\Omega \lvert z\rvert^2+2 \lVert\eta\rVert^2_{L^{\infty}(\Omega)})\mathcal{E}^{\theta}(\xi (\hat{u}_k^{[\delta]}-\hat{u}_{\ell}^{[\delta]}),\xi (\hat{u}_k^{[\delta]}-\hat{u}_{\ell}^{[\delta]})).
    \end{align*}
    Letting $\delta\to 0$, we obtain
    \begin{align*}
        \mathcal{E}^{\theta}(\eta (u_k-u_{\ell}), \eta (u_k-u_{\ell})) \leq (8 C \sup_\Omega |z|^2+2 \lVert\eta\rVert^2_{L^{\infty}(\Omega)}) \mathcal{E}^{\theta}(\xi (\hat{u}_k - \hat{u}_{\ell}), \xi (\hat{u}_k - \hat{u}_{\ell})).
    \end{align*}
    Therefore, if $(\xi \hat{u}_k)_{k \in \N}$ is $\mathcal{E}^{\theta}$-Cauchy, then $(\eta u_k)_{k \in \N}$ is $\mathcal{E}^{\theta}$-Cauchy. Repeating the same computation by replacing $u_{\ell}$ with $u$, we have
    \begin{align*}
        \mathcal{E}^{\theta}(\eta (u_k-u), \eta (u_k-u)) \leq (8 C \sup_\Omega |z|^2+2 \lVert\eta\rVert^2_{L^{\infty}(\Omega)}) \mathcal{E}^{\theta}(\xi (\hat{u}_k - \hat{u}), \xi (\hat{u}_k - \hat{u})).
    \end{align*}
    Therefore, if $\xi \hat{u}_k$ converges to $\xi \hat{u}$ in $\mathcal{E}^{\theta}$, then $\eta u_k$ converges to $\eta u$ in $\mathcal{E}^{\theta}$.
\end{proof}

We now establish a theorem on $\mathcal{E}^\theta_1$-convergence of monotone sequences of $m$-subharmonic functions, generalizing \cite[Theorem 2.3]{FO87} to the $m$-subharmonic setting.
\begin{theorem}\label{Theorem:Cauchy for locally uniformly bounded and monotone}
    Let $\theta\in\mathcal{M}_m(\Omega)$, and $(\mathcal{F}^{\theta}, \mathcal{E}^{\theta})$ be the Dirichlet space for an admissible pair $(\theta,\mu)$. Assume that $(u_k)_{k\in\N}\subset SH_m(\Omega)$ is monotone, i.e., increasing or decreasing, locally uniformly bounded, and $\lim_{k\to\infty}u_k=u$ $\mu$-a.e. on $\Omega$ for some $u\in SH_m(\Omega)\cap L_{loc}^\infty(\Omega)$. Then $\eta u_k$ converges to $\eta u$ in $\mathcal{E}_1^{\theta}$ for every $\eta \in \mathscr{C}_c^{\infty}(\Omega,\R)$.
\end{theorem}

\begin{proof}
    It suffices to prove that $(\eta u_k)_{k\in\N}$ is $\mathcal{E}^\theta$-Cauchy, since $\eta u_k$ converges to $\eta u$ in $L^2(\Omega;\mu)$. Moreover, due to the partition of unity, it suffices to prove that every $z\in\Omega$ admits an open neighborhood $D\subset\Omega$ such that $(\eta u_k)_{k\in\N}$ is $\mathcal{E}^\theta$-Cauchy for every $\eta\in \mathscr{C}_c^\infty(D,\R)$. Let $(\mathcal{F}^{\theta|_D},\mathcal{E}^{\theta|_D})$ be the Dirichlet space for the admissible pair $(\theta|_D,\mu|_D)$. By Theorem \ref{bounded msh in F loc} (i), for every $\eta\in \mathscr{C}_c^\infty(D,\R)$ we have
    \begin{align}\label{equation: admissible pair and its restriction}
        \mathcal{E}^\theta(\eta u_k -\eta u_\ell,\eta u_k -\eta u_\ell)&=\lim_{\delta\to0}\mathcal{E}^\theta(\eta u_k^{[\delta]} -\eta u_\ell^{[\delta]},\eta u_k^{[\delta]} -\eta u_\ell^{[\delta]})\\
        &=\lim_{\delta\to0}\mathcal{E}^{\theta|_D}(\eta u_k^{[\delta]} -\eta u_\ell^{[\delta]},\eta u_k^{[\delta]} -\eta u_\ell^{[\delta]})\notag\\
        &=\mathcal{E}^{\theta|_D}(\eta u_k -\eta u_\ell,\eta u_k -\eta u_\ell).\notag
    \end{align}
    Therefore, we can assume that $\Omega$ is a bounded strictly $m$-pseudoconvex domain, and $(u_k)_{k\in\N}$, $u$ are uniformly bounded on $\Omega$ with a uniform negative upper bound.
    
    Fix an arbitrary $\eta\in \mathscr{C}_c^{\infty}(\Omega,\R)$. Apply Lemma \ref{Lemma: reduction process for admissible pairs} to $(u_k)_{k \in \N}$, $u$, and $\eta$. Then we obtain $U$, $K$, $(\hat{u}_k)_{k \in \N}$, and $\hat{u}$ that satisfy the property (i) of Lemma \ref{Lemma: reduction process for admissible pairs}. Take $\xi \in \mathscr{C}_c^{\infty}(\Omega, \R)$ so that $\xi = 1$ on an open neighborhood of $K$. 
    
    Assume that $(\hat{u}_k)_{k\in\N}$ is increasing. For every $k <\ell$, by Lemma \ref{Lemma: weak monotonicity for admissible pairs} we have 
    \begin{align*}
        \mathcal{E}^\theta(\xi \hat{u}_k, \xi \hat{u}_\ell) \geq \mathcal{E}^\theta(\xi \hat{u}_\ell , \xi \hat{u}_\ell).
    \end{align*}
    Hence,
    \begin{align*}
        0\leq \mathcal{E}^\theta(\xi \hat{u}_k - \xi \hat{u}_\ell, \xi \hat{u}_k - \xi \hat{u}_\ell) \leq \mathcal{E}^\theta(\xi \hat{u}_k, \xi \hat{u}_k) - \mathcal{E}^\theta(\xi \hat{u}_\ell , \xi \hat{u}_\ell).
    \end{align*}
    We deduce that $\left(\mathcal{E}^{\theta}(\xi \hat{u}_k, \xi \hat{u}_k)\right)_{k\in\N}$ is decreasing and hence $(\xi \hat{u}_k)_{k\in\N}$ is $\mathcal{E}^\theta$-Cauchy. 
    
    Assume that $(\hat{u}_k)_{k\in\N}$ is decreasing. Then $v := \lim_{k \rightarrow \infty}\hat{u}_k\in SH_m(\Omega) \cap L^{\infty}(\Omega)$. By Lemma \ref{Lemma: weak monotonicity for admissible pairs}, for every $k <\ell$ we have
    \begin{align*}
        \mathcal{E}^\theta(\xi \hat{u}_k, \xi \hat{u}_\ell) \geq \mathcal{E}^\theta(\xi \hat{u}_k, \xi \hat{u}_k),\quad \mathcal{E}^\theta(\xi \hat{u}_k, \xi \hat{u}_k)\leq \mathcal{E}^\theta(\xi v, \xi v).
    \end{align*}
    Hence,
    \begin{align*}
        0\leq\mathcal{E}^\theta(\xi \hat{u}_k - \xi \hat{u}_\ell, \xi \hat{u}_k - \xi \hat{u}_\ell) \leq \mathcal{E}^\theta(\xi \hat{u}_\ell, \xi \hat{u}_\ell) - \mathcal{E}^\theta(\xi \hat{u}_k, \xi \hat{u}_k).
    \end{align*}
    We deduce that $\left(\mathcal{E}^{\theta}(\xi \hat{u}_k, \xi \hat{u}_k)\right)_{k\in\N}$ is increasing, bounded from above by $\mathcal{E}^\theta(\xi v, \xi v)$, and hence it has a finite limit. Therefore, $(\xi \hat{u}_k)_{k\in\N}$ is $\mathcal{E}^\theta$-Cauchy.

    By the property (ii) of Lemma \ref{Lemma: reduction process for admissible pairs}, $(\eta u_k)_{k\in\N}$ is $\mathcal{E}^\theta$-Cauchy.
\end{proof}

\begin{remark}\label{Remark: after theorem on Cauchy for locally uniformly bounded and monotone}
    Let $w:=\lim_{k\to\infty}u_k$. Then $u=w~ \mu$-a.e. on $\Omega$, and $u_k,u$ are $\mathcal{E}^\theta$-q.c. By \cite[Theorem 2.1.4 (ii)]{FOT11}, $w$ is $\mathcal{E}^\theta$-q.c. Then \cite[Lemma 2.1.4]{FOT11} implies that $u=w$ $\mathcal{E}^\theta$-q.e. on $\Omega$.
\end{remark}

We now generalize \cite[Theorem 2.4]{FO87} to the $m$-subharmonic setting. Let $E\subset\Omega$ be a subset. We say that $\mu$ does not charge $E$ if $E$ is contained in a Borel set of zero measure $\mu$. For example, $V$ does not charge any $m$-negligible set.
\begin{theorem}\label{Theorem:m-negligible is E theta-polar}
    Let $\theta\in\mathcal{M}_m(\Omega)$, and $(\mathcal{F}^{\theta}, \mathcal{E}^{\theta})$ be the Dirichlet space for an admissible pair $(\theta, \mu)$. Then $\mu$ does not charge any $m$-negligible set if and only if every $m$-negligible set is $\mathcal{E}^{\theta}$-polar.
\end{theorem}
\begin{proof}
    (i) Assume that an $m$-negligible set $E$ is $\mathcal{E}^{\theta}$-polar. Then
    \begin{align*}
         0&=\mathrm{Cap}_\theta(E)\\
         &= \inf\left\{\mathrm{Cap}_\theta(F): E \subset F\subset\Omega, F~\text{Borel}\right\}\\
         &\geq \inf\left\{\mu(F): E \subset F\subset\Omega, F~\text{Borel}\right\}.
    \end{align*}
    Therefore, $\mu$ does not charge $E$.

    (ii) Assume that $\mu$ does not charge any $m$-negligible set. Let $E$ be an $m$-negligible set. By Choquet's lemma, there exists an increasing sequence $(u_k)_{k\in\N}\subset SH_m(\Omega)$ such that $(u_k)_{k\in\N}$ is locally uniformly bounded from above and $E\subset\{u<u^*\}$, where $u:=\lim_{k\to\infty}u_k$.

    Let $u_{k,r}:=\max(u_k,r)$ for every $k\in\N$ and $r\in\Q$. Then for every $r\in\Q$, $(u_{k,r})_{k\in\N}\subset SH_m(\Omega)$ is increasing and locally uniformly bounded. Let $u_r:=\lim_{k\to\infty}u_{k,r}=\max(u,r)$. Then $(u_r)^*\in SH_m(\Omega)\cap L^\infty_{loc}(\Omega)$. The $m$-negligible set $\{u_r<(u_r)^*\}$ is Borel and hence has zero measure $\mu$. Theorem \ref{Theorem:Cauchy for locally uniformly bounded and monotone} and Remark \ref{Remark: after theorem on Cauchy for locally uniformly bounded and monotone} imply that $\{u_r<(u_r)^*\}$ is $\mc E^\theta$-polar.

    Since $\max(u^*,r)$ is u.s.c. and $\max(u^*,r)\geq u_r$, we have $\max(u^*,r)\geq (u_r)^*$. On the other hand, we have $\max(u^*,r)=u_r\leq (u_r)^*$ outside $\{u<u^*\}$, which is of zero measure $V$. We deduce that $\max(u^*,r)$ and $(u_r)^*$ are equal outside a set of zero measure $V$ and hence they are equal everywhere by the strong upper semi-continuity of $m$-subharmonic functions. We then have
    \begin{align*}
        E\subset\left\{u<u^*\right\}\subset \bigcup_{r\in\Q} \left\{u_r<\max(u^*,r)\right\}=\bigcup_{r\in\Q} \left\{u_r<(u_r)^*\right\}.
    \end{align*}
    Therefore, $\mathrm{Cap}_\theta (E)=0$.
\end{proof}

By Theorem \ref{Theorem:example of admissible pair} and Theorem \ref{Theorem:m-negligible is E theta-polar}, we immediately have the following result.
\begin{corollary}\label{Corollary: m-negligible are Etheta-polar when mu is V}
    Let $\theta\in\mathcal{M}^+_m(\Omega)$, and $(\mathcal{F}^{\theta}, \mathcal{E}^{\theta})$ be the Dirichlet space for the admissible pair $(\theta, V)$. Then every $m$-negligible set is $\mc E^\theta$-polar.
\end{corollary}

We conclude this section by discussing how Theorem \ref{Theorem:example of admissible pair} applies to the study of symmetric forms induced by bounded $m$-subharmonic functions.

Let $w_1,\dots,w_{m-1}\in SH_m(\Omega)\cap L^\infty(\Omega)$ and $\delta>0$. Set
\begin{align*}
    \theta&:= dd^c w_1\wedge \dots\wedge dd^c w_{m-1}\wedge\omega^{n-m},\\
    \theta_\delta&:= dd^c (w_1+ \delta |z|^2)\wedge \dots\wedge dd^c (w_{m-1}+ \delta |z|^2)\wedge\omega^{n-m}.
\end{align*}
Let $(\mathcal{F}^{\theta_\delta}, \mathcal{E}^{\theta_\delta})$ be the Dirichlet space for the admissible pair $(\theta_\delta,V)$, whose admissibility follows from Theorem \ref{Theorem:example of admissible pair}. In general, $(\theta,V)$ need not be an admissible pair. Notice that
\begin{align*}
    \mathcal{E}^\theta(\varphi,\varphi)\leq \mathcal{E}^{\theta_\delta} (\varphi,\varphi),\quad \forall \varphi\in \mathscr{C}_c^\infty(\Omega,\R).
\end{align*}
Hence, the bilinear form $\mathcal{E}^\theta$ on $\mathscr{C}_c^\infty(\Omega,\R)$ can be extended to a bilinear form on $\mathcal{F}^{\theta_\delta}$. More precisely, for every $u,v\in\mathcal{F}^{\theta_\delta}$, let $(\varphi_k)_{k\in\N},(\psi_k)_{k\in\N}\subset \mathscr{C}_c^\infty(\Omega,\R)$ so that $\varphi_k$ converges to $u$ in $\mathcal{E}^{\theta_\delta}_1$ and $\psi_k$ converges to $v$ in $\mathcal{E}^{\theta_\delta}_1$. Then, both $(\varphi_k)_{k\in\N}$ and $(\psi_k)_{k\in\N}$ are $\mathcal{E}^{\theta_\delta}$-Cauchy and hence are also $\mathcal{E}^{\theta}$-Cauchy. Since
\begin{align*}
    |\mathcal{E}^\theta(\varphi_k,\psi_k)-\mathcal{E}^\theta(\varphi_\ell,\psi_\ell)|\leq | \mathcal{E}^\theta(\varphi_k,\psi_k-\psi_\ell)|+| \mathcal{E}^\theta(\varphi_k-\varphi_\ell,\psi_\ell)|,
\end{align*}
$\left(\mathcal{E}^\theta(\varphi_k,\psi_k)\right)_{k\in\N}$ is a Cauchy sequence. We then define
\begin{align*}
    \mathcal{E}^\theta(u,v):=\lim_{k\to\infty}\mathcal{E}^\theta(\varphi_k,\psi_k).
\end{align*}
This definition does not depend on the choice of $(\varphi_k)_{k\in\N}$ and $(\psi_k)_{k\in\N}$. Therefore, $\mathcal{E}^\theta$ is extended to be a symmetric form on $L^2(\Omega)$ with its domain $\mathcal{F}^{\theta_\delta}$. We will still denote this specific extension by $\mathcal{E}^\theta$. By the above definition, we also have
\begin{align}\label{Equation:Etheta bounded by Etheta_delta}
    \mathcal{E}^\theta(u,u)\leq \mathcal{E}^{\theta_\delta} (u,u),\quad \forall u\in\mathcal{F}^{\theta_\delta}.
\end{align}

The extension $\mathcal{E}^\theta$ constructed above has the following compatibility property. On the one hand, we have $\mathcal{F}^{\theta_{\delta'}}\subset\mathcal{F}^{\theta_\delta}$ for every $0<\delta<\delta'$, and the definitions of $\mathcal{E}^\theta(u,v)$ using $\delta$ and $\delta'$ agree whenever $u,v\in\mathcal{F}^{\theta_{\delta'}}$.
On the other hand, assume that $(\theta,V)$ is an admissible pair and let $(\mathcal{F}^\theta,\mathcal{G}^\theta)$ be the Dirichlet space for it. Then for every $\delta>0$, we have $\mathcal{F}^{\theta_{\delta}}\subset\mathcal{F}^\theta$ and
\begin{align*}
    \mathcal{E}^{\theta}(u,v)= \mathcal{G}^{\theta} (u,v),\quad \forall u,v\in\mathcal{F}^{\theta_\delta}.
\end{align*}

Now we fix an arbitrary $\delta_0>0$. By  Theorem \ref{bounded msh in F loc} (i) and the above discussion,  $\mathcal{F}^{\theta_{\delta_0}}$ contains the linear span of 
$$\{\eta u: \eta\in \mathscr{C}_c^\infty(\Omega,\R), u\in SH_m(\Omega)\cap L_{loc}^\infty(\Omega) \},$$
$\eta u^{[\delta]}$ converges to $\eta u$ in $\mathcal{E}^{\theta}$ as $\delta \downarrow 0$ for every $\eta\in\mathscr{C}_c^\infty(\Omega,\R)$, and hence
\begin{align}\label{Equation: Etheta convergence (for bounded functions) of approximation sequence}
    \mathcal{E}^\theta(\eta u,\xi v)=\lim_{\delta\to 0}\mathcal{E}^\theta(\eta u^{[\delta]},\xi v^{[\delta]})=\lim_{\delta\to 0}\lim_{\delta'\to 0}\mathcal{E}^\theta(\eta u^{[\delta]},\xi v^{[\delta']})
\end{align}
for every $u,v\in SH_m(\Omega)\cap L_{loc}^\infty(\Omega)$ and $\eta,\xi\in \mathscr{C}_c^\infty(\Omega,\R)$. Moreover, if we assume in addition that 
$\xi=1$ on an open neighborhood of $\supp \eta$, then
\begin{align}\label{Equation:formula of Etheta when xi equals 1 near support of eta}
    \mathcal{E}^\theta(\eta u,\xi v)
    =\lim_{\delta\to 0}\lim_{\delta'\to 0}\mathcal{E}^\theta(\eta u^{[\delta]},\xi v^{[\delta']})
    =-\int_\Omega \eta u dd^c v\wedge\theta.
\end{align}
By Theorem \ref{Theorem:Cauchy for locally uniformly bounded and monotone}, if $(u_k)_{k\in\N}\subset SH_m(\Omega)$ is monotone, locally uniformly bounded, and $\lim_{k\to\infty}u_k=u~V$-a.e. on $\Omega$ for some $u\in SH_m(\Omega)\cap L_{loc}^\infty(\Omega)$, then for every $\eta\in \mathscr{C}_c^\infty(\Omega,\R)$ we have
\begin{align}\label{equation:approximation of Etheta under monotone sequence}
    \lim_{k\to\infty}\mathcal{E}^\theta(\eta u_k- \eta u,\eta u_k- \eta u)=0.
\end{align}

In what follows, we use both of the following terminologies.
\begin{itemize}
    \item[(i)] Dirichlet space $(\mathcal{F}^\theta,\mathcal{E}^\theta)$ for an admissible pair $(\theta,\mu)$.
    \item[(ii)] Extended symmetric form $\mathcal{E}^\theta$ induced by $u_1,\ldots,u_{m-1}\in SH_m(\Omega)\cap L^\infty(\Omega)$.
\end{itemize}
These should be distinguished despite the shared notation.
\section{$\mu$-quasi-monotone convergence and complex Hessian operators}
In this section, we apply the Dirichlet form methods developed in the previous section to establish several results on complex Hessian operators. We also study the relationship between Dirichlet forms and the $\mu$-quasi-monotone convergence of $m$-subharmonic functions.

We first establish the commutativity of complex Hessian operators, inspired by \cite[Lemma 3.1]{FO87}.
\begin{lemma}\label{lemma:communtativity of complex hessian operators}
    For every $u_1,\dots,u_m\in SH_m(\Omega)\cap L^\infty_{loc}(\Omega)$ and every permutation $(j_1,j_2,\dots,j_m)$ of $(1,2,\dots,m)$, we have
    \begin{align*}
        dd^c u_1\wedge\dots\wedge dd^c u_m\wedge\omega^{n-m}=dd^c u_{j_1}\wedge\dots\wedge dd^c u_{j_m}\wedge\omega^{n-m}.
    \end{align*}
\end{lemma}
\begin{proof}
    Fix
    an arbitrary $\eta \in \mathscr{C}_c^{\infty}(\Omega,\R)$. We will prove by induction on $k=2,\dots,m$ that for every $u_1,\dots,u_m\in SH_m(\Omega)\cap L^\infty_{loc}(\Omega)$ and every permutation $(j_1,j_2,\dots,j_k)$ of $(1,2,\dots,k)$, we have
    \begin{align*}
        \int_{\Omega}\eta dd^c u_1\wedge\dots\wedge dd^c u_k\wedge \theta=\int_{\Omega}\eta dd^c u_{j_1}\wedge\dots\wedge dd^c u_{j_k}\wedge \theta,
    \end{align*}
    where $\theta:=dd^c u_{k+1}\wedge\dots\wedge dd^c u_m\wedge\omega^{n-m}$.

    We first prove the desired assertion for $k=2$. Set
    \begin{align*}
        a_{\ell,r}:=\int_\Omega \eta dd^c u_1^{[\frac{1}{\ell}]} \wedge dd^c u_2^{[\frac{1}{r}]} \wedge \theta
    \end{align*}
    for every $\ell,r\in\N$ large enough. Notice that
    \begin{align*}
        \lim_{\ell\to\infty}\lim_{r\to\infty}a_{\ell,r}  = \lim_{\ell\to\infty}\int_\Omega \eta dd^c u_1^{[\frac{1}{\ell}]} \wedge dd^c u_2 \wedge \theta  = \int_\Omega \eta dd^c u_1 \wedge dd^c u_2 \wedge \theta.
    \end{align*}
    Thus, it suffices to show that 
    \begin{align*}
        \lim_{\ell\to\infty}a_{\ell,\ell}=\lim_{\ell\to\infty}\lim_{r\to\infty}a_{\ell,r},
    \end{align*}
    for then
    \begin{align*}
        \int_\Omega \eta dd^c u_1 \wedge dd^c u_2 \wedge \theta 
        &=\lim_{\ell\to\infty}\int_\Omega \eta dd^c u_1^{[\frac{1}{\ell}]} \wedge dd^c u_2^{[\frac{1}{\ell}]} \wedge \theta\\
        &=\lim_{\ell\to\infty}\int_\Omega \eta dd^c u_2^{[\frac{1}{\ell}]} \wedge dd^c u_1^{[\frac{1}{\ell}]} \wedge \theta\\
        &=\int_\Omega \eta dd^c u_2 \wedge dd^c u_1 \wedge \theta.
    \end{align*}
    Let $\xi \in \mathscr{C}_c^{\infty}(\Omega,\R)$ so that $\xi = 1$ on an open neighborhood of $\supp\eta$, and $C_1>0$ be a constant so that $C_1 dd^c |z|^2\geq dd^c \eta$ on $\Omega$. Then $\rho := C_1\lvert z\rvert^2 - \eta \in SH_m(\Omega) \cap L^{\infty}(\Omega)$. Since both $\theta_1 := C_1dd^c\lvert z\rvert^2\wedge\theta$ and $\theta_2 := dd^c\rho \wedge \theta$ belong to $\mathcal{M}_m(\Omega)$, Lemma \ref{Lemma:crucial Cauchy sequence} implies that both $(\xi u_1^{[\frac{1}{\ell}]})_{\ell\in\N}$ and $(\xi u_2^{[\frac{1}{\ell}]})_{\ell\in\N}$ are $\mathcal{E}^{\theta_1}$-Cauchy and $\mathcal{E}^{\theta_2}$-Cauchy. Hence $\left(\mathcal{E}^{\theta_1}(\xi u_1^{[\frac{1}{\ell}]},\xi u_1^{[\frac{1}{\ell}]})\right)_{\ell\in\N}$ and $\left(\mathcal{E}^{\theta_2}(\xi u_1^{[\frac{1}{\ell}]},\xi u_1^{[\frac{1}{\ell}]})\right)_{\ell\in\N}$ are bounded from above by a constant $C_2>0$. For every $\varepsilon>0$, there exists $N\in\N$ such that for every $\ell,r\geq N$,
    \begin{align*}
        \mathcal{E}^{\theta_1}(\xi u_2^{[\frac{1}{\ell}]}-\xi u_2^{[\frac{1}{r}]},\xi u_2^{[\frac{1}{\ell}]}-\xi u_2^{[\frac{1}{r}]})<\varepsilon,\quad\mathcal{E}^{\theta_2}(\xi u_2^{[\frac{1}{\ell}]}-\xi u_2^{[\frac{1}{r}]},\xi u_2^{[\frac{1}{\ell}]}-\xi u_2^{[\frac{1}{r}]})<\varepsilon.
    \end{align*}
    Then we have
    \begin{align*}
        |a_{\ell,r}-a_{\ell,\ell}|
        &= \left\vert \int_{\Omega}\eta dd^c(\xi u_1^{[\frac{1}{\ell}]}) \wedge dd^c(\xi u_2^{[\frac{1}{r}]}-\xi u_2^{[\frac{1}{\ell}]})\wedge \theta \right\vert\\
        &= \left\vert \int_{\Omega}d(\xi u_1^{[\frac{1}{\ell}]}) \wedge d^c(\xi u_2^{[\frac{1}{r}]}-\xi u_2^{[\frac{1}{\ell}]}) \wedge dd^c\eta \wedge \theta \right\vert\\
        &= \left\vert \int_{\Omega}d(\xi u_1^{[\frac{1}{\ell}]}) \wedge d^c(\xi u_2^{[\frac{1}{r}]}-\xi u_2^{[\frac{1}{\ell}]}) \wedge (\theta_1-\theta_2) \right\vert\\
        &\leq \left\vert \mathcal{E}^{\theta_1}(\xi u_1^{[\frac{1}{\ell}]}, \xi u_2^{[\frac{1}{r}]}-\xi u_2^{[\frac{1}{\ell}]}) \right\vert + \left\vert \mathcal{E}^{\theta_2}(\xi u_1^{[\frac{1}{\ell}]}, \xi u_2^{[\frac{1}{r}]}-\xi u_2^{[\frac{1}{\ell}]}) \right\vert\\
        &< 2\sqrt{C_2\varepsilon}.
    \end{align*}
    Hence $|\lim_{r\to \infty}a_{\ell,r}-a_{\ell,\ell}|\leq 2 \sqrt{C_2\varepsilon}$ for every $\ell\geq\N$, which implies that $\lim_{\ell\to\infty}a_{\ell,\ell}$ exists and is equal to $\lim_{\ell\to\infty}\lim_{r\to\infty}a_{\ell,r}$.

    Assume that the desired assertion is proved for some $2\leq k\leq m-1$. For proving the case of $k+1$, we only need to show that
    \begin{align*}
        \int_{\Omega}\eta dd^c u_1\wedge\dots\wedge dd^c u_k\wedge dd^c u_{k+1}\wedge \theta=\int_{\Omega}\eta dd^c u_1\wedge\dots\wedge dd^c u_{k+1}\wedge dd^c u_k\wedge \theta.
    \end{align*}
    Indeed,
    \begin{align*}
        &\int_{\Omega}\eta dd^c u_1\wedge\dots\wedge dd^c u_k\wedge dd^c u_{k+1}\wedge \theta\\
        &=\lim_{\ell\to\infty} \int_{\Omega}\eta dd^c u_1^{[\frac{1}{\ell}]}\wedge dd^c u_2\wedge\dots\wedge dd^c u_k\wedge dd^c u_{k+1}\wedge \theta\\
        &=\lim_{\ell\to\infty} \int_{\Omega}\eta dd^c u_2\wedge\dots\wedge dd^c u_k\wedge dd^c u_{k+1}\wedge dd^c u_1^{[\frac{1}{\ell}]}\wedge \theta\\
        &=\lim_{\ell\to\infty} \int_{\Omega}\eta dd^c u_2\wedge\dots\wedge dd^c u_{k+1}\wedge dd^c u_k\wedge dd^c u_1^{[\frac{1}{\ell}]}\wedge \theta\\
        &=\lim_{\ell\to\infty} \int_{\Omega}\eta dd^c u_1^{[\frac{1}{\ell}]}\wedge dd^c u_2\wedge\dots\wedge dd^c u_{k+1}\wedge dd^c u_k\wedge  \theta\\
        &=\int_{\Omega}\eta dd^c u_1\wedge\dots\wedge dd^c u_{k+1}\wedge dd^c u_k\wedge \theta.
    \end{align*}
    This finishes the proof.
\end{proof}

Motivated by the computation in \cite{FO87}, we prove a lemma on $\mathcal{E}^{\theta}$-convergence that will be used in subsequent proofs.
\begin{lemma}\label{Etheta-convergence}
    For every $k\in\N$, let $\theta_k$ be a closed positive $(n-1,n-1)$-current on $\Omega$, $(\theta_k,\mu)$ be an admissible pair, and $(\mathcal{F}^{\theta_k}, \mathcal{E}^{\theta_k})$ be the Dirichlet space for it. Let $u_k\in \mathcal{F}^{\theta_k}_{loc}$ for every $k\in\N$, $u\in \bigcap_{k\in\N}\mathcal{F}^{\theta_k}_{loc}$, and $\eta\in \mathscr{C}_c^\infty(\Omega,\R)$. Assume that one of the following conditions is satisfied.
    \begin{itemize}
        \item[(i)] $\limsup_{k\to\infty} \left(\mathcal{E}^{\theta_k}(\eta u, \eta u)-\mathcal{E}^{\theta_k}(\eta u_k, \eta u_k)\right)\leq 0$,\\
        $\limsup_{k \rightarrow \infty}\mathcal{E}^{\theta_k}\left(\eta (u_k-u), \eta u_k\right) \leq 0$.
        \item[(ii)] $\limsup_{k \rightarrow \infty}\left(\mathcal{E}^{\theta_k}(\eta u_k, \eta u_k)-\mathcal{E}^{\theta_k}(\eta u, \eta u)\right) \leq 0$,\\
        $\limsup_{k \rightarrow \infty}\mathcal{E}^{\theta_k}\left(\eta (u-u_k), \eta u\right) \leq 0$.
    \end{itemize}
    Then $\lim_{k \to \infty}\mathcal{E}^{\theta_k}(\eta u_k-\eta u, \eta u_k-\eta u) = 0$.
\end{lemma}
\begin{proof}
    It follows from Proposition \ref{proposition:characterization for local domain of Dirichlet space} that $\eta u_k, \eta u \in \mathcal{F}^{\theta_k}$ for every $k\in\N$.
    
    Assume that condition (i) is satisfied. By direct computation,
    \begin{align*}
        \mathcal{E}^{\theta_k}(\eta (u_k-u), \eta (u_k-u)) = \mathcal{E}^{\theta_k}(\eta u, \eta u)-\mathcal{E}^{\theta_k}(\eta u_k, \eta u_k) + 2\mathcal{E}^{\theta_k}(\eta (u_k-u), \eta u_k).
    \end{align*}
    Then we have
    \begin{align*}
        0 &\leq \liminf_{k \rightarrow \infty}\mathcal{E}^{\theta_k}(\eta (u_k-u), \eta (u_k-u)) \\
        &\leq \limsup_{k \rightarrow \infty}\mathcal{E}^{\theta_k}(\eta (u_k-u), \eta (u_k-u)) \\
        &\leq \limsup_{k \rightarrow \infty}\left(\mathcal{E}^{\theta_k}(\eta u, \eta u) - \mathcal{E}^{\theta_k}(\eta u_k, \eta u_k)\right) + 2\limsup_{k \rightarrow \infty}\mathcal{E}^{\theta_k}(\eta (u_k-u), \eta u_k) \\
        &\leq 0.
    \end{align*}
    Therefore, $\lim_{k\to\infty} \mathcal{E}^{\theta_k}(\eta u_k-\eta u,\eta u_k-\eta u)=0$.
    
    Assume that condition (ii) is satisfied. The proof is similar to the one above.
\end{proof}

\begin{remark}\label{remark: Etheta-convergence for bounded functions}
    For every $k\in\N$, assume that $\theta_k=dd^c w_{k,1}\wedge\cdots\wedge dd^c w_{k,m-1}\wedge\omega^{n-m}$ for $w_{k,1},\dots,w_{k,m-1}\in SH_m(\Omega)\cap L^\infty(\Omega)$, and let $\mathcal{E}^{\theta_k}$ be the  extended symmetric form induced by $w_{k,1},\dots,w_{k,m-1}$. If $u_k, u\in SH_m(\Omega)\cap L^\infty_{loc}(\Omega)$, then Lemma \ref{Etheta-convergence} remains valid under the same condition.
\end{remark}

We now establish the complete form of the  monotonicity for Dirichlet forms.

\begin{theorem}[Monotonicity for Dirichlet forms]\label{Theorem: Monotonicity for admissible pairs}
    Let $(\mathcal{F}^\theta,\mathcal{E}^\theta)$ and $(\mathcal{F}^{\Theta},\mathcal{E}^{\Theta})$ be the Dirichlet spaces for admissible pairs $(\theta,\mu)$ and $(\Theta,\nu)$, respectively, where
    \begin{align*}
        \theta=dd^c w_1\wedge\dots\wedge dd^c w_{m-1}\wedge\omega^{n-m},\quad \Theta=dd^c h_1\wedge\dots\wedge dd^c h_{m-1}\wedge\omega^{n-m},
    \end{align*}
    with $w_1,\dots,w_{m-1}, h_1,\dots,h_{m-1}\in SH_m(\Omega)\cap L_{loc}^{\infty}(\Omega)$. Let $u_1,u_2,v_1,v_2\in SH_m(\Omega)\\\cap L^{\infty}_{loc}(\Omega)$ and $K\subset\Omega$ be a compact set. Assume that
    \begin{align*}
        &w_k\leq h_k~ \text{on} ~K, \quad w_k=h_k ~\text{on} ~\Omega\setminus K,\quad \forall k=1,\dots,m-1,\\
        &u_k\leq v_k~ \text{on} ~K, \quad u_k=v_k ~\text{on} ~\Omega\setminus K,\quad\forall k=1,2.
    \end{align*}
    Then for every $\xi,\gamma \in \mathscr{C}_c^{\infty}(\Omega,\R)$ such that $\xi = \gamma=1$ on an open neighborhood of $K$, we have
    \begin{align*}
        \mathcal{E}^{\theta}(\xi u_1, \gamma u_2) \geq \mathcal{E}^{\Theta}(\xi v_1, \gamma v_2).
    \end{align*}
\end{theorem}
\begin{proof}
    By Lemma \ref{Lemma: weak monotonicity for admissible pairs} we have $\mathcal{E}^{\Theta}(\xi u_1, \gamma u_2)\geq \mathcal{E}^{\Theta}(\xi v_1, \gamma v_2)$, hence we only need to prove $\mathcal{E}^{\theta}(\xi u_1, \gamma u_2) \geq \mathcal{E}^{\Theta}(\xi u_1, \gamma u_2)$. Then, by Theorem \ref{bounded msh in F loc} (i), it suffices to prove that $\mathcal{E}^{\theta}(\xi u_1^{[\delta]}, \gamma u_2^{[\delta]}) \geq \mathcal{E}^{\Theta}(\xi u_1^{[\delta]}, \gamma u_2^{[\delta]})$
    for every $\delta>0$ small enough, i.e.,
    \begin{align}\label{equation: goal in proof of monotonicity}
        \int_{\Omega}\xi u_1^{[\delta]}dd^c(\gamma u_2^{[\delta]}) \wedge (\Theta-\theta)\geq0.
    \end{align}
    Let $E\Subset\Omega$ be an open neighborhood of $K$ such that $\xi=\gamma=1$ on $E$. Then
    \begin{align*}
        &\int_{\Omega}\xi u_1^{[\delta]}dd^c(\gamma u_2^{[\delta]}) \wedge dd^c (h_1-w_1) \wedge dd^c w_2 \wedge \cdots \wedge dd^c w_{m-1} \wedge \omega^{n-m} \\
        &= \int_{\Omega}(h_1-w_1) dd^c(\xi u_1^{[\delta]})\wedge dd^c(\gamma u_2^{[\delta]}) \wedge dd^c w_2 \wedge \cdots \wedge dd^c w_{m-1} \wedge \omega^{n-m} \\
        &= \int_E(h_1-w_1) dd^c u_1^{[\delta]}\wedge dd^c u_2^{[\delta]}\wedge dd^c w_2 \wedge \cdots \wedge dd^c w_{m-1} \wedge \omega^{n-m} \\
        &\quad + \int_{\Omega\setminus E}(h_1-w_1) dd^c(\xi u_1^{[\delta]})\wedge dd^c(\gamma u_2^{[\delta]}) \wedge dd^c w_2 \wedge \cdots \wedge dd^c w_{m-1} \wedge \omega^{n-m}\\
        &\geq 0.
    \end{align*}
    Applying Lemma \ref{lemma:communtativity of complex hessian operators} and repeating the above process $m-2$ times yields (\ref{equation: goal in proof of monotonicity}).
\end{proof}

Using the same proofs of Lemma \ref{Lemma: weak monotonicity for admissible pairs} and Theorem \ref{Theorem: Monotonicity for admissible pairs}, we can establish the monotonicity property for the extended symmetric forms induced by bounded $m$-subharmonic functions. We omit the details of the proof.
\begin{theorem}[Monotonicity for extended symmetric forms]\label{Theorem: Monotonicity for bounded functions}
     Let $w_1,\dots,w_{m-1},\\ h_1,\dots,h_{m-1}\in SH_m(\Omega)\cap L^{\infty}(\Omega)$, and denote
    \begin{align*}
        \theta:=dd^c w_1\wedge\dots\wedge dd^c w_{m-1}\wedge\omega^{n-m},\quad \Theta:=dd^c h_1\wedge\dots\wedge dd^c h_{m-1}\wedge\omega^{n-m}.
    \end{align*}
    Let $\mathcal{E}^{\theta}$ and $\mathcal{E}^\Theta$ be the extended symmetric forms induced by $w_1,\dots,w_{m-1}$ and $h_1,\dots,h_{m-1}$, respectively.
    Let $u_1,u_2,v_1,v_2\in SH_m(\Omega)\cap L^{\infty}_{loc}(\Omega)$ and $K\subset\Omega$ be a compact set. Assume that
    \begin{align*}
        &w_k\leq h_k~ \text{on} ~K, \quad w_k=h_k ~\text{on} ~\Omega\setminus K,\quad \forall k=1,\dots,m-1,\\
        &u_k\leq v_k~ \text{on} ~K, \quad u_k=v_k ~\text{on} ~\Omega\setminus K,\quad \forall k=1,2.
    \end{align*}
    Then for every $\xi,\gamma \in \mathscr{C}_c^{\infty}(\Omega,\R)$ such that $\xi = \gamma=1$ on an open neighborhood of $K$, we have
    \begin{align*}
        \mathcal{E}^{\theta}(\xi u_1, \gamma u_2) \geq \mathcal{E}^{\Theta}(\xi v_1, \gamma v_2).
    \end{align*}
\end{theorem}

We next prove the weak convergence of complex Hessian operators under monotone sequences of $m$-subharmonic functions, generalizing \cite[Theorem 3.2]{FO87}.
\begin{theorem}\label{Theorem:convergence of complex Hessian operatiors under monotone sequences}
    Assume that $u_{j,k},u_j\in SH_m(\Omega)$ for $j=0,\dots,m$ and $k\in\N$ are locally uniformly bounded and satisfy the following conditions.
    \begin{itemize}
        \item[(i)] Either $\lim_{k\to\infty}u_{0,k}=u_0$ on $\Omega$, or $(u_{0,k})_{k\in\N}$ is increasing with $\lim_{k\to\infty}u_{0,k}\\=u_0$ $V$-a.e. on $\Omega$.
        \item[(ii)] Either $(u_{j,k})_{k\in\N}$ is decreasing for every $j=1,\dots,m$ with $\lim_{k\to\infty}u_{j,k}=u_j$ on $\Omega$, or $(u_{j,k})_{k\in\N}$ is increasing for every $j=1,\dots,m$ with $\lim_{k\to\infty}u_{j,k}=u_j$ $V$-a.e. on $\Omega$.
    \end{itemize}
    Then we have 
    \begin{align*}
        u_{0,k} dd^c u_{1,k}\wedge\dots\wedge dd^c u_{m,k}\wedge\omega^{n-m}\to u_0 dd^c u_1\wedge\dots\wedge dd^c u_m\wedge\omega^{n-m},\quad k\to\infty
    \end{align*}
    in the sense of currents.
\end{theorem}
\begin{proof}
    Since the desired assertion is local, we can assume that $\Omega$ is a bounded strictly $m$-pseudoconvex domain, $-B_2<u_{0,k},u_0<B_2$, $-B_2<u_{j,k},u_j<-B_1$, for every $j=1,\dots,m$ and $k\in\N$, where $B_2>B_1>0$ are two constants.
    
    For proving the desired assertion, we only need to show
    \begin{align}\label{Equation:first desired inequality in convergence thm for monotone sequence}
        \lim_{k\to\infty} \int_\Omega \eta (u_{0,k} - u_0) dd^c u_1\wedge\dots\wedge dd^c u_m\wedge\omega^{n-m}=0
    \end{align}
    and
    \begin{align}\label{Equation:second desired inequality in convergence thm for monotone sequence}
        \lim_{k\to\infty} \int_\Omega& \eta u_{0,k} dd^c u_{1,k}\wedge\dots\wedge dd^c u_{j-1,k} \wedge dd^c (u_{j,k}-u_j) \\
        &\quad\wedge dd^c u_{j+1} \wedge\cdots\wedge dd^c u_m\wedge\omega^{n-m}=0\notag
    \end{align}
    for every $\eta\in \mathscr{C}_c^\infty(\Omega,\R)$ and $j=1,\dots,m$.

    Under the former condition of (i), (\ref{Equation:first desired inequality in convergence thm for monotone sequence}) is deduced by the dominated convergence theorem. Under the latter condition of (i), denote $\theta:=dd^c u_2\wedge\dots\wedge dd^c u_m\wedge\omega^{n-m}$. Let $\mathcal{E}^\theta$ be the extended symmetric form induced by $u_2,\dots,u_m$. Let $\xi\in \mathscr{C}^\infty_c(\Omega,\R)$ so that $\xi=1$ on an open neighborhood of $\supp\eta$. By (\ref{Equation:formula of Etheta when xi equals 1 near support of eta}) and (\ref{equation:approximation of Etheta under monotone sequence}), we have
    \begin{align*}
        \lim_{k \rightarrow \infty}\int_\Omega \eta (u_{0,k} - u_0) dd^c u_1\wedge\cdots\wedge dd^c u_m\wedge\omega^{n-m}=\lim_{k \rightarrow \infty}\mathcal{E}^\theta(\eta(u_0 - u_{0,k}),\xi u_1)=0.
    \end{align*}
    
    We now prove (\ref{Equation:second desired inequality in convergence thm for monotone sequence}). Fix an arbitrary $\eta\in \mathscr{C}_c^\infty(\Omega,\R)$. By Proposition \ref{Proposition:localization principle}, there exist $\hat{u}_{j,k}, \hat{u}_j\in SH_m(\Omega)\cap L^\infty(\Omega)$ for $j=1,\dots,m$ and $k\in\N$, an open set $U$, and a compact set $K$ satisfying the following properties.
    \begin{itemize}
        \item[(i')] $\supp{\eta} \subset U \subset K \subset \Omega$, $(\hat{u}_{j,k})_{k\in\N}$ is uniformly bounded for every $j$.
        \item[(ii')] $\hat{u}_j = u_j$ on $U$, $\hat{u}_{j, k} = u_{j, k}$ on $U$, and $\hat{u}_j = \hat{u}_{j, k}$ on $\Omega \setminus K$ for every $j$ and $k$.
        \item[(iii')] Either $(\hat{u}_{j,k})_{k\in\N}$ is decreasing for every $j$ with $\lim_{k\to\infty}\hat{u}_{j,k}=\hat{u}_j$ on $\Omega$, or $(\hat{u}_{j,k})_{k\in\N}$ is increasing for every $j$ with $\lim_{k\to\infty}\hat{u}_{j,k}=\hat{u}_j$ $V$-a.e. on $\Omega$.
    \end{itemize}

    Let $\xi\in \mathscr{C}^\infty_c(\Omega,\R)$ so that $\xi=1$ on an open neighborhood of $K$. Fix an arbitrary $j\in\{1,\dots,m\}$. For every $k\in\N$, set
    \begin{align*}
        \theta_{j,k}:=dd^c \hat{u}_{1,k}\wedge\dots\wedge dd^c \hat{u}_{j-1,k}\wedge dd^c \hat{u}_{j+1}\wedge\dots\wedge dd^c \hat{u}_m\wedge\omega^{n-m},
    \end{align*}
    let $\mathcal{E}^{\theta_{j,k}}$ be the extended symmetric form induced by $\hat{u}_{1,k},\dots, \hat{u}_{j-1,k},\hat{u}_{j+1},\dots,\hat{u}_m$.
    By Lemma \ref{lemma:communtativity of complex hessian operators}, (\ref{Equation:formula of Etheta when xi equals 1 near support of eta}), and (\ref{equation:approximation of Etheta under monotone sequence}), the integral in equation (\ref{Equation:second desired inequality in convergence thm for monotone sequence}) is equal to
    \begin{align*}
        -\mathcal{E}^{\theta_{j,k}}(\eta u_{0,k}, \xi \hat{u}_{j,k} - \xi \hat{u}_j)=-\lim_{\ell\to\infty}\mathcal{E}^{\theta_{j,k}}(\eta u_{0,k}, \xi \hat{u}_{j,k} - \xi \hat{u}_{j,\ell}).
    \end{align*}
    Hence, we only need to show $$\lim_{k\to\infty}\lim_{\ell\to\infty}\mathcal{E}^{\theta_{j,k}}(\eta u_{0,k}, \xi \hat{u}_{j,k} - \xi \hat{u}_{j,\ell})=0,$$
    which can be implied by
    \begin{align*}
        \sup_{k\in\N}\mathcal{E}^{\theta_{j,k}}(\eta u_{0,k},\eta u_{0,k})<\infty
    \end{align*}
    and
    \begin{align*}
        \lim_{k,\ell\to\infty}\mathcal{E}^{\theta_{j,k}}(\xi \hat{u}_{j,k} - \xi \hat{u}_{j,\ell},\xi \hat{u}_{j,k} - \xi \hat{u}_{j,\ell})=0.
    \end{align*}
    
   By (\ref{Equation: Etheta convergence (for bounded functions) of approximation sequence}), we have $\mathcal{E}^{\theta_{j,k}}(\eta u_{0,k},\eta u_{0,k})=\lim_{\delta\to 0}\mathcal{E}^{\theta_{j,k}}(\eta u^{[\delta]}_{0,k},\eta u^{[\delta]}_{0,k})$. Then
    \begin{align*}
        &\mathcal{E}^{\theta_{j,k}}(\eta u^{[\delta]}_{0,k},\eta u^{[\delta]}_{0,k})\\
        &\leq 2 \int_\Omega \eta^2 d u^{[\delta]}_{0,k}\wedge d^c u^{[\delta]}_{0,k}\wedge \theta_{j,k}+ 2 \int_\Omega (u^{[\delta]}_{0,k})^2 d\eta\wedge d^c \eta\wedge \theta_{j,k}\\
        &\leq \int_\Omega \eta^2 dd^c (u^{[\delta]}_{0,k}+B_2)^2\wedge\theta_{j,k}+2 B_2^2 \int_\Omega d\eta\wedge d^c \eta\wedge \theta_{j,k}\\
        &= \int_\Omega (u^{[\delta]}_{0,k}+B_2)^2 dd^c (\eta^2)\wedge\theta_{j,k}+2 B_2^2  \int_\Omega d\eta\wedge d^c \eta\wedge \theta_{j,k}\\
        &\leq 6B_2^2 B_3\int_{\supp\eta} \theta_{j,k}\wedge dd^c|z|^2,
    \end{align*}
    where $B_3>0$ is a constant depending only on $n$ and $\eta$. By the Chern--Levine--Nirenberg inequality,
    \begin{align*}
        \sup_{k\in\N}\mathcal{E}^{\theta_{j,k}}(\eta u_{0,k},\eta u_{0,k})\leq6B_2^2 B_3\sup_{k\in\N}\int_{\supp\eta} \theta_{j,k}\wedge dd^c|z|^2\leq 6B_2^2 B_3 B_4<\infty,
    \end{align*}
    where $B_4 > 0$ is a constant independent of $k\in\N$.

    We finally prove $\lim_{k,\ell\to\infty}\mathcal{E}^{\theta_{j,k}}(\xi \hat{u}_{j,k} - \xi \hat{u}_{j,\ell},\xi \hat{u}_{j,k} - \xi \hat{u}_{j,\ell})=0$. Under the former property of (iii'), denote $\theta_j:=dd^c \hat{u}_1\wedge\dots\wedge dd^c \hat{u}_{j-1}\wedge dd^c \hat{u}_{j+1}\wedge\dots\wedge dd^c \hat{u}_m\wedge\omega^{n-m}$, and let $\mathcal{E}^{\theta_j}$ be the extended symmetric form induced by $\hat{u}_1,\dots,\hat{u}_{j-1},\hat{u}_{j+1},\dots, \hat{u}_m$. When $k<\ell$, by Theorem \ref{Theorem: Monotonicity for bounded functions} we have
    \begin{align*}
        & \mathcal{E}^{\theta_{j,k}}(\xi \hat{u}_{j,k},\xi \hat{u}_{j,k})\leq \mathcal{E}^{\theta_{j,k}}(\xi \hat{u}_{j,k},\xi \hat{u}_{j,\ell}),\\ &\mathcal{E}^{\theta_{j,k}}(\xi \hat{u}_{j,\ell},\xi \hat{u}_{j,\ell})\leq \mathcal{E}^{\theta_{j,\ell}}(\xi \hat{u}_{j,\ell},\xi \hat{u}_{j,\ell}),\\
        &\mathcal{E}^{\theta_{j,k}}(\xi \hat{u}_{j,k},\xi \hat{u}_{j,k})\leq \mathcal{E}^{\theta_j}(\xi \hat{u}_{j,k},\xi \hat{u}_{j,k})\leq \mathcal{E}^{\theta_j}(\xi \hat{u}_j,\xi \hat{u}_j).
    \end{align*}
    Hence, when $k<\ell$ we have
    \begin{align*}
        0&\leq \mathcal{E}^{\theta_{j,k}}(\xi \hat{u}_{j,k} - \xi \hat{u}_{j,\ell},\xi \hat{u}_{j,k} - \xi \hat{u}_{j,\ell}) \\
        &\leq \mathcal{E}^{\theta_{j,k}}(\xi \hat{u}_{j,\ell},\xi \hat{u}_{j,\ell})-\mathcal{E}^{\theta_{j,k}}(\xi \hat{u}_{j,k},\xi \hat{u}_{j,k})\\
        &\leq \mathcal{E}^{\theta_{j,\ell}}(\xi \hat{u}_{j,\ell},\xi \hat{u}_{j,\ell})-\mathcal{E}^{\theta_{j,k}}(\xi \hat{u}_{j,k},\xi \hat{u}_{j,k}),
    \end{align*}
    which implies that $\left(\mathcal{E}^{\theta_{j,k}}(\xi \hat{u}_{j,k},\xi \hat{u}_{j,k})\right)_{k\in\N}$ is increasing and bounded from above by $\mathcal{E}^{\theta_j}(\xi \hat{u}_j,\xi \hat{u}_j)$. Hence, $\left(\mathcal{E}^{\theta_{j,k}}(\xi \hat{u}_{j,k},\xi \hat{u}_{j,k})\right)_{k\in\N}$ has a finite limit. When $k>\ell$ we compute that
    \begin{align*}
        0&\leq \mathcal{E}^{\theta_{j,k}}(\xi \hat{u}_{j,k} - \xi \hat{u}_{j,\ell},\xi \hat{u}_{j,k} - \xi \hat{u}_{j,\ell}) \\
        &\leq \mathcal{E}^{\theta_{j,k}}(\xi \hat{u}_{j,k},\xi \hat{u}_{j,k})-\mathcal{E}^{\theta_{j,k}}(\xi \hat{u}_{j,\ell},\xi \hat{u}_{j,\ell})\\
        &\leq \mathcal{E}^{\theta_{j,k}}(\xi \hat{u}_{j,k},\xi \hat{u}_{j,k})-\mathcal{E}^{\theta_{j,\ell}}(\xi \hat{u}_{j,\ell},\xi \hat{u}_{j,\ell}).
    \end{align*}
    Therefore, $\lim_{k,\ell\to\infty}\mathcal{E}^{\theta_{j,k}}(\xi \hat{u}_{j,k} - \xi \hat{u}_{j,\ell},\xi \hat{u}_{j,k} - \xi \hat{u}_{j,\ell})=0$. Under the latter property of (iii'), when $k<\ell$, by Theorem \ref{Theorem: Monotonicity for bounded functions} we have
    \begin{align*}
        &\mathcal{E}^{\theta_{j,k}}(\xi \hat{u}_{j,k},\xi \hat{u}_{j,\ell})\geq \mathcal{E}^{\theta_{j,k}}(\xi \hat{u}_{j,\ell},\xi \hat{u}_{j,\ell})\geq \mathcal{E}^{\theta_{j,\ell}}(\xi \hat{u}_{j,\ell},\xi \hat{u}_{j,\ell}).
    \end{align*}
    Hence, when $k<\ell$ we have
    \begin{align*}
        0&\leq \mathcal{E}^{\theta_{j,k}}(\xi \hat{u}_{j,k} - \xi \hat{u}_{j,\ell},\xi \hat{u}_{j,k} - \xi \hat{u}_{j,\ell})\\&\leq \mathcal{E}^{\theta_{j,k}}(\xi \hat{u}_{j,k},\xi \hat{u}_{j,k})-\mathcal{E}^{\theta_{j,k}}(\xi \hat{u}_{j,\ell},\xi \hat{u}_{j,\ell})\\
        &\leq\mathcal{E}^{\theta_{j,k}}(\xi \hat{u}_{j,k},\xi \hat{u}_{j,k})-\mathcal{E}^{\theta_{j,\ell}}(\xi \hat{u}_{j,\ell},\xi \hat{u}_{j,\ell}),
    \end{align*}
    which implies that $\left(\mathcal{E}^{\theta_{j,k}}(\xi \hat{u}_{j,k},\xi \hat{u}_{j,k})\right)_{k\in\N}$ is decreasing and hence has a finite limit. When $k>\ell$ we have
    \begin{align*}
        0&\leq \mathcal{E}^{\theta_{j,k}}(\xi \hat{u}_{j,k} - \xi \hat{u}_{j,\ell},\xi \hat{u}_{j,k} - \xi \hat{u}_{j,\ell})\\
        &\leq \mathcal{E}^{\theta_{j,k}}(\xi \hat{u}_{j,\ell},\xi \hat{u}_{j,\ell})- \mathcal{E}^{\theta_{j,k}}(\xi \hat{u}_{j,k},\xi \hat{u}_{j,k})\\
        &\leq \mathcal{E}^{\theta_{j,\ell}}(\xi \hat{u}_{j,\ell},\xi \hat{u}_{j,\ell})- \mathcal{E}^{\theta_{j,k}}(\xi \hat{u}_{j,k},\xi \hat{u}_{j,k}).
    \end{align*}
    Therefore, $\lim_{k,\ell\to\infty}\mathcal{E}^{\theta_{j,k}}(\xi \hat{u}_{j,k} - \xi \hat{u}_{j,\ell},\xi \hat{u}_{j,k} - \xi \hat{u}_{j,\ell})=0$. 
\end{proof}

We now show that $\mu$-quasi-monotone convergence implies $\mathcal{E}_1^{\theta}$-convergence, generalizing Theorem \ref{Theorem:Cauchy for locally uniformly bounded and monotone}. This result illustrates the compatibility of $\mu$-quasi-monotone convergence with our Dirichlet form approach. 
\begin{theorem}\label{Theorem:Etheta convergence under quasi-monotone convergence}
    Let $\theta\in\mathcal{M}_m(\Omega)$, and $(\mc F^\theta,\mc E^\theta)$ be the Dirichlet space for an admissible pair $(\theta,\mu)$. Assume that $(u_k)_{k\in\N}\subset SH_m(\Omega)$ is locally uniformly bounded and converges $\mu$-quasi-monotonically to $u\in SH_m(\Omega)\cap L_{loc}^\infty(\Omega)$. Then $\eta u_k$ converges to $\eta u$ in $\mathcal{E}_1^{\theta}$ for every $\eta \in \mathscr{C}_c^{\infty}(\Omega,\R)$.
\end{theorem}

\begin{proof}
    Since $u_k$ converges to $u$ $\mu$-a.e. on $\Omega$, $\eta u_k$ converges to $\eta u$ in $L^2(\Omega;\mu)$ for every $\eta\in \mathscr{C}_c^\infty(\Omega,\R)$. We then only need to prove that $\eta u_k$ converges to $\eta u$ in $\mathcal{E}^\theta$. Combining Remark \ref{Remark: mu-quasi-monotone convergence for locally bounded functions} with the same argument as the first paragraph in the proof of Theorem \ref{Theorem:Cauchy for locally uniformly bounded and monotone}, we can assume that $\Omega$ is a bounded strictly $m$-pseudoconvex domain, $(u_k)_{k\in\N},u$ are uniformly bounded on $\Omega$, and there exist uniformly bounded sequences $(\varphi_k)_{k\in\N},(\psi_k)_{k\in\N}\subset SH_m(\Omega)$ satisfying the following properties.
        \begin{itemize}
            \item [(i)] $\varphi_k \leq u_k \leq \psi_k$ for every $k\in\N$.
            \item [(ii)] $(\varphi_k)_{k \in \N}$ is increasing and $(\psi_k)_{k \in \N}$ is decreasing.
            \item [(iii)] $\varphi_k$ converges to $u$ $\mu$-a.e. on $\Omega$ and $\psi_k$ converges to $u$ $\mu$-a.e. on $\Omega$.
        \end{itemize}
    Moreover, we can assume that $u$, $(u_k)_{k\in\N}$, $(\varphi_k)_{k\in\N}$, and $(\psi_k)_{k\in\N}$ have a uniform negative upper bound. Fix an arbitrary $\eta\in \mathscr{C}_c^{\infty}(\Omega,\R)$. By Lemma \ref{Lemma: reduction process for admissible pairs}, there exist uniformly bounded functions $\hat{u} \in SH_m(\Omega)$ and  $(\hat{u}_k)_{k \in \N}, (\hat{\varphi}_k)_{k \in \N},(\hat{\psi}_k)_{k \in \N}  \subset SH_m(\Omega)$, an open set $U$, and a compact set $K$ that satisfy the following properties.
    \begin{itemize}
        \item[(i')] $\supp\eta \subset U \subset K \subset \Omega$.
        \item[(ii')] On $U$ and for every $k \in \N$, we have $\hat{u} = u$, $\hat{u}_k = u_k$, $\hat{\varphi}_k = \varphi_k$, and $\hat{\psi}_k = \psi_k$.
        \item[(iii')] On $\Omega \setminus K$ and for every $k\in \N$, we have $\hat{u}= \hat{u}_k=\hat{\varphi}_k=\hat{\psi}_k$.
        \item[(iv')] $\hat{\varphi}_k \leq\hat{u}_k \leq \hat{\psi}_k$ for every $k\in\N$.
        \item [(v')] $(\hat{\varphi}_k)_{k \in \N}$ is increasing and $(\hat{\psi}_k)_{k \in \N}$ is decreasing.
        \item [(vi')] $\hat{\varphi}_k$ converges to $\hat{u}$ $\mu$-a.e. on $\Omega$ and $\hat{\psi}_k$ converges to $\hat{u}$ $\mu$-a.e. on $\Omega$.
    \end{itemize}

    Let $\xi \in \mathscr{C}_c^{\infty}(\Omega, \R)$ so that $\xi = 1$ on an open neighborhood of $K$. By Lemma \ref{Lemma: reduction process for admissible pairs}, we only need to prove that $\xi\hat{u}_k$ converges to $\xi\hat{u}$ in $\mathcal{E}^\theta$. Then by Lemma \ref{Etheta-convergence}, we only need to prove
    \begin{align*}
        \lim_{k \rightarrow \infty}\mathcal{E}^{\theta}(\xi \hat{u}_k, \xi \hat{u}_k) = \mathcal{E}^{\theta}(\xi \hat{u}, \xi \hat{u}),\quad
        \lim_{k \rightarrow \infty}\mathcal{E}^{\theta}(\xi (\hat{u}_k - \hat{u}), \xi \hat{u}) = 0.
    \end{align*}

    By Theorem \ref{Theorem: Monotonicity for admissible pairs},
    \[\mathcal{E}^\theta(\xi \hat\psi_k, \xi \hat\psi_k)\leq \mathcal{E}^\theta(\xi \hat u_k, \xi \hat u_k) \leq \mathcal{E}^\theta(\xi \hat\varphi_k, \xi \hat\varphi_k).\]
    By Theorem \ref{Theorem:Cauchy for locally uniformly bounded and monotone},
    \[\lim_{k\to\infty}\mc E^\theta(\xi\hat\varphi_k,\xi\hat\varphi_k)=\mc E^\theta(\xi\hat u,\xi\hat u),\quad\lim_{k\to\infty}\mc E^\theta(\xi\hat\psi_k,\xi\hat\psi_k)=\mc E^\theta(\xi\hat u,\xi\hat u).\]
    Therefore, 
    \[\lim_{k\to\infty}\mc E^\theta(\xi\hat u_k,\xi\hat u_k)=\mc E^\theta(\xi\hat u,\xi\hat u).\]
    By Theorem \ref{Theorem: Monotonicity for admissible pairs} again,
    \[\mc E^\theta(\xi(\hat{\psi}_k-\hat u),\xi\hat{u})\leq \mc E^\theta(\xi(\hat u_k-\hat u),\xi\hat u)\leq \mc E^\theta(\xi(\hat{\varphi}_k-\hat u),\xi\hat{u}).\]
    By Theorem \ref{Theorem:Cauchy for locally uniformly bounded and monotone} again,
    \begin{align*}
        \lim_{k\to\infty}\mc E^\theta(\xi(\hat{\psi}_k-\hat u),\xi\hat{u})=\lim_{k\to\infty} \mc E^\theta(\xi(\hat{\varphi}_k-\hat u),\xi\hat{u})=0.
    \end{align*}
    Therefore,
    \begin{align*}
        \lim_{k \rightarrow \infty}\mathcal{E}^{\theta}(\xi (\hat{u}_k - \hat{u}), \xi \hat{u}) = 0.
    \end{align*}
    This finishes the proof.
\end{proof}

\begin{remark}\label{remark: Etheta convergence for bounded functions and V-quasi-monotone sequence}
    Let $\mathcal{E}^\theta$ be the extended symmetric form induced by $w_1,\dots,w_{m-1}\in SH_m(\Omega)\cap L^\infty(\Omega)$. Assume that $(u_k)_{k\in\N}\subset SH_m(\Omega)$ is locally uniformly bounded and $V$-quasi-monotonically converges to $u\in SH_m(\Omega)\cap L_{loc}^\infty(\Omega)$. Then for every $\eta\in \mathscr{C}_c^\infty(\Omega)$,  Theorem \ref{Theorem:Etheta convergence under quasi-monotone convergence} and (\ref{Equation:Etheta bounded by Etheta_delta}) imply that
    \begin{align*}
    \lim_{k\to\infty}\mathcal{E}^\theta(\eta u_k- \eta u,\eta u_k- \eta u)=0.
    \end{align*}
\end{remark}

\begin{remark}\label{Remark: Etheat convergence and BT capacity}
    In Section 5 we will review the Bedford--Taylor capacity $\mathrm{Cap}_{BT, m}$ for $m$-subharmonic functions, and show in Proposition \ref{Proposition: equvialence of convergence quasi-monotonically and in BT capacity} that, for every locally uniformly bounded sequence of $m$-subharmonic functions, $V$-quasi-monotone convergence and convergence in $\mathrm{Cap}_{BT, m}$ are equivalent up to a subsequence. Let $\theta\in\mathcal{M}_m(\Omega)$, and  $(\mc F^\theta,\mc E^\theta)$ be the Dirichlet space for an admissible pair $(\theta, V)$.  Assume that $(u_k)_{k\in\N}\subset SH_m(\Omega)$ is locally uniformly bounded and converges to $u\in SH_m(\Omega)\cap L_{loc}^\infty(\Omega)$ in $\mathrm{Cap}_{BT, m}$. Theorem \ref{Theorem:Etheta convergence under quasi-monotone convergence} and Proposition \ref{Proposition: equvialence of convergence quasi-monotonically and in BT capacity} imply that there exists a subsequence $(u_{j_k})_{k\in\N}$ of $(u_k)_{k\in\N}$ such that $\eta u_{j_k}$ converges to $\eta u$ in $\mathcal{E}^{\theta}_1$ for every $\eta \in \mathscr{C}_c^{\infty}(\Omega, \R)$. Without considering Proposition \ref{Proposition: equvialence of convergence quasi-monotonically and in BT capacity}, we can directly show that $\eta u_k$ converges to $\eta u$ in $\mathcal{E}_1^{\theta}$ as follows. 
    
    Since $u_k$ converges to $u$ in $\mathrm{Cap}_{BT,m}$, it also converges in $L_{loc}^1(\Omega)$. Moreover, the locally uniformly bounded condition implies that $\eta u_k$ converges to $\eta u$ in $L^2(\Omega)$. Set 
    $$\psi_k := \left(\sup_{\ell \geq k}u_{\ell}\right)^*$$
    for every $k\in\N$. Then $(\psi_k)_{k\in\N}$ is decreasing and converges to $u$. Since convergence in $\mathrm{Cap}_{BT, m}$ relative to $\Omega$ implies convergence in $\mathrm{Cap}_{BT, m}$ relative to every relatively compact open subset of $\Omega$, we can use the same argument and notation as in the proof of Theorem \ref{Theorem:Etheta convergence under quasi-monotone convergence}. Let $E\Subset\Omega$ be an open subset containing $K$, $\xi\in \mathscr{C}_c^\infty(\Omega,\R)$ so that $\xi=1$ on $E$, and $\gamma \in \mathscr{C}_c^{\infty}(E,\R)$ so that $\gamma = 1$ on $K$. Since we have already shrunk $\Omega
    $, there exist $v_1,\dots,v_{m-1}\in SH_m(\Omega)\cap L^\infty(\Omega)$ such that $$\theta=dd^c v_1\wedge\cdots\wedge dd^c v_{m-1}\wedge\omega^{n-m}.$$
    
    We only need to prove that $\xi\hat{u}_k$ converges to $\xi\hat{u}$ in $\mathcal{E}^\theta$. Moreover, by Lemma \ref{Etheta-convergence}, we only need to prove that
    \begin{align*}
        \liminf_{k \rightarrow \infty}\mathcal{E}^{\theta}(\xi \hat{u}_k, \xi \hat{u}_k) \geq \mathcal{E}^{\theta}(\xi \hat{u}, \xi \hat{u}),\quad
        \lim_{k \rightarrow \infty}\mathcal{E}^{\theta}(\xi (\hat{u}_k - \hat{u}), \xi \hat{u}_k) = 0.
    \end{align*}
    
    It follows from Theorem \ref{Theorem: Monotonicity for admissible pairs} and Theorem \ref{Theorem:Cauchy for locally uniformly bounded and monotone} that
    \begin{align*}
        \liminf_{k \rightarrow \infty}\mathcal{E}^{\theta}(\xi \hat u_k, \xi \hat{u}_k) \geq \lim_{k \rightarrow \infty} \mathcal{E}^{\theta}(\xi \hat{\psi}_k, \xi \hat{\psi}_k) = \mathcal{E}^{\theta}(\xi \hat{u}, \xi \hat{u}).
    \end{align*}
    Since $\hat{u}_k=\hat{u}$ on $\Omega\setminus K$, Theorem \ref{bounded msh in F loc} (iii) implies that
    \begin{align*}
        \mathcal{E}^{\theta}(\xi (\hat{u}_k-\hat{u}), \xi \hat{u}_k) &= \mathcal{E}^{\theta}(\gamma \xi (\hat{u}_k-\hat{u}), \xi \hat{u}_k) \\
        &= \mathcal{E}^{\theta}(\gamma\xi \hat{u}_k, \xi \hat{u}_k) - \mathcal{E}^{\theta}(\gamma\xi \hat{u}, \xi\hat
        u_k) \\
        &= -\int_\Omega\gamma \xi \hat{u}_k dd^c\hat{u}_k \wedge \theta  + \int_\Omega\gamma \xi \hat{u}dd^c\hat{u}_k \wedge \theta \\
        &= \int_\Omega\xi (\hat{u} - \hat{u}_k)dd^c\hat{u}_k \wedge \theta.
    \end{align*}
    For every $\delta>0$, by the Chern--Levine--Nirenberg inequality we have
    \begin{align*}
        &\int_{\supp\xi}|\hat{u}_k -\hat{u}|dd^c \hat{u}_k\wedge\theta\\
        &=\int_{\supp\xi\cap\{|\hat{u}_k -\hat{u}|>\delta\}}|\hat{u}_k -\hat{u}|dd^c \hat{u}_k\wedge\theta + \int_{\supp\xi\cap\{|\hat{u}_k -\hat{u}|\leq\delta\}}|\hat{u}_k -\hat{u}|dd^c \hat{u}_k\wedge\theta\\
        &\leq \int_{\supp\xi\cap\{|\hat{u}_k -\hat{u}|>\delta\}}|\hat{u}_k -\hat{u}|\left(dd^c (\hat{u}_k+v_1+\cdots +v_{m-1})\right)^m\wedge\omega^{n-m} \\
        &\quad + \delta\int_{\supp\xi\cap\{|\hat{u}_k -\hat{u}|\leq\delta\}}dd^c \hat{u}_k\wedge\theta
        \\ 
        &\leq C\left( \mathrm{Cap}_{BT,m}(\supp\xi\cap\{|\hat{u}_k -\hat{u}|>\delta\},\Omega)+\delta\right),
    \end{align*}
    where $C>0$ is a constant independent of $k$. Finally, $|\hat{u}_k - \hat{u}|\leq |u_k -u|$ implies that $\hat{u}_k$ converges to $\hat{u}$ in $\mathrm{Cap}_{BT, m}$, and hence we have
    \begin{align*}
        \lim_{k \rightarrow \infty}\int_{\supp\xi}|\hat{u}_k -\hat{u}|dd^c \hat{u}_k\wedge\theta = 0.
    \end{align*}
    This finishes the proof of the remark.
\end{remark}

We now establish the weak convergence of complex Hessian operators under $V$-quasi-monotone sequences of $m$-subharmonic functions, using Theorem \ref{Theorem: Monotonicity for bounded functions}, Theorem \ref{Theorem:convergence of complex Hessian operatiors under monotone sequences}, and Theorem \ref{Theorem:Etheta convergence under quasi-monotone convergence}.
\begin{theorem}\label{theorem:continuity of hessian operator}
    Assume that $u_{j,k},u_j\in SH_m(\Omega)$ for $j=0,\dots,m$ and $k\in\N$ are locally uniformly bounded and satisfy the following conditions. 
    \begin{itemize}
        \item[(i)] Either $\lim_{k\to\infty}u_{0,k}=u_0$ on $\Omega$, or $u_{0,k}$ $V$-quasi-monotonically converges to $u_0$  on $\Omega$.
        \item[(ii)] $u_{j,k}$ $V$-quasi-monotonically converges to $u_j$  on $\Omega$ for every $j=1,\dots,m$.
    \end{itemize}
    Then we have 
    \begin{align*}
        u_{0,k} dd^c u_{1,k}\wedge\dots\wedge dd^c u_{m,k}\wedge\omega^{n-m}\to u_0 dd^c u_1\wedge\dots\wedge dd^c u_m\wedge\omega^{n-m},\quad k\to\infty
    \end{align*}
    in the sense of currents.
\end{theorem}
\begin{proof}
    For proving the desired assertion, we only need to show
    \begin{align}\label{Equation:first desired inequality in convergence thm for quasi-monotone sequence}
        \lim_{k\to\infty} \int_\Omega \eta (u_{0,k} - u_0) dd^c u_1\wedge\dots\wedge dd^c u_m\wedge\omega^{n-m}=0
    \end{align}
    and    
    \begin{align}\label{Equation:second desired inequality in convergence thm under quasi-monotone convergence}
        \lim_{k\to\infty} \int_\Omega& \eta u_{0,k} dd^c u_{1,k}\wedge\dots\wedge dd^c u_{j-1,k} \wedge dd^c (u_{j,k}-u_j) \\
        &\quad\wedge dd^c u_{j+1} \wedge\cdots\wedge dd^c u_m\wedge\omega^{n-m}=0\notag
    \end{align}
    for every $\eta\in \mathscr{C}_c^\infty(\Omega,\R)$ and $j=1,\dots,m$.

    (\ref{Equation:first desired inequality in convergence thm for quasi-monotone sequence}) follows by the same argument as (\ref{Equation:first desired inequality in convergence thm for monotone sequence}), using Remark \ref{remark: Etheta convergence for bounded functions and V-quasi-monotone sequence}.

    We now prove (\ref{Equation:second desired inequality in convergence thm under quasi-monotone convergence}). The desired assertion is local. By Remark \ref{Remark: mu-quasi-monotone convergence for locally bounded functions}, we can assume that $\Omega$ is a bounded strictly $m$-pseudoconvex domain, $(u_{j,k})_{k\in\N}$ and $u_j$ are uniformly bounded on $\Omega$ for every $j=0,\dots,m$, and there exist uniformly bounded sequences
    $(\varphi_{j,k})_{k\in\N},(\psi_{j,k})_{k\in\N}\subset SH_m(\Omega)$ for $j=1,\dots,m$ satisfying the following properties.
    \begin{itemize}
        \item [(i)] $\varphi_{j,k} \leq u_{j,k} \leq \psi_{j,k}$ for every $j,k$.
        \item [(ii)] $(\varphi_{j,k})_{k \in \N}$ is increasing and $(\psi_{j,k})_{k \in \N}$ is decreasing for every $j$.
        \item [(iii)] $\varphi_{j,k}$ converges to $u_j$ $V$-a.e. on $\Omega$ and $\psi_{j,k}$ converges to $u_j$ on $\Omega$ for every $j$.
    \end{itemize}
    Moreover, we can assume that $u_{j,k}$, $u_{j}$, $\varphi_{j,k}$, and $\psi_{j,k}$ have a uniform negative upper bound for all $j=1,\dots,m$ and $k\in\N$.
    Fix an arbitrary $\eta\in \mathscr{C}_c^\infty(\Omega,\R)$. By Proposition \ref{Proposition:localization principle}, there exist uniformly bounded functions $\hat{u}_{j,k}, \hat{u}_j, \hat{\varphi}_{j,k},\hat{\psi}_{j,k}\in SH_m(\Omega)$ for $j=1,\dots,m$ and $k\in\N$, an open set $U$, and a compact set $K$ satisfying the following properties.
    \begin{itemize}
        \item[(i')] $\supp{\eta} \subset U \subset K \subset \Omega$.
        \item[(ii')] $\hat{u}_j = u_j$, $\hat{u}_{j,k} = u_{j,k}$, $\hat{\varphi}_{j,k} = \varphi_{j,k}$, and $\hat{\psi}_{j,k} = \psi_{j,k}$ on $U$ for every $j,k$.
        \item[(iii')] $\hat{u}_j= \hat{u}_{j,k}=\hat{\varphi}_{j,k}=\hat{\psi}_{j,k}$ on $\Omega \setminus K$ for every $j,k$.
        \item[(iv')] $\hat{\varphi}_{j,k} \leq\hat{u}_{j,k} \leq \hat{\psi}_{j,k}$ for every $j,k$.
        \item [(v')] $(\hat{\varphi}_{j,k})_{k \in \N}$ is increasing and $(\hat{\psi}_{j,k})_{k \in \N}$ is decreasing for every $j$.
        \item [(vi')] $\hat{\varphi}_{j,k}$ converges to $\hat{u}_j$ $V$-a.e. on $\Omega$ and $\hat{\psi}_{j,k}$ converges to $\hat{u}_j$ on $\Omega$ for every $j$.
    \end{itemize}
    Notice that $\hat{\varphi}_{j,k} \leq\hat{u}_{j} \leq \hat{\psi}_{j,k}$ for every $j,k$.
    Let $E\Subset\Omega$ be an open subset containing $K$, $\xi\in \mathscr{C}_c^\infty(\Omega,\R)$ so that $\xi=1$ on $E$, and $\gamma \in \mathscr{C}_c^{\infty}(E,\R)$ so that $\gamma = 1$ on $K$. Fix an arbitrary $j\in\{1,\dots,m\}$. Denote
    \begin{align*}
        \theta_j &:= dd^c\hat{u}_1 \wedge \cdots \wedge dd^c\hat{u}_{j-1} \wedge dd^c\hat{u}_{j+1} \wedge \cdots \wedge dd^c\hat{u}_{m} \wedge \omega^{n-m},\\
        \theta_{j,k}&:=dd^c \hat{u}_{1,k}\wedge\dots\wedge dd^c \hat{u}_{j-1,k}\wedge dd^c \hat{u}_{j+1}\wedge \cdots \wedge dd^c \hat{u}_m\wedge\omega^{n-m},\\
        \alpha_{j, k}& := dd^c\hat{\varphi}_{1,k} \wedge \cdots \wedge dd^c\hat{\varphi}_{j-1,k} \wedge dd^c\hat{u}_{j+1} \wedge \cdots \wedge dd^c\hat{u}_{m} \wedge \omega^{n-m}, \\
        \beta_{j, k} &:= dd^c \hat{\psi}_{1, k} \wedge \cdots \wedge dd^c\hat{\psi}_{j-1, k} \wedge dd^c\hat{u}_{j+1} \wedge \cdots \wedge dd^c\hat{u}_m \wedge \omega^{n-m},
    \end{align*}
    and let $\mathcal{E}^{\theta_j},\mathcal{E}^{\theta_{j,k}},\mathcal{E}^{\alpha_{j,k}},\mathcal{E}^{\beta_{j,k}}$ be the extended symmetric forms induced by their corresponding functions.
    By Lemma \ref{lemma:communtativity of complex hessian operators} and (\ref{Equation:formula of Etheta when xi equals 1 near support of eta}), the integral in (\ref{Equation:second desired inequality in convergence thm under quasi-monotone convergence}) is equal to
    \begin{align*}
        -\mathcal{E}^{\theta_{j,k}}(\eta u_{0,k}, \xi \hat{u}_{j,k} - \xi \hat{u}_j).
    \end{align*}
    We now prove $\lim_{k\to\infty}\mathcal{E}^{\theta_{j,k}}(\eta u_{0,k}, \xi \hat{u}_{j,k} - \xi \hat{u}_{j})=0$,
    which can be implied by
    \begin{align}\label{Equation: last two equations need to prove in the weak convergence for V-quasi-monotone}
        \sup_{k\in\N}\mathcal{E}^{\theta_{j,k}}(\eta u_{0,k},\eta u_{0,k})<\infty,\quad\lim_{k\to\infty}\mathcal{E}^{\theta_{j,k}}(\xi \hat{u}_{j,k} - \xi \hat{u}_{j},\xi \hat{u}_{j,k} - \xi \hat{u}_{j})=0.
    \end{align}
    The former assertion in (\ref{Equation: last two equations need to prove in the weak convergence for V-quasi-monotone}) follows by the same argument used in the proof of Theorem \ref{Theorem:convergence of complex Hessian operatiors under monotone sequences}. By Remark \ref{remark: Etheta-convergence for bounded functions}, the latter assertion in (\ref{Equation: last two equations need to prove in the weak convergence for V-quasi-monotone}) can be implied by
    \begin{align*}
        &\limsup_{k\to\infty} \left(\mathcal{E}^{\theta_{j,k}}(\xi \hat{u}_{j,k}, \xi \hat{u}_{j,k})-\mathcal{E}^{\theta_{j,k}}(\xi \hat{u}_{j}, \xi \hat{u}_j)\right)\leq 0,\\
        &\limsup_{k \rightarrow \infty}\mathcal{E}^{\theta_{j,k}}(\xi (\hat{u}_j-\hat{u}_{j,k}), \xi \hat{u}_j)\leq 0.
    \end{align*}
    By Theorem \ref{Theorem: Monotonicity for bounded functions}, we have
    \begin{align*}
        &\mathcal{E}^{\theta_{j,k}}(\xi \hat{u}_{j,k}, \xi \hat{u}_{j,k})-\mathcal{E}^{\theta_{j,k}}(\xi \hat{u}_{j}, \xi \hat{u}_j)\leq \mathcal{E}^{\alpha_{j,k}}(\xi \hat{\varphi}_{j,k}, \xi \hat{\varphi}_{j,k})-\mathcal{E}^{\beta_{j,k}}(\xi \hat{\psi}_{j,k}, \xi \hat{\psi}_{j,k}),\\
        &\mathcal{E}^{\theta_{j,k}}(\xi (\hat{u}_j-\hat{u}_{j,k}), \xi \hat{u}_j)\leq \mathcal{E}^{\alpha_{j,k}}(\xi \hat{\varphi}_{j,k}, \xi \hat{\varphi}_{j,k})-\mathcal{E}^{\beta_{j,k}}(\xi \hat{\psi}_{j,k}, \xi \hat{\psi}_{j,k}).
    \end{align*}
    Therefore, we only need to prove that
    \begin{align*}
        \lim_{k\to\infty}\mathcal{E}^{\alpha_{j,k}}(\xi \hat{\varphi}_{j,k}, \xi \hat{\varphi}_{j,k})=\mathcal{E}^{\theta_{j}}(\xi \hat{u}_{j}, \xi \hat{u}_{j})=\lim_{k\to\infty}\mathcal{E}^{\beta_{j,k}}(\xi \hat{\psi}_{j,k}, \xi \hat{\psi}_{j,k}).
    \end{align*}
    We prove the first equality above, and the proof of the second equality is the same.
    On the one hand,
    \begin{align*}
        &\mathcal{E}^{\alpha_{j, k}}(\xi \hat{\varphi}_{j,k}, \xi \hat{\varphi}_{j,k})-\mathcal{E}^{\theta_{j}}(\xi \hat{\varphi}_{j,k}, \xi \hat{\varphi}_{j,k})\\
        &= \lim_{\delta\to 0}\left(\mathcal{E}^{\alpha_{j, k}}(\xi \hat{\varphi}_{j,k}^{[\delta]}, \xi \hat{\varphi}_{j,k}^{[\delta]})-\mathcal{E}^{\theta_{j}}(\xi \hat{\varphi}_{j,k}^{[\delta]}, \xi \hat{\varphi}_{j,k}^{[\delta]})\right)\\
        &= -\lim_{\delta\to 0}\int_\Omega \xi \hat{\varphi}^{[\delta]}_{j,k} dd^c(\xi \hat{\varphi}^{[\delta]}_{j,k})\wedge(\alpha_{j,k}-\theta_j)\\
        &= -\lim_{\delta\to 0}\sum_{r=1}^{j-1}\int_\Omega \xi \hat{\varphi}^{[\delta]}_{j,k} dd^c(\xi \hat{\varphi}^{[\delta]}_{j,k})\wedge dd^c \hat{\varphi}_{1,k}\wedge\cdots\wedge dd^c\hat{\varphi}_{r-1,k}\wedge dd^c (\hat{\varphi}_{r,k}-\hat{u}_r)\wedge\\
        &\quad dd^c \hat{u}_{r+1}\wedge\cdots\wedge dd^c \hat{u}_{j-1} \wedge dd^c \hat{u}_{j+1}\wedge\cdots\wedge dd^c \hat{u}_m\wedge\omega^{n-m}\\
        &= -\lim_{\delta\to 0}\sum_{r=1}^{j-1}\int_\Omega (\hat{\varphi}_{r,k}-\hat{u}_r) \left(dd^c(\xi \hat{\varphi}^{[\delta]}_{j,k})\right)^2\wedge dd^c \hat{\varphi}_{1,k}\wedge\cdots\wedge dd^c\hat{\varphi}_{r-1,k}\wedge\\
        & \quad dd^c \hat{u}_{r+1}\wedge\cdots\wedge dd^c \hat{u}_{j-1} \wedge dd^c \hat{u}_{j+1}\wedge\cdots\wedge dd^c \hat{u}_m\wedge\omega^{n-m}\\
        &= -\lim_{\delta\to 0}\sum_{r=1}^{j-1}\int_E \gamma(\hat{\varphi}_{r,k}-\hat{u}_r) \left(dd^c \hat{\varphi}^{[\delta]}_{j,k}\right)^2\wedge dd^c \hat{\varphi}_{1,k}\wedge\cdots\wedge dd^c\hat{\varphi}_{r-1,k}\wedge\\
        & \quad dd^c \hat{u}_{r+1}\wedge\cdots\wedge dd^c \hat{u}_{j-1} \wedge dd^c \hat{u}_{j+1}\wedge\cdots\wedge dd^c \hat{u}_m\wedge\omega^{n-m}\\
        &= -\sum_{r=1}^{j-1}\int_E \gamma(\hat{\varphi}_{r,k}-\hat{u}_r) (dd^c \hat{\varphi}_{j,k})^2\wedge dd^c \hat{\varphi}_{1,k}\wedge\cdots\wedge dd^c\hat{\varphi}_{r-1,k}\wedge\\
        & \quad dd^c \hat{u}_{r+1}\wedge\cdots\wedge dd^c \hat{u}_{j-1} \wedge dd^c \hat{u}_{j+1}\wedge\cdots\wedge dd^c \hat{u}_m\wedge\omega^{n-m},
    \end{align*}
    which tends to $0$ as $k\to\infty$ by Theorem \ref{Theorem:convergence of complex Hessian operatiors under monotone sequences}. On the other hand,
    by (\ref{equation:approximation of Etheta under monotone sequence}) we have $$\lim_{k\to\infty}\mathcal{E}^{\theta_{j}}(\xi \hat{\varphi}_{j,k}, \xi \hat{\varphi}_{j,k})=\mathcal{E}^{\theta_{j}}(\xi \hat{u}_{j}, \xi \hat{u}_{j}).$$
    Therefore, we have
    \begin{align*}
        \lim_{k\to\infty}\mathcal{E}^{\alpha_{j,k}}(\xi \hat{\varphi}_{j,k}, \xi \hat{\varphi}_{j,k})=\mathcal{E}^{\theta_{j}}(\xi \hat{u}_{j}, \xi \hat{u}_{j}).
    \end{align*}
    This finishes the proof.
\end{proof}

\begin{remark}\label{Remark:convergence of complex Hessian under BT capacity convergence}
    Theorem \ref{theorem:continuity of hessian operator} remains valid when $V$-quasi-monotone convergence is replaced by convergence in $\mathrm{Cap}_{BT,m}$. The proof follows the same argument, and uses the result in Remark \ref{Remark: Etheat convergence and BT capacity}. We sketch it below, retaining the notation of that proof. For $j=1,\dots,m$ and $k\in\N$, set
    $$\psi_{j,k}:=\left(\sup_{\ell\geq k}u_{j,\ell}\right)^*,$$
    which is decreasing in $k$ and converges to $u_j$ as $k\to \infty$.
 
    We only need to show
    \begin{align*}
        \lim_{k\to\infty}\mathcal{E}^{\theta_{j,k}}(\xi \hat{u}_{j,k} - \xi \hat{u}_{j},\xi \hat{u}_{j,k} - \xi \hat{u}_{j})=0,
    \end{align*}
    which can be implied by
    \begin{align}\label{goal 1}
        \limsup_{k\to\infty} \left(\mathcal{E}^{\theta_{j,k}}(\xi \hat{u}_{j}, \xi \hat{u}_j)-\mathcal{E}^{\theta_{j,k}}(\xi \hat{u}_{j,k}, \xi \hat{u}_{j,k})\right)\leq 0
    \end{align}
    and
    \begin{align}\label{goal 2}
        \limsup_{k \rightarrow \infty}\mathcal{E}^{\theta_{j,k}}(\xi (\hat{u}_{j,k}-\hat{u}_j), \xi \hat{u}_{j,k})\leq 0.
    \end{align}

    Following the same argument as in Remark \ref{Remark: Etheat convergence and BT capacity}, we have
    \begin{align*}
        \lim_{k\to\infty}\mathcal{E}^{\theta_{j,k}}(\xi (\hat{u}_{j,k}-\hat{u}_j), \xi \hat{u}_{j,k})=\lim_{k\to\infty}\int_\Omega \xi (\hat{u}_{j}-\hat{u}_{j,k}) dd^c \hat{u}_{j,k}\wedge \theta_{j,k}=0,
    \end{align*}
    since $\hat{u}_{j,k}$ converges to $\hat{u}_{j}$ in $\mathrm{Cap}_{BT,m}$.
   
    We now prove (\ref{goal 1}). We have
    \begin{align*}
       \mathcal{E}^{\theta_{j,k}}(\xi \hat{u}_{j}, \xi \hat{u}_j)-\mathcal{E}^{\theta_{j,k}}(\xi \hat{u}_{j,k}, \xi \hat{u}_{j,k})\leq \mathcal{E}^{\theta_{j,k}}(\xi \hat{u}_{j}, \xi \hat{u}_j)-\mathcal{E}^{\beta_{j, k}}(\xi \hat{\psi}_{j,k}, \xi \hat{\psi}_{j,k}).
    \end{align*}
    Hence, (\ref{goal 1}) follows from the following equalities.
    \begin{align*}
        \lim_{k\to\infty}\mathcal{E}^{\theta_{j,k}}(\xi \hat{u}_{j}, \xi \hat{u}_j)
        &=\mathcal{E}^{\theta_{j}}(\xi \hat{u}_{j}, \xi \hat{u}_j)\\
        &=\lim_{k\to\infty}\mathcal{E}^{\theta_{j}}(\xi \hat{\psi}_{j,k}, \xi \hat{\psi}_{j,k})\\
        &=\lim_{k\to\infty}\mathcal{E}^{\beta_{j, k}}(\xi \hat{\psi}_{j,k}, \xi \hat{\psi}_{j,k}).
    \end{align*}
    The second equality follows from (\ref{equation:approximation of Etheta under monotone sequence}). As for the third equality, by the same computation as in the proof of Theorem \ref{theorem:continuity of hessian operator}, we have
    \begin{align*}
        &\mathcal{E}^{\beta_{j, k}}(\xi \hat{\psi}_{j,k}, \xi \hat{\psi}_{j,k})-\mathcal{E}^{\theta_{j}}(\xi \hat{\psi}_{j,k}, \xi \hat{\psi}_{j,k})\\
        &= -\sum_{r=1}^{j-1}\int_E \gamma(\hat{\psi}_{r,k}-\hat{u}_r) (dd^c \hat{\psi}_{j,k})^2\wedge dd^c \hat{\psi}_{1,k}\wedge\cdots\wedge dd^c\hat{\psi}_{r-1,k}\wedge\\
        & \quad dd^c \hat{u}_{r+1}\wedge\cdots\wedge dd^c \hat{u}_{j-1} \wedge dd^c \hat{u}_{j+1}\wedge\cdots\wedge dd^c \hat{u}_m\wedge\omega^{n-m},
    \end{align*}
    which tends to $0$ as $k\to\infty$ by Theorem \ref{Theorem:convergence of complex Hessian operatiors under monotone sequences}. As for the first  equality,
    \begin{align*}
        &\mathcal{E}^{\theta_{j,k}}(\xi \hat{u}_{j}, \xi \hat{u}_j)-\mathcal{E}^{\theta_{j}}(\xi \hat{u}_{j}, \xi \hat{u}_j)\\
        &= -\sum_{r=1}^{j-1}\int_E \gamma(\hat{u}_{r,k}-\hat{u}_r) (dd^c \hat{u}_{j})^2\wedge dd^c \hat{u}_{1,k}\wedge\cdots\wedge dd^c\hat{u}_{r-1,k}\wedge\\
        & \quad dd^c \hat{u}_{r+1}\wedge\cdots\wedge dd^c \hat{u}_{j-1} \wedge dd^c \hat{u}_{j+1}\wedge\cdots\wedge dd^c \hat{u}_m\wedge\omega^{n-m},
    \end{align*}
    which tends to $0$ as $k\to\infty$, since $\hat{u}_{r,k}$ converges to $\hat{u}_r$ in $\mathrm{Cap}_{BT,m}$ for every $r$.
\end{remark}

\begin{remark}
    The commutativity and weak convergence results we proved in this section can be extended from top bi-degree to arbitrary bi-degree. 
    
    Indeed, every test $(p,p)$-form on $\Omega$ is a finite sum of terms of the form
    $$\eta dd^c |f_1|^2\wedge\dots\wedge dd^c|f_p|^2,$$
    where $\eta\in \mathscr{C}_c^\infty(\Omega,\C)$ and $f_1,\dots,f_p$ are $\C$-linear functions. To test either assertion against a test form, it suffices to consider each term as above. After wedging the currents with $dd^c|f_1|^2\wedge\cdots\wedge dd^c|f_p|^2$, we obtain currents of top bi-degree, and hence the preceding results apply, since $|f_1|^2,\dots,|f_p|^2$ are bounded and $m$-subharmonic.
\end{remark}

\section{Dirichlet forms and stochastic characterizations of $m$-polar sets}
In this section, we establish the new characterizations of $m$-polar sets in terms of Dirichlet forms and stochastic processes. We first review the Bedford--Taylor capacity $\mathrm{Cap}_{BT,m}$ for $m$-subharmonic functions, revisit the $\mathcal{E}^\theta$-quasi-continuity, and show the equivalence between $V$-quasi-monotone convergence and convergence in $\mathrm{Cap}_{BT,m}$. We then compare $\mathrm{Cap}_{BT,m}$ with $\mathrm{Cap}_\theta$ and express the Choquet capacity $\mathrm{C}_{\sharp,m}$ in terms of the conformal diffusion associated with the Dirichlet form. Finally, we establish the characterizations of $m$-polar sets. 

\begin{definition}
    For every Borel subset $E\subset\Omega$, define
    $$\mathrm{Cap}_{BT,m}(E,\Omega):=\sup\left\{\int_E(dd^c u)^m\wedge\omega^{n-m}:u\in SH_m(\Omega), -1\leq u\leq 0\right\}.$$
    For every subset $E\subset\Omega$, define
    $$\mathrm{Cap}^*_{BT,m}(E,\Omega):=\inf \left\{ \mathrm{Cap}_{BT,m}(D,\Omega): E\subset D\subset \Omega, D~\text{open} \right\}.$$
\end{definition}

$\mathrm{Cap}_{BT,m}^*$ is a Choquet capacity on $\Omega$ with countably subadditive property, and
$$\mathrm{Cap}_{BT,m}^*(E,\Omega)=\mathrm{Cap}_{BT,m}(E,\Omega)$$
for every Borel subset $E\subset\Omega$. We call $\mathrm{Cap}_{BT,m}^*$ the Bedford--Taylor capacity for $m$-subharmonic functions relative to $\Omega$. From now on, we omit the superscript $*$ and write $\mathrm{Cap}_{BT,m}$ for this capacity.
The following result explains the relationship between $\mathrm{Cap}_{BT, m}$ and the $m$-relative extremal function. 

\begin{theorem}\cite[Theorem 1.4.12]{Lu12}, \cite[Theorem 2.20]{Lu15JMAA}
\label{Theorem:Lu}
    Let $\Omega$ be a bounded strictly $m$-pseudoconvex domain. Let $E\Subset \Omega$ be a Borel subset and $K\subset\Omega$ be a compact subset. Then
    \begin{align*}
        &\mathrm{Cap}_{BT,m}(E,\Omega)=\int_\Omega(dd^cu_{E,\Omega,m}^*)^m\wedge\omega^{n-m},\\
        &\mathrm{Cap}_{BT,m}(K,\Omega)=\int_K(dd^cu_{K,\Omega,m}^*)^m\wedge\omega^{n-m}.
    \end{align*}
\end{theorem}

We recall the following classical characterization of $m$-polar sets.
\begin{theorem}\cite[Theorems 1.6.6 and 1.6.7, Corollaries 1.6.4 and 1.6.10]{Lu12}
\label{theorem: Lu's characterization}
    Let $\Omega\subset\C^n$ be a bounded open subset and $E\subset\Omega$. Then the following statements are equivalent.
    \begin{itemize}
        \item[(i)] $u^*_{E,\Omega,m}\equiv0$.
        \item[(ii)] $E$ is $m$-polar.
        \item[(iii)] $E$ is $m$-negligible.
        \item[(iv)] $\mathrm{Cap}_{BT,m}(E,\Omega)=0$.
    \end{itemize}
\end{theorem}

\begin{definition}\label{def:BT-capacity-convergence}
    Let $\Omega\subset\C^n$ be a bounded open subset, $(u_k)_{k\in\N}\subset SH_m(\Omega)$, and
    $u\in SH_m(\Omega)$. We say that $u_k$ converges to $u$
    in $\mathrm{Cap}_{BT,m}$ relative to $\Omega$ if, for every compact subset
    $K\subset\Omega$ and every constant $\delta>0$,
    \[
        \lim_{k\to\infty}
        \mathrm{Cap}_{BT,m}
        \left(K\cap\{|u_k-u|>\delta\},\Omega\right)=0.
    \]
    Here $u_k-u$ may be defined arbitrarily on
    $\{u_k=-\infty\}\cup\{u=-\infty\}$, since this set has zero capacity $\mathrm{Cap}_{BT,m}$ and therefore does not affect the definition.
\end{definition}

We now revisit the $\mathcal{E}^\theta$-quasi-continuity of locally bounded $m$-subharmonic functions established in Theorem \ref{bounded msh in F loc}, extending the result to unbounded $m$-subharmonic functions and providing equivalent characterizations.

\begin{theorem}\label{Theorem for quasi continuity}
    Let $\theta\in \mathcal{M}_m(\Omega)$, $(\mathcal{F}^{\theta}, \mathcal{E}^{\theta})$ be the Dirichlet space for an admissible pair $(\theta, \mu)$, and $u \in SH_m(\Omega)$. Then the following statements are equivalent.
    \begin{itemize}
        \item [(i)] $u$ is $\mathcal{E}^{\theta}$-q.c.
        \item [(ii)] $\mu\left(\{u = -\infty\}\right) = 0$.
        \item [(iii)] $\mathrm{Cap}_{\theta}(\{u = -\infty\}) = 0$.
    \end{itemize}
\end{theorem}

\begin{proof}
    (i) $\Rightarrow$ (iii) It follows from the definition of $\mathcal{E}^\theta$-q.c.
    
    (iii) $\Rightarrow$ (ii) $\left\{u=-\infty \right\}\subset\Omega$ is Borel, hence $\mu\left(\{u = -\infty\}\right)\leq\mathrm{Cap}_{\theta}(\{u = -\infty\})$.

    (ii) $\Rightarrow$ (i) Fix an arbitrary open subset $D \Subset \Omega$. Let $\xi \in \mathscr{C}_c^{\infty}(\Omega, \R)$ so that $\xi = 1$ on $D$. Let $G \Subset \Omega$ be an open neighborhood of $\supp{\xi}$. Let $M$ be a constant so that $u \leq M$ on $G$. We first show that
    \begin{align}\label{convergence of sublevel set with respect to cap}
        \lim_{k \to \infty}\mathrm{Cap}_{\theta}\left(D \cap \{u - M < -k\}\right) = 0.
    \end{align}

    For every $k \in \N$, we have
    \begin{align*}
        u_k := \max \left(\frac{u-M}{k}, -1\right)\in SH_m(\Omega) \cap L^{\infty}_{loc}(\Omega).
    \end{align*}
    It follows from Theorem \ref{bounded msh in F loc} that $-\xi u_k \in \mathcal{F}^{\theta}$. Since $D \cap \{u - M < -k\}$ is an open subset of $\Omega$ and $-\xi u_k = 1$ on $D \cap \{u - M < -k\}$,
    \begin{align*}
        \mathrm{Cap}_{\theta}(D \cap \{u-M < -k\})\leq \mathcal{E}_1^{\theta}(-\xi u_k, -\xi u_k)=\mathcal{E}_1^{\theta}(\xi u_k, \xi u_k).
    \end{align*}
    Denote by $(\mathcal{F}^{\theta|_G},\mathcal{E}^{\theta|_G})$ the Dirichlet space for the admissible pair $(\theta|_G,\mu|_G)$. Since $\supp\xi \subset G$, we have $\xi u_k\in\mathcal{F}^{\theta|_G}$. By the same computation as (\ref{equation: admissible pair and its restriction}), we have
    $$\mathcal{E}_1^{\theta}(\xi u_k, \xi u_k)=\mathcal{E}_1^{\theta|_G}(\xi u_k, \xi u_k).$$
    Notice that $(u_k)_{k\in\N}$ is uniformly bounded on $G$, increasing on $G$, and converges to $0$ on $G\setminus \{u=-\infty\}$ where $\{u=-\infty\}$ is of zero measure $\mu$. It follows from Theorem \ref{Theorem:Cauchy for locally uniformly bounded and monotone} that $\lim_{k\to\infty}\mathcal{E}_1^{\theta|_G}(\xi u_k, \xi u_k)=0$. Therefore, we obtain (\ref{convergence of sublevel set with respect to cap}).
    
    For every $\varepsilon>0$, there exists $k\in\N$ such that 
    $$\mathrm{Cap}_{\theta}(D \cap \{u-M < -k\}) < \varepsilon.$$
    Since $u_k$ is $\mathcal{E}^\theta$-q.c., there exists an open subset $E\subset D$ such that $\mathrm{Cap}_{\theta}(E) < \varepsilon$, and $u_k$ is finite-valued and continuous on $D\setminus E$. Set 
    $$U:= E \cup \left(D \cap \{u-M < -k\} \right).$$
    $U$ is open in $D$ and $\mathrm{Cap}_\theta (U)<2\varepsilon$. Notice that $u=k u_k+M$ on $D\setminus U$. Hence $u$ is finite-valued and continuous on $D\setminus U$. Finally, proceeding with the last part of the proof of Theorem \ref{bounded msh in F loc} (ii), we obtain the desired assertion.
\end{proof}

Since $m$-polar sets are globally $m$-polar \cite[Theorem 1.6.5]{Lu12}, Theorem \ref{Theorem for quasi continuity} and Theorem \ref{Theorem:example of admissible pair} imply the following result. Since $m$-polar sets and $m$-negligible sets are equivalent, the following result generalizes Theorem \ref{Theorem:m-negligible is E theta-polar} and Corollary \ref{Corollary: m-negligible are Etheta-polar when mu is V}.
\begin{corollary}\label{corollart:every msh function is Etheta qc}
    Let $\theta \in \mathcal{M}_m(\Omega)$ and $(\mathcal{F}^{\theta}, \mathcal{E}^{\theta})$ be the Dirichlet space for an admissible pair $(\theta, \mu)$. Then the following statements are equivalent.
    \begin{itemize}
        \item [(i)] Every function in $SH_m(\Omega)$ is $\mathcal{E}^\theta$-q.c.
        \item [(ii)] $\mu$ does not charge any $m$-polar set.
        \item [(iii)] Every $m$-polar set is $\mathcal{E}^\theta$-polar.
    \end{itemize}
    In particular, if $\theta\in \mathcal{M}_m^+(\Omega)$ and $\mu=V$, then every function in $SH_m(\Omega)$ is $\mathcal{E}^\theta$-q.c., and every $m$-polar set is $\mathcal{E}^\theta$-polar.
\end{corollary}
Now we show the equivalence between $V$-quasi-monotone convergence and convergence in $\mathrm{Cap}_{BT,m}$, whose proof
is inspired by \cite[Theorem 2.4]{GT23}.

\begin{lemma}\label{lemma:small capacity}
    Let $\Omega\subset\C^n$ be a bounded open subset, and $K\subset\Omega$ be a compact subset. Then for every $\varepsilon>0$, there exists $\delta>0$ such that for every subset $E\subset K$,
    $$\mathrm{Cap}_{BT,m}(E,\Omega)<\delta\Longrightarrow\mathrm{C}_{\sharp,m}(E)=\int_\Omega(-u^*_{E,\Omega,m}) dV<\varepsilon.$$
\end{lemma}

\begin{proof}
    We can assume that $E$ is Borel since $\mathrm{Cap}_{BT,m}$ and $\mathrm{C}_{\sharp,m}$ are outer regular.

    Assume first that $\Omega=B$ is a ball. Suppose that the desired assertion fails. Then there exist $\varepsilon_0>0$ and Borel subsets $E_k\subset K$ for $k\in\N$ such that
    $$\mathrm{Cap}_{BT,m}(E_k,B)\leq\frac{1}{2^k},\quad\int_B(-u^*_{E_k,B,m})dV\geq\varepsilon_0.$$
    For every $k\in\N$,
    $$F_k:=\bigcup_{\ell\geq k}E_\ell\subset K.$$
    Then
    $$\mathrm{Cap}_{BT,m}(F_k,B)\leq\frac{1}{2^{k-1}},~ ~ (F_k)_{k\in\N} \downarrow,~ ~u_{F_k,B,m}\uparrow,~ ~u:=\left(\sup_{k\in\N}u^*_{F_k,B,m}\right)^*\in SH_m(B).$$
    By Theorem \ref{Theorem:convergence of complex Hessian operatiors under monotone sequences} and Theorem \ref{Theorem:Lu}, for every nonnegative function $\eta\in\mathscr{C}_c^\infty (B,\mathbb{R})$,
    \begin{align*}
        0&\leq\int_B \eta (dd^cu)^m\wedge\omega^{n-m}\\
        &=\lim_{k\to\infty}\int_B\eta(dd^cu^*_{F_k,B,m})^m\wedge\omega^{n-m}\\
        &\leq \lVert \eta \rVert_{L^\infty(B)}\lim_{k\to\infty}\int_B(dd^cu^*_{F_k,B,m})^m\wedge\omega^{n-m}\\
        &= \lVert \eta \rVert_{L^\infty(B)} \lim_{k\to\infty}\mathrm{Cap}_{BT,m}(F_k,B)\\
        &=0,
    \end{align*}
    which means that $(dd^cu)^m\wedge\omega^{n-m}=0$ on $B$ in the sense of currents. Notice that
    $$u^*_{K,B,m}\leq u^*_{F_k,B,m}\leq u\leq0, \quad \lim_{z\to\partial B}u^*_{K,B,m}(z)=0,$$
    hence $\lim_{z\to\partial B}u(z)=0.$ By the comparison principle \cite[Corollary 1.15 (b)]{Nguyen13}, $u\equiv0$ on $B$, and hence
    $$\lim_{k\to\infty}\int_B(-u^*_{F_k,B,m})dV=0.$$
    But we also have
    $$\varepsilon_0\leq\int_B(-u^*_{E_k,B,m}) dV\leq\int_B(-u^*_{F_k,B,m}) dV.$$
    This is a contradiction. For an arbitrary bounded open subset $\Omega\subset\C^n$, choose a ball $B$ containing $\Omega$. For every Borel subset $E\subset K$,
    \begin{align*}
        &\mathrm{Cap}_{BT,m}(E,B)\leq \mathrm{Cap}_{BT,m}(E,\Omega),\\
        &\int_\Omega(-u^*_{E,\Omega,m})dV\leq \int_{\Omega}(-u^*_{E,B,m})dV\leq\int_B(-u^*_{E,B,m})dV.
    \end{align*}
    The desired assertion follows from the case of balls.
\end{proof}

\begin{proposition}\label{Proposition: equvialence of convergence quasi-monotonically and in BT capacity}
    Let $(u_k)_{k \in \N} \subset SH_m(\Omega)$ and $u \in SH_m(\Omega)$.
    \begin{itemize}
        \item [(i)] If $u_k$ converges $V$-quasi-monotonically to $u$, then it converges in $\operatorname{Cap}_{BT, m}$. 
        \item [(ii)] If $(u_k)_{k \in \N}$ is locally uniformly bounded and converges to $u$ in $\operatorname{Cap}_{BT, m}$, then a subsequence converges $V$-quasi-monotonically to $u$.
    \end{itemize}
\end{proposition}

\begin{proof}
    (i) Let $K\subset\Omega$ be a compact subset and $\delta>0$. Take an open subset $D\Subset\Omega$ containing $K$. Then there exist sequences $(\phi_k)_{k \in \N}, (\psi_k)_{k \in \N} \subset SH_m(D)$ that satisfy the following properties.
    \begin{itemize}
        \item [(i')] $\phi_k \leq u_k \leq \psi_k$ on $D$ for every $k\in\N$.
        \item [(ii')] $(\phi_k)_{k \in \N}$ is increasing and $(\psi_k)_{k \in \N}$ is decreasing.
        \item [(iii')] Both $\phi_k$ and $\psi_k$ converge to $u$ in $L^1_{loc}(D)$.
    \end{itemize}
    Moreover, $\phi_k$ converges to $u~V$-a.e. on $D$, $\psi_k$ converges to $u$ on $D$, and $\phi_k\leq u\leq \psi_k$ on $D$ for every $k\in\N$. Outside the $m$-polar set
    $\cup_{k=1}^\infty\{u_k=-\infty\}\cup\{u=-\infty\}$,
    $$K\cap \{|u_k - u|>\delta\}\subset\left( K\cap\{|\phi_k - u|>\delta\} \right)\cup \left( K\cap\{|\psi_k - u|>\delta\} \right).$$
    Since monotone convergence implies convergence in $\mathrm{Cap}_{BT,m}$, both $\phi_k$ and $\psi_k$ converge to $u$ in $\mathrm{Cap}_{BT,m}$ relative to $D$. Then
    \begin{align*}
        &\mathrm{Cap}_{BT,m}\left(K\cap \{|u_k - u|>\delta\},\Omega\right)\\
        &\leq \mathrm{Cap}_{BT,m}\left(K\cap \{|u_k - u|>\delta\},D\right)\\
        &\leq \mathrm{Cap}_{BT,m}\left(K\cap \{|\phi_k - u|>\delta\},D\right) + \mathrm{Cap}_{BT,m}\left(K\cap \{|\psi_k - u|>\delta\},D\right),
    \end{align*}
    which tends to $0$ as $k\to\infty$.

    (ii) By the conditions, $u$ is locally bounded and $u_k$ converges to $u$ in $L_{loc}^1(\Omega)$.

    Fix an arbitrary open subset $D\Subset\Omega$. Our first goal is to show that there exist a subsequence $(u_{j_k})_{k\in\N}$ of $(u_k)_{k\in\N}$ and sequences $(\phi_k)_{k\in\N},(\psi_k)_{k\in\N}\subset SH_m(D)$ that satisfy the following properties.
    \begin{itemize}
        \item[(i')] $\phi_k\leq u_{j_k}\leq \psi_k$ on $D$ for every $k\in\N$.
        \item[(ii')] $(\phi_k)_{k \in \N}$ is increasing and $(\psi_k)_{k \in \N}$ is decreasing.
        \item[(iii')]  Both $\phi_k$ and $\psi_k$ converge to $u$ in $L^1_{loc}(D)$.
    \end{itemize} 

    We can assume that $-1\leq u_k,u\leq 0$ on $\overline{D}$ for every $k\in\N$. By Lemma \ref{lemma:small capacity}, there exists a constant $\delta_k>0$ for every $k\in\N$ such that for every Borel subset $E\subset\overline{D}$,
    $$\mathrm{Cap}_{BT,m}(E,\Omega)<\delta_k\Longrightarrow\int_\Omega(-u^*_{E,\Omega,m}) dV<\frac{1}{2^k}.$$
    Since $u_k$ converges to $u$ in $\mathrm{Cap}_{BT,m}$, there exists a subsequence $(u_{j_k})_{k\in\N}$ of $(u_k)_{k\in\N}$ such that
    $$\mathrm{Cap}_{BT,m}\left(\overline{D}\cap\left\{|u_{j_k}-u_{j_{k+1}}|>\frac{1}{2^k}\right\},\Omega\right)<\delta_k.$$
    For every $k\in\N$, set
    $$E_k:=\overline{D}\cap\left\{u_{j_{k+1}}<u_{j_k}-\frac{1}{2^k}\right\},\quad H_k:= \sum_{\ell=k}^\infty u^*_{E_\ell,\Omega,m}.$$
    Then 
    $$\int_\Omega(-u^*_{E_k,\Omega,m}) dV<\frac{1}{2^k},\quad \int_\Omega |H_k|dV\leq \frac{1}{2^{k-1}}.$$
    Hence, $H_k\in SH_m(\Omega)$, $(H_k)_{k\in\N}$ is increasing, and $H_k$ converges to $0~V$-a.e. on $\Omega$. For every $k\in\N$, set
    $$\phi_k:= u_{j_k}+H_k -\frac{1}{2^{k-1}}.$$
    Then $\phi_k \in SH_m(\Omega)$, $\phi_k\leq u_{j_k}$, and $\phi_k$ converges to $u$ in $L_{loc}^1(\Omega)$. For every $k\in \N$, outside the zero Lebesgue measure set $\{H_k=-\infty\}\cup\{H_{k+1}=-\infty\}$ we have
    $$\phi_{k+1}-\phi_k=u_{j_{k+1}}-u_{j_k}-u^*_{E_k,\Omega,m}+\frac{1}{2^k}.$$
    On $D\setminus E_k$, 
    $$u_{j_{k+1}}-u_{j_k}-u^*_{E_k,\Omega,m}+\frac{1}{2^k}\geq -u^*_{E_k,\Omega,m}\geq 0.$$
    On $D\cap E_k$,
    $$u_{j_{k+1}}-u_{j_k}-u^*_{E_k,\Omega,m}+\frac{1}{2^k}\geq -1 - u^*_{E_k,\Omega,m}.$$
    Since $u^*_{E_k,\Omega,m}=-1~V$-a.e. on $E_k$, $u_{j_{k+1}}-u_{j_k}-u^*_{E_k,\Omega,m}+\frac{1}{2^k}\geq 0~V$-a.e. on $D\cap E_k$.
    Therefore, $\phi_{k+1}\geq \phi_k~V$-a.e. on $D$, and hence $\phi_{k+1}\geq \phi_k$ on $D$. For every $k\in\N$,
    $$u_{j_k}\leq\left(\sup_{\ell\geq k} u_{j_\ell}\right)^*=:\psi_k\in SH_m(\Omega).$$
    $(\psi_k)_{k\in\N}$ is decreasing. Since $(u_{j_k})_{k\in\N}$ converges to $u$ in $L_{loc}^1(\Omega)$, $\psi_k$ converges to $u$ on $\Omega$. Therefore, $(\psi_k)_{k\in\N}$ converges to $u$ in $L^1_{loc}(\Omega)$. We achieved the first goal.

    For every $j\in\N$, let $D_j\Subset\Omega$ be an open subset so that $D_1\Subset D_2\Subset\dots$ and $\Omega=\cup_{j=1}^\infty D_j$. Applying our achieved first goal on $D_1,D_2,\dots$ each time to the subsequence obtained in the preceding step, we obtain a subsequence $I_j=\{i_{j,1},i_{j,2},\dots\}$ of $\N$ for every $j\in\N$, such that every $I_{j+1}$ is a subsequence of $I_j$, and obtain $(\phi_\ell^j)_{\ell\in\N},(\psi_\ell^j)_{\ell\in\N}\subset SH_m(D_j)$ that satisfy the following properties.
    \begin{itemize}
        \item[(i'')] $\phi_\ell^j\leq u_{i_{j,\ell}}\leq \psi_\ell^j$ on $D_j$ for every $\ell\in\N$.
        \item[(ii'')] $(\phi_\ell^j)_{\ell \in \N}$ is increasing and $(\psi_\ell^j)_{\ell \in \N}$ is decreasing.
        \item[(iii'')]  Both $\phi_\ell^j$ and $\psi_\ell^j$ converge to $u$ in $L^1_{loc}(D_j)$.
    \end{itemize} 
    Take a subsequence $(s_k)_{k\in\N}$ of $\N$ so that $s_k\in I_k$ for every $k$. We show that $u_{s_k}$ converges $V$-quasi-monotonically to $u$ on $\Omega$. Since every compact subset of $\Omega$ is contained in some $D_j$, we only need to check the behavior of $(u_{s_k})_{k\in\N}$ on every $D_j$.
    
    Fix an arbitrary $j\in N$. For every $k\geq j$, we have $s_k\in I_k \subset I_j$, so there exists  $r_j(k)\in \N$ such that $s_k=i_{j,r_j(k)}$. Notice that $\left(r_j(k)\right)_{k\geq j}$ is increasing. Set
    $$\Phi_k^j:=\phi^j_{r_j(k)},\quad \Psi_k^j:=\psi^j_{r_j(k)},\quad \forall k\geq j.$$
    Then the functions $(\Phi_k^j)_{k\geq j}, (\Psi_k^j)_{k\geq j}\subset SH_m(D_j)$ satisfy the following properties.
    \begin{itemize}
        \item[(a)] $\Phi_k^j\leq u_{s_k}\leq \Psi_k^j$ on $D_j$ for every $k\geq j$.
        \item[(b)] $(\Phi_k^j)_{k\geq j}$ is increasing and $(\Psi_k^j)_{k\geq j}$ is decreasing.
        \item[(c)] Both $\Phi_k^j$ and $\Psi_k^j$ converge to $u$ in $L^1_{loc}(D_j)$.  
    \end{itemize}
    We set $\Phi_k^j$ and $\Psi_k^j$ for $1\leq k\leq j-1$. Notice that (a) implies that $\Phi^j_{j}$ is bounded from above on $D_j$. Let $C_j$ be a constant so that $\Phi^j_{j},u_{s_1},\dots,u_{s_{j-1}}\leq C_j$ on $D_j$. Set
    \begin{align*}
        &\Phi^j_1:=\dots:=\Phi^j_{j-1}:=\Phi^j_j+ (u_{s_1}-C_j)+\dots+(u_{s_{j-1}}-C_j),\\
        &\Psi^j_1:=\dots:=\Psi^j_{j-1}:=\max(\Psi^j_j,u_{s_1},\dots,u_{s_{j-1}}).
    \end{align*}
    Then $\Phi^j_1,\dots,\Phi^j_{j-1},\Psi^j_1,\dots,\Psi^j_{j-1}\in SH_m(D_j)$ and satisfy the following properties.
    \begin{itemize}
        \item[(d)] $\Phi_k^j\leq u_{s_k}\leq \Psi_k^j$ on $D_j$ for every $1\leq k\leq j-1$.
        \item[(e)] $\Phi^j_1=\dots=\Phi^j_{j-1}\leq \Phi^j_j$ and $\Psi^j_1=\dots=\Psi^j_{j-1}\geq \Psi^j_j$.
    \end{itemize}
    This finishes the proof.
\end{proof}

We then give the following result as a variant of (\ref{Equation:corollary of Schwarz inequality}) which will be used later.
Let $\theta\in \mathcal{M}_m(\Omega)$ and $u \in SH_m(\Omega)\cap L_{loc}^{\infty}(\Omega)$. Recall the definition that
\begin{align}\label{Equation: define du wedge dcu}
    d u\wedge d^c u\wedge \theta:= \frac{1}{2} dd^c (u-M)^2\wedge\theta - (u-M)dd^c u\wedge\theta,
\end{align}
where $M$ is a local lower bound of $u$. Then $d u\wedge d^c u\wedge\theta$ is a well-defined closed positive $(n,n)$-current on $\Omega$. In particular, by Theorem \ref{Theorem:convergence of complex Hessian operatiors under monotone sequences} we have
$$d u^{[\frac{1}{k}]}\wedge d^c u^{[\frac{1}{k}]}\wedge\theta \to d u \wedge d^c u \wedge \theta,\quad k\to \infty$$
on every relatively compact open subset of $\Omega$.

\begin{lemma}\label{Lemma: variant of Schwarz inequality for bdd m-sh functions}
    Let $\theta\in \mathcal{M}_m(\Omega)$, and $(\mathcal{F}^{\theta}, \mathcal{E}^{\theta})$ be the Dirichlet space for an admissible pair $(\theta, \mu)$. Let $u \in SH_m(\Omega)\cap L_{loc}^{\infty}(\Omega)$ and $\eta \in \mathscr{C}_c^{\infty}(\Omega,\R)$. Then we have
    \begin{align*}
        \mathcal{E}^\theta(\eta u, \eta u)\leq 2\int_\Omega \eta^2 d u \wedge d^c u \wedge \theta +  2\int_\Omega u^2 d\eta \wedge d^c\eta \wedge \theta.
    \end{align*}
\end{lemma}
\begin{proof}
    Let $D\Subset\Omega$ be an open subset containing $\supp\eta$. For $k\in\N$ large enough, $u^{[\frac{1}{k}]}\in SH_m(D)\cap \mathscr{C}^\infty(D,\R)$. Then by (\ref{Equation:corollary of Schwarz inequality}) we have
    \begin{align*}
        \mathcal{E}^\theta(\eta u^{[\frac{1}{k}]},\eta u^{[\frac{1}{k}]})&=\int_D d(\eta u^{[\frac{1}{k}]})\wedge d^c(\eta u^{[\frac{1}{k}]})\wedge\theta\\
        &\leq 2 \int_D \eta^2 d u^{[\frac{1}{k}]}\wedge d^c u^{[\frac{1}{k}]}\wedge\theta + 2\int_D (u^{[\frac{1}{k}]})^2 d\eta\wedge d^c \eta\wedge\theta.
    \end{align*}
    Letting $k\to\infty$ finishes the proof.
\end{proof}

We now establish the quantitative comparison between $\mathrm{Cap}_{BT,m}$ and $\mathrm{Cap}_{\theta}$ for $\theta\in\mathcal{M}^+_m(\Omega)$, where $\mathrm{Cap}_{\theta}$ is the capacity induced by the Dirichlet space $(\mathcal{F}^\theta, \mathcal{E}^\theta)$ for the admissible pair $(\theta,V)$. Set $\gamma:=\sup_{\Omega}|z|^2$. For every $\theta\in\mathcal{M}_m^+(\Omega)$, write
\[\theta=dd^c(v_1+\delta|z|^2)\wedge\cdots\wedge dd^c(v_{m-1}+\delta|z|^2)\wedge\omega^{n-m},\]
where $v_1,\dots,v_{m-1}\in SH_m(\Omega)\cap L^\infty(\Omega)$ and $\delta>0$ is a constant. By Lemma \ref{Poincare type inequality}, we obtain the following estimate similar to \cite[(2.11)]{FO87}. For every $\varphi\in \mathscr{C}_c^\infty(\Omega,\R)$, 
\begin{align}\label{poincare inequality for lebesgue measure}
    &\int_\Omega\varphi^2 dV=\frac{\pi^n}{2^n n!}\int_\Omega\varphi^2 \omega^{n}\leq\frac{8\gamma\pi^n}{2^n n!}\mathcal{E}^{\omega^{n-1}}(\varphi,\varphi)\leq\frac{8\gamma\pi^n}{2^n \delta^{m-1}n!}\mathcal{E}^\theta(\varphi,\varphi).
\end{align}
Since $\mathscr{C}_c^\infty(\Omega,\R)$ is dense in the Hilbert space $(\mathcal{F}^\theta, \mathcal{E}^\theta_1)$, we have
\begin{align}\label{equation: Etheta and Etheta1}
    \mathcal{E}^\theta_1(u,u)\leq \left(1+ \frac{8\gamma\pi^n}{2^n\delta^{m-1} n!}\right)\mathcal{E}^\theta(u,u),\quad \forall u\in \mathcal{F}^\theta.
\end{align}
The following result is inspired by
\cite[Lemma 7.8]{Nguyen26}, whose argument extends to
our setting and yields the quantitative comparison between
$\mathrm{Cap}_{BT,m}$ and $\mathrm{Cap}_\theta$.

\begin{lemma}\label{equivalence of caps}
    Let $\Omega\subset\C^n$ be a bounded strictly $m$-pseudoconvex domain.
    \begin{itemize}
        \item[(i)] Let $K\subset\Omega$ be a compact subset and $\delta>0$ be a constant. Then
        \begin{align*}
            \mathrm{Cap}_{BT,m}(K, \Omega)\leq 8 \mathrm{Cap}_\theta(K),
        \end{align*}
        where $\theta:=\left(dd^c (u_{K,\Omega,m}^*+\delta|z|^2)\right)^{m-1}\wedge\omega^{n-m}\in \mathcal{M}_m^+(\Omega)$.
        \item[(ii)] Let $D\Subset\Omega$ be an open subset and $B,\delta>0$ be constants. Then there exists a constant $C$ that depends only on $\Omega,D,B,\delta,m$, such that for every open subset $G\subset D$ and every $v\in SH_m(\Omega)$ with $\lVert v\rVert_{L^\infty(\Omega)}\leq B$, we have
        \[\mathrm{Cap}_\theta(G)\leq C\left(\mathrm{Cap}_{BT,m}(G, \Omega)\right)^{\frac{1}{m}},\]
        where $\theta:=\left(dd^c(v+\delta|z|^2)\right)^{m-1}\wedge\omega^{n-m}\in \mathcal{M}_m^+(\Omega)$.
    \end{itemize}
\end{lemma}
\begin{proof}
    (i) Since $\mathscr{C}_c^\infty(\Omega,\R)$ is a special standard core of $\mathcal{E}^\theta$, it follows from \cite[Lemma 2.2.7 (ii)]{FOT11} that
    \begin{align*}
        \mathrm{Cap}_\theta(K)=\inf\left\{\mathcal{E}^\theta_1(\xi,\xi): \xi\in \mathscr{C}_c^\infty(\Omega,\R), \xi|_K\geq 1\right\}.
    \end{align*}
    For every $\xi\in \mathscr{C}_c^\infty(\Omega,\R)$ with $\xi|_K\geq 1$, by Theorem \ref{Theorem:Lu} and Lemma \ref{Poincare type inequality} we have
    \begin{align*}
        \mathrm{Cap}_{BT,m}(K, \Omega)&=\int_K(dd^cu^*_{K,\Omega,m})^m\wedge\omega^{n-m}\\
        &\leq\int_\Omega\xi^2(dd^cu^*_{K,\Omega,m})^m\wedge\omega^{n-m}\\
        &\leq \int_\Omega\xi^2 dd^cu^*_{K,\Omega,m}\wedge\theta\\
        &\leq 8 \lVert u^*_{K,\Omega,m}\rVert_{L^\infty(\Omega)}\mathcal{E}^\theta(\xi,\xi)\\
        &\leq 8 \mathcal{E}^\theta_1(\xi,\xi).
    \end{align*}
    Therefore,
    \begin{align*}
        \mathrm{Cap}_{BT,m}(K, \Omega)\leq 8 \mathrm{Cap}_\theta(K).
    \end{align*}

    (ii) Fix an arbitrary open subset $G\subset D$ and $v\in SH_m(\Omega)$ with $\lVert v\rVert_{L^\infty(\Omega)}\leq B$. Let $\xi\in \mathscr{C}_c^\infty(\Omega,\R)$ so that $0\leq\xi\leq 1$ and $\xi|_D=1$. By Theorem \ref{bounded msh in F loc} (i) we have $-\xi u^*_{G,\Omega,m}\in \mathcal{F}^\theta$. Notice that $-\xi u^*_{G,\Omega,m}=1$ on $G$. Hence, by (\ref{equation: Etheta and Etheta1}) we have
    \begin{align*}
        \mathrm{Cap}_\theta(G)\leq \mathcal{E}^\theta_1(-\xi u^*_{G,\Omega,m},-\xi u^*_{G,\Omega,m})\leq C_1 \mathcal{E}^\theta(-\xi u^*_{G,\Omega,m},-\xi u^*_{G,\Omega,m}),
    \end{align*}
    where $C_1=1+ \frac{8\gamma\pi^n}{2^n\delta^{m-1} n!}$. Then by Lemma \ref{Lemma: variant of Schwarz inequality for bdd m-sh functions} we have
    \begin{align}\label{bounded CapTheta by I1 and I2}
        \frac{1}{2C_1}\mathrm{Cap}_\theta(G)&\leq \int_\Omega\xi^2 du^*_{G,\Omega,m}\wedge d^cu^*_{G,\Omega,m}\wedge\theta+\int_\Omega (u^*_{G,\Omega,m})^2 d\xi\wedge d^c\xi\wedge\theta\\
        &=:I_1+I_2.\notag
    \end{align}
    
    We make the following preparations before estimating $I_1$ and $I_2$.
    
    Let $\rho$ be a $\mathscr{C}^2$ defining function of the bounded strictly $m$-pseudoconvex domain $\Omega$ and $C_2>0$ be a constant so that $(dd^c\rho)^k\wedge\omega^{n-k}\geq C_2 \omega^n$ on $\overline{\Omega}$ for $1\leq k\leq m$. We have the following estimates of two functions in $SH_m(\Omega)$.
    \begin{itemize}
        \item[(i')] $-2B-(2\delta+1)\gamma-1\leq v+ \delta|z|^2-B-(\delta+1)\gamma-1\leq -1$.
        \item[(ii')] $-2B-(2\delta+1)\gamma-1\leq |z|^2-B-(\delta+1)\gamma-1\leq -1$.
    \end{itemize}
    Hence, by Proposition \ref{Proposition:localization principle}, there exist an open set $U$, a compact set $K$, and a constant $C_3>0$ depending only on $\Omega, D, B,\delta$, that satisfy the following properties.
    \begin{itemize}
        \item[(i'')] $D\subset\supp\xi\subset U\subset K\subset\Omega$.
        \item[(ii'')] $\hat{v}:=\max\left(v+\delta|z|^2-B-(\delta+1)\gamma-1,C_3\rho\right)\in SH_m(\Omega)\cap L^\infty(\Omega)$,\\
        $\hat{w}:=\max\left(|z|^2-B-(\delta+1)\gamma-1,C_3\rho\right)\in SH_m(\Omega)\cap L^\infty(\Omega)$.
        \item[(iii'')] $\hat{v}=v+\delta|z|^2-B-(\delta+1)\gamma-1$ on $U$,\\
        $\hat{w}=|z|^2-B-(\delta+1)\gamma-1$ on $U$,\\
        $\hat{v}=\hat{w}=C_3\rho$ on $\Omega\setminus K$.
    \end{itemize}
    By Theorem \ref{Theorem:Lu} and the Chern--Levine--Nirenberg inequality,
    \begin{align*}
        \int_\Omega(dd^cu_{G,\Omega,m}^*)^m\wedge\omega^{n-m}=\mathrm{Cap}_{BT,m}(G,\Omega)<\infty.
    \end{align*}
    Since $\lim_{z\to\partial\Omega}u_{G,\Omega,m}^*(z)=\lim_{z\to\partial\Omega} \hat{v}(z)=\lim_{z\to\partial\Omega} \hat{w}(z)=0$, Cegrell's inequality \cite[Lemma 2.3]{KolodziejNguyen} and the comparison principle \cite[Corollary 1.15 (c)]{Nguyen13} imply
    \begin{align}\label{first finite estimate}
        &\int_\Omega dd^c u^*_{G,\Omega,m}\wedge(dd^c\hat{v})^{m-1}\wedge\omega^{n-m}\\
        &\leq \left(\int_\Omega(dd^cu_{G,\Omega,m}^*)^m\wedge\omega^{n-m}\right)^{\frac{1}{m}}\left(\int_\Omega(dd^c\hat{v})^m\wedge\omega^{n-m}\right)^{\frac{m-1}{m}}\notag\\
        &\leq C_3^{m-1}\left(\mathrm{Cap}_{BT,m}(G, \Omega)\right)^{\frac{1}{m}}\left(\int_\Omega(dd^c\rho)^m\wedge\omega^{n-m}\right)^{\frac{m-1}{m}}\notag\\
        &<\infty,\notag
    \end{align}
    \begin{align}\label{second finite estimate}
        &\int_\Omega dd^c \hat{w}\wedge(dd^c\hat{v})^{m-1}\wedge\omega^{n-m}\\
        &\leq \left(\int_\Omega(dd^c\hat{w})^m\wedge\omega^{n-m}\right)^{\frac{1}{m}}\left(\int_\Omega(dd^c\hat{v})^m\wedge\omega^{n-m}\right)^{\frac{m-1}{m}}\notag\\
        &\leq C_3^{m}\int_\Omega(dd^c\rho)^m\wedge\omega^{n-m}\notag\\
        &<\infty\notag.
    \end{align}
    
    We estimate $I_1$.
    \begin{align*}
        I_1&=\frac{1}{2}\int_\Omega \xi^2dd^c (u^*_{G,\Omega,m}+1)^2\wedge\theta-\int_\Omega \xi^2(u^*_{G,\Omega,m}+1) dd^c u^*_{G,\Omega,m}\wedge\theta\\
        &\leq \frac{1}{2}\int_\Omega \xi^2 dd^c (u^*_{G,\Omega,m}+1)^2\wedge\theta\\
        &=\frac{1}{2}\int_\Omega \xi^2 dd^c (u^*_{G,\Omega,m}+1)^2\wedge (dd^c \hat{v})^{m-1}\wedge \omega^{n-m}\\
        &\leq \frac{1}{2}\int_\Omega dd^c\left((u^*_{G,\Omega,m}+1)^2-1\right)\wedge(dd^c\hat{v})^{m-1}\wedge\omega^{n-m}.
    \end{align*}
    The two $m$-subharmonic functions satisfy $2u^*_{G,\Omega,m}\leq(u^*_{G,\Omega,m}+1)^2-1\leq 0$ and $\lim_{z\to\partial\Omega}\left((u^*_{G,\Omega,m}(z)+1)^2-1\right)=\lim_{z\to\partial\Omega}2u^*_{G,\Omega,m}(z)=0$. It follows from the comparison principle that
    \begin{align*}
        \int_\Omega \left(dd^c\left((u^*_{G,\Omega,m}+1)^2-1\right)\right)^m\wedge\omega^{n-m}\leq 2^m \int_\Omega (dd^c u^*_{G,\Omega,m})^m\wedge\omega^{n-m}<\infty.
    \end{align*}
    Hence, it follows from Cegrell's inequality and the comparison principle that
    \begin{align*}
        I_1&\leq\frac{1}{2}\left(\int_\Omega \left(dd^c\left((u^*_{G,\Omega,m}+1)^2-1\right)\right)^m\wedge\omega^{n-m}\right)^{\frac{1}{m}} \left(\int_\Omega(dd^c\hat{v})^m\wedge\omega^{n-m}\right)^{\frac{m-1}{m}}\\
        &\leq \left(\int_\Omega (dd^c u^*_{G,\Omega,m})^m\wedge\omega^{n-m}\right)^{\frac{1}{m}}\left(\int_\Omega(dd^c\hat{v})^m\wedge\omega^{n-m}\right)^{\frac{m-1}{m}}\\
        &\leq C_3^{m-1}\left(\mathrm{Cap}_{BT,m}(G, \Omega)\right)^{\frac{1}{m}}\left(\int_\Omega(dd^c\rho)^m\wedge\omega^{n-m}\right)^{\frac{m-1}{m}}.
    \end{align*}

    We estimate $I_2$. Let $C_4>0$ be a constant so that $d\xi\wedge d^c\xi\leq C_4 dd^c|z|^2$. Choose $\eta\in \mathscr{C}_c^\infty(U,\R)$ so that $0\leq \eta\leq 1$ and $\eta=1$ on an open neighborhood of $\supp\xi$. Notice that $(u^*_{G,\Omega,m})^2\leq -u^*_{G,\Omega,m}$. Then we have
    \begin{align*}
        I_2&\leq \int_\Omega (-u^*_{G,\Omega,m}) d\xi\wedge d^c\xi\wedge\theta\\ 
        &= \int_\Omega (-\eta u^*_{G,\Omega,m}) d\xi\wedge d^c\xi\wedge\theta\\
        &= \int_\Omega (-\eta u^*_{G,\Omega,m}) d\xi\wedge d^c\xi\wedge(dd^c\hat{v})^{m-1}\wedge\omega^{n-m}\\
        &\leq C_4 \int_\Omega (-\eta u^*_{G,\Omega,m}) d d^c|z|^2\wedge(dd^c\hat{v})^{m-1}\wedge\omega^{n-m}\\
        &= C_4 \int_\Omega (-\eta u^*_{G,\Omega,m}) d d^c \hat{w}\wedge(dd^c\hat{v})^{m-1}\wedge\omega^{n-m}\\
        &\leq C_4 \int_\Omega (- u^*_{G,\Omega,m}) d d^c \hat{w}\wedge(dd^c\hat{v})^{m-1}\wedge\omega^{n-m}.
    \end{align*}
    By the finite estimates of (\ref{first finite estimate}) and (\ref{second finite estimate}), we can then apply integration by parts \cite[Theorem 2.11]{Lu15JMAA}.
    \begin{align*}
        I_2&\leq C_4 \int_\Omega (- u^*_{G,\Omega,m}) d d^c \hat{w}\wedge(dd^c\hat{v})^{m-1}\wedge\omega^{n-m}\\
        &=C_4 \int_\Omega (-\hat{w})d d^c u^*_{G,\Omega,m} \wedge(dd^c\hat{v})^{m-1}\wedge\omega^{n-m}\\
        &\leq C_3 C_4 \lVert \rho\rVert_{L^\infty(\Omega)} \int_\Omega d d^c u^*_{G,\Omega,m} \wedge(dd^c\hat{v})^{m-1}\wedge\omega^{n-m}\\
        &\leq C_3^m C_4 \lVert \rho\rVert_{L^\infty(\Omega)}\left(\mathrm{Cap}_{BT,m}(G, \Omega)\right)^{\frac{1}{m}} \left(\int_\Omega(dd^c\rho)^m\wedge\omega^{n-m}\right)^{\frac{m-1}{m}}.
    \end{align*}

    Combining (\ref{bounded CapTheta by I1 and I2}) and the estimates of $I_1,I_2$, we obtain
    \begin{align*}
        \mathrm{Cap}_\theta(G)\leq C\left(\mathrm{Cap}_{BT,m}(G, \Omega)\right)^{\frac{1}{m}},
    \end{align*}
    where the constant
    \begin{align*}
        C=2C_1C_3^{m-1}\left(1+ C_3C_4 \lVert \rho\rVert_{L^\infty(\Omega)}\right)\left(\int_\Omega(dd^c\rho)^m\wedge\omega^{n-m}\right)^{\frac{m-1}{m}}
    \end{align*}
    depends only on $\Omega,D,B,\delta,m$.
\end{proof}

\begin{corollary}\label{Corollary: zero BT capacity iff zero capacity CapTheta}
    Let $\Omega\subset\C^n$ be a bounded strictly $m$-pseudoconvex domain, and $E\subset\Omega$ be a Borel subset. Then $\mathrm{Cap}_{BT,m}(E, \Omega) = 0$ if and only if $\mathrm{Cap}_\theta(E) = 0$ for every $\theta \in \mathcal{M}_m^+(\Omega)$.
\end{corollary}
\begin{proof}
    Since both $\mathrm{Cap}_{BT,m}$ and $\mathrm{Cap}_\theta$ are Choquet capacities, we can assume that $E=K$ is a compact set.
    
    If $\mathrm{Cap}_\theta(K) = 0$ for every $\theta \in \mathcal{M}_m^+(\Omega)$, by Lemma \ref{equivalence of caps} (i), $\mathrm{Cap}_{BT,m}(K, \Omega) = 0$.
    
    Assume that $\mathrm{Cap}_{BT,m}(K, \Omega) = 0$. Fix an arbitrary $\theta \in \mathcal{M}_m^+(\Omega)$ and write
    $$\theta=dd^c (v_1 +\delta |z|^2)\wedge\dots\wedge dd^c (v_{m-1}+\delta |z|^2)\wedge\omega^{n-m},$$
    where $v_1,\dots,v_{m-1}\in SH_m(\Omega)\cap L^\infty(\Omega)$ and $\delta>0$ is a constant. Set
    $$\Theta:= \left(dd^c\left(v_1+\dots+ v_{m-1}+\delta |z|^2\right)\right)^{m-1}\wedge\omega^{n-m}.$$
    Since $\theta\leq \Theta$, \cite[Lemma 2.2.7 (ii)]{FOT11} implies that $\mathrm{Cap}_\theta(K)\leq \mathrm{Cap}_\Theta(K)$. Let $D\Subset\Omega$ be an open subset containing $K$. By Lemma \ref{equivalence of caps} (ii), there exists a constant $C>0$ depending only on $\Omega,D,\Theta,m$, such that
    $$\mathrm{Cap}_\Theta(G)\leq C \left(\mathrm{Cap}_{BT,m}(G, \Omega)\right)^{\frac{1}{m}}$$
    for every open subset $G\subset D$.
    For every $\varepsilon>0$, there exists an open subset $G_{\varepsilon}\subset D$ such that $K\subset G_\varepsilon$ and $\mathrm{Cap}_{BT,m}(G_\varepsilon, \Omega) < \varepsilon$.  
    Therefore,
    \begin{align*}
        \mathrm{Cap}_\theta(K)\leq \mathrm{Cap}_\Theta(K)\leq \mathrm{Cap}_\Theta(G_\varepsilon)\leq C\left(\mathrm{Cap}_{BT,m}(G_\varepsilon, \Omega)\right)^{\frac{1}{m}} < C \varepsilon^{\frac{1}{m}}.
    \end{align*}
    By letting $\varepsilon\to 0$, we obtain $\mathrm{Cap}_\theta(K)=0$.
\end{proof}
\begin{remark}
    The proof of Corollary \ref{Corollary: zero BT capacity iff zero capacity CapTheta} is based on the quantitative estimates in Lemma \ref{equivalence of caps}. The implication from zero capacity $\mathrm{Cap}_{BT,m}$
    to zero capacity $\mathrm{Cap}_\theta$ can also be obtained by the qualitative characterization in Corollary \ref{corollart:every msh function is Etheta qc}. Indeed, let $\Omega\subset\C^n$ be a bounded open subset, and $E\subset\Omega$ be a subset. If $\mathrm{Cap}_{BT,m}(E,\Omega)=0$, then $E$ is $m$-polar. Since $V$ does not charge any $m$-polar set, Corollary \ref{corollart:every msh function is Etheta qc} implies
    \[
        \mathrm{Cap}_\theta(E)=0,\quad\forall\theta\in\mathcal{M}_m^+(\Omega).
    \]
\end{remark}

Our next goal is to give a stochastic interpretation of the Choquet capacity $\mathrm C_{\sharp,m}$. Let $\Omega\subset\C^n$ be a bounded open subset. Fix an arbitrary $\theta\in\mathcal{M}_m^+(\Omega)$. Let $(\mathcal{F}^\theta,\mathcal{E}^\theta)$ be the Dirichlet space for the admissible pair $(\theta,V)$. Since $\mathcal{E}^\theta$ is regular, it follows from \cite[Theorem~7.2.1]{FOT11} that  there exists a $V$-symmetric Hunt process associated with $\mathcal{E}^\theta$. This Dirichlet form $\mathcal{E}^\theta$ is also strong local by \cite[Exercise 3.1.1]{FOT11}. Thus, by \cite[Theorem 4.5.3]{FOT11}, we can choose such a process such that it is a diffusion  with no killing inside for every starting point in $\Omega$. Moreover, any such diffusion is a conformal martingale by \cite[Theorem 1]{FO84}. We fix any such conformal diffusion and denote it by
\begin{align*}
    M^{\theta}=\left((Z_t^\theta)_{t\in[0,\infty]},(P_z^\theta)_{z\in\Omega_\Delta},\zeta^\theta\right).
\end{align*}

We briefly explain some probabilistic concepts, notations and properties about this stochastic process, which will be used in our results. For more details including precise definition of Hunt process, we refer the reader to \cite[Appendix A.2]{FOT11}. Here $\Omega_\Delta$ is the state space of $M^{(\theta,V)}$,  where $\Delta\notin\Omega$ is a cemetery point and 
$\Omega_\Delta:=\Omega\cup\{\Delta\}$ endowed with the one-point
compactification topology. We can realize this process $M^\theta$ on a measurable space $(\mathsf{W},\mathcal{A})$, while for every fixed $\varpi \in\mathsf{W}$, the map $[0, \infty] \ni t \mapsto Z_t^{\theta}(\varpi)$ is a path in $\Omega_{\Delta}$. For every fixed $t$, the map
$$Z_t^\theta:(\mathsf{W}, \mathcal{A})\to\left(\Omega_\Delta,\mathcal{B}(\Omega_\Delta)\right)$$
is measurable. Here $\mathcal{B}(\Omega_\Delta)$ is the Borel $\sigma$-algebra of $\Omega_\Delta$. For every $z\in\Omega_\Delta$, $P_z^\theta$ is a probability measure on $(\mathsf{W},\mathcal{A})$ satisfying $P_z^\theta(Z_0^\theta=z)=1$. A statement holds $P_z^\theta$-almost surely ($P_z^\theta$-a.s. for short) if it fails only on some $A\in\mathcal{A}$ which is of zero measure $P_z^\theta$.

Let $X$ be a random variable on $\mathsf{W}$. Its expectation is $$E_z^\theta[X]:=\int_{\mathsf W}X dP_z^\theta,$$ whenever the integral above is well-defined.

The lifetime $\zeta^{\theta}$ of $M^{\theta}$ is defined by 
$$\zeta^\theta(\varpi):=\inf\{t\geq0:Z_t^\theta(\varpi)=\Delta\}.$$
For every $z\in\Omega$, we have
$P_z^\theta(\zeta^\theta>0)=1$. The definition of Hunt process gives the following properties.
For every $\varpi\in\mathsf{W}$, we have $Z^\theta_\infty(\varpi)=\Delta$ and $Z^\theta_t(\varpi)=\Delta$ when $t\geq \zeta^\theta(\varpi)$.
For every $\varpi\in\mathsf{W}$, the path
$t\mapsto Z_t^\theta(\varpi)$ is right-continuous at every time
$t\in[0,\infty)$ and has a left limit at every time $t\in(0,\infty)$. For every $t\in(0,\infty)$, we write
$Z_{t-}^\theta(\varpi):=\lim_{s\uparrow t}Z_s^\theta(\varpi)$.
These limits are taken in $\Omega_\Delta$.

The diffusion property means that
\begin{align}\label{equaton:diffusion}
    P^\theta_z(t\mapsto Z^\theta_t ~\text{is continuous on}~[0,\zeta^\theta))=1,\quad\forall z\in\Omega.
\end{align}
No killing inside for every starting point in $\Omega$ means that
\begin{align}\label{equation: no killing inside}
    P_z^\theta(Z^\theta_{\zeta^\theta-}\in\Omega,0<\zeta^\theta<\infty)=0,\quad\forall z\in\Omega.
\end{align}

Let $U\Subset\Omega$ be an open subset. Its first exit time is
$\tau_U:=\inf\{t\geq0:Z_t^\theta\notin U\}$. Since every path is right continuous, $Z^\theta_{\tau_U}\notin U$ and $\tau_U\leq\zeta^\theta$. Then
(\ref{equation: no killing inside}) implies that
\begin{align}\label{equation: no killing inside corollary}
    P_z^\theta(\tau_U<\zeta^\theta<\infty)=P_z^\theta(\zeta^\theta<\infty),\quad\forall z\in\Omega.
\end{align}
Let $U_j\Subset\Omega$ be an open subset of $\Omega$ for every $j\in\N$ so that
$\overline{U_j}\subset U_{j+1}$ for every $j\in\N$ and $\cup_{j\in\N}U_j=\Omega$. Then (\ref{equaton:diffusion}) implies that
\begin{align}\label{equation: diffusion corollary}
    P_z^{\theta}\left(\lim_{j \rightarrow \infty}\tau_{U_{j}} = \zeta^{\theta}\right) = 1,\quad \forall z\in\Omega.
\end{align}
Let $E\subset\Omega$ be a Borel subset. Its first hitting time is $\sigma_E:=\inf\{t>0:Z_t^\theta\in E\}$. Notice that $\sigma_E<\zeta^\theta$ on $\{\sigma_E<\infty\}$, and $Z_{\sigma_K}^\theta\in K$ on
$\{\sigma_K<\infty\}$ for every compact subset $K\subset\Omega$. The exit and hitting times depend on $\theta$, but we omit
this dependence from the notation. Both of them are stopping times \cite[Theorem A.2.3]{FOT11}

We recall the $\mathcal{E}^\theta$-superharmonicity introduced in \cite[Section 9]{FO87}. Let $u\in\mathcal{F}_{loc}^\theta$ and $\xi,\eta\in C_c^\infty(\Omega,\mathbb R)$. By Proposition \ref{proposition:characterization for local domain of Dirichlet space}, $\xi u\in\mathcal{F}^\theta$. Assume that $\xi=1$ on an open neighborhood of $\supp\eta$. Since $\mathcal E^\theta$ has local property, we can define
$$\mathcal E^\theta(u,\eta):=\mathcal E^\theta(\xi u,\eta).$$
We say a function $u\in \mathcal{F}_{loc}^\theta$ is $\mathcal{E}^\theta$-superharmonic on an open subset $U\subset\Omega$ if $\mathcal{E}^\theta(u,\eta)\geq 0$ for every non-negative function $\eta\in\mathscr{C}_c^\infty(U,\R)$. By Proposition \ref{proposition:characterization for local domain of Dirichlet space} and Theorem \ref{bounded msh in F loc}, for every $u\in SH_m(\Omega)\cap L_{loc}^\infty(\Omega)$, $-u$ is $\mathcal{E}^\theta$-superharmonic on $U$.

We now generalize \cite[Lemma 4]{FO84} to the setting of $m$-subharmonic functions. In the setting of plurisubharmonic functions, the proof of \cite[Lemma 4]{FO84} directly used Fukushima's decomposition theorem and the optional sampling theorem. Here we provide an alternative proof based on the super-mean-value inequality for $\mathcal{E}^\theta$-superharmonic functions in \cite{FO87}. In the proofs below, we write $a\wedge b :=\min(a,b)$.

\begin{lemma}\label{Lemma 5.13}
    Let $\theta \in \mathcal{M}_m^+(\Omega)$, and $M^\theta$ be a conformal diffusion associated with the Dirichlet form $\mathcal{E}^\theta$. Let $E \subset \Omega$ be a Borel subset. Then
    \begin{align*}
        \mathrm{C}_{\sharp, m}(E) \geq \int_\Omega P_z^{\theta}(\sigma_{E} < \infty)dV(z).
    \end{align*}
\end{lemma}

\begin{proof}
    Assume first $E=K$ is compact. Fix an arbitrary $u\in SH_m(\Omega)$ so that $-1\leq u\leq 0$ and $u|_K=-1$. By Theorem \ref{bounded msh in F loc}, $-u\in\mathcal F_{\mathrm{loc}}^\theta$ and is $\mathcal{E}^\theta$-quasi-continuous. Moreover, $-u$ is $\mathcal{E}^\theta$-superharmonic on $\Omega$. Let $U\Subset\Omega$ be an arbitrary open subset containing $K$. The Dirichlet form $\mathcal{E}^\theta$ satisfies (\ref{poincare inequality for lebesgue measure}), by \cite[Theorem 9.4 (ii)]{FO87},
    \begin{align*}
        -u(z) \geq 
        E_z^{\theta}\left(-u(Z_{\tau_U \wedge \sigma_K}^{\theta})\right)=\int_{\mathsf{W}}-u\left(Z_{\tau_U \wedge \sigma_K}^{\theta}\right)dP^\theta_z, \quad \mathcal{E}^\theta\text{-q.e. on~} \Omega.
    \end{align*}
    
    Although $Z^\theta_{\tau_U \wedge \sigma_K}\in \Omega_\Delta$ and $u$ is defined only on $\Omega$, the integral above is well-defined $\mathcal{E}^\theta$-q.e. on $\Omega$. Indeed, $V(\Omega)<\infty$, by \cite[Theorem 9.4 (i)]{FO87} we have
    \begin{align*}
        \int_{\mathsf{W}}\zeta^\theta d P_z^\theta<\infty,\quad\mathcal{E}^\theta\text{-q.e. on~} \Omega,
    \end{align*}
    which implies that
    \begin{align*}
        \tau_U\leq \zeta^\theta<\infty,\quad P^\theta_z\text{-a.s.,}\quad\mathcal{E}^\theta\text{-q.e. on~} \Omega.
    \end{align*}
    Combining this with (\ref{equation: no killing inside corollary}), we obtain
    \begin{align*}
        \tau_U\wedge\sigma_K\leq \tau_U< \zeta^\theta<\infty,\quad P^\theta_z\text{-a.s.,}\quad\mathcal{E}^\theta\text{-q.e. on~} \Omega.
    \end{align*}

    Since $V$ is bounded by $\mathrm{Cap}_\theta$ on every Borel subset of $\Omega$ and $\mathrm{Cap}_\theta$ is outer regular,
    \begin{align*}
        -u(z) \geq \int_{\mathsf{W}}-u\left(Z_{\tau_U \wedge \sigma_K}^{\theta}\right)dP^\theta_z\geq \int_{\{\sigma_K<\infty\}}-u\left(Z_{\tau_U \wedge \sigma_K}^{\theta}\right)dP^\theta_z,\quad V\text{-a.e. on~} \Omega.
    \end{align*}
    Let $U_j\Subset\Omega$ be an open subset of $\Omega$ for every $j\in\N$ so that $\overline{U_j}\subset U_{j+1}$ for every $j\in\N$, $\cup_{j\in\N}U_j=\Omega$, and $K\subset U_1$. Applying the above inequality to each $U_j$. Then there exists a subset $A\subset\Omega$ of zero measure $V$ such that
    \begin{align}\label{superharmonic inequality for exhausion}
        -u(z) \geq \int_{\{\sigma_K<\infty\}}-u\left(Z_{\tau_{U_j} \wedge \sigma_K}^{\theta}\right)dP^\theta_z,\quad \forall j\in\N, \quad \forall z\in \Omega\setminus A.
    \end{align}
    Fix an arbitrary $z\in\Omega\setminus A$. On $\{\sigma_K<\infty\}$, we have $\sigma_K<\zeta^\theta$. Then by (\ref{equation: diffusion corollary}),
    \begin{align*}
        \lim_{j\to\infty}\tau_{U_j}=\zeta^\theta,\quad P^\theta_z\text{-a.s. on~} \{\sigma_K<\infty\}.
    \end{align*}
    Therefore, for every $\varpi$ on $\{\sigma_K<\infty\}$ but outside a subset of zero measure $P^\theta_z$, there exists $j(\varpi)\in\N$ such that $\tau_{U_j}(\varpi)>\sigma_K(\varpi)$ whenever $j\geq j(\varpi)$. It then follows from the dominated convergence theorem that for every $z\in \Omega\setminus A$,
    \begin{align*}
         -u(z) &\geq \lim_{j\to\infty}\int_{\{\sigma_K<\infty\}}-u\left(Z_{\tau_{U_j} \wedge \sigma_K}^{\theta}\right)dP^\theta_z\\
         &= \int_{\{\sigma_K<\infty\}}\lim_{j\to\infty}-u\left(Z_{\tau_{U_j} \wedge \sigma_K}^{\theta}\right)dP^\theta_z\\
         &= \int_{\{\sigma_K<\infty\}}-u\left(Z_{\sigma_K}^{\theta}\right)dP^\theta_z\\
         &= P^\theta_z (\sigma_K<\infty).
    \end{align*}
    Integrating the inequality above on $\Omega\setminus A$ gives
    \begin{align}\label{equation: u and P}
        -\int_\Omega u(z) dV(z)\geq \int_\Omega P_z^{\theta}(\sigma_{K} < \infty)dV(z).
    \end{align}
    By Choquet's lemma, there exists an increasing sequence $(u_\ell)_{\ell\in\N}\subset SH_m(\Omega)$ such that $-1\leq u_\ell\leq 0$ and $u_\ell|_K=-1$ for every $\ell\in\N$, and $\lim_{\ell\to\infty}u_\ell=u^*_{K,\Omega,m}~V$-a.e. on $\Omega$. Applying (\ref{equation: u and P}) to each $u_\ell$ and then letting $\ell\to\infty$, we obtain 
    \begin{align*}
        \mathrm{C}_{\sharp, m}(K) \geq \int_{\Omega}P_z^{\theta}(\sigma_K < \infty)dV(z).
    \end{align*}   
    
    Now assume that $E\subset\Omega$ is a Borel subset. By \cite[Theorem A.2.4 (i)]{FOT11}, there exists an increasing sequence $(K_j)_{j \in\N}$ of compact subsets of $E$ such that 
    \begin{align*}
        V(\Omega)=\int_\Omega P^\theta_z\left( \lim_{j\to\infty }\sigma_{K_j}=\sigma_E\right)dV(z),
    \end{align*}
    which implies that 
    \begin{align*}
        \lim_{j\to\infty }\sigma_{K_j}=\sigma_E,\quad P^\theta_z\text{-a.s.,}\quad V\text{-a.e. on~}\Omega.
    \end{align*}
    The dominated convergence theorem and the case of compact subsets give
    \begin{align*}
        \int_\Omega P^\theta_z(\sigma_E<\infty)dV(z)&=\int_\Omega P^\theta_z\left(\cup_{j\in\N}\{\sigma_{K_j}<\infty\}\right)dV(z)\\
        &= \int_\Omega  \lim_{j\to\infty}P^\theta_z\left(\sigma_{K_j}<\infty\right) dV(z)\\
        &= \lim_{j\to\infty}\int_\Omega P^\theta_z\left(\sigma_{K_j}<\infty\right) dV(z)\\
        &\leq \lim_{j\to\infty} \mathrm{C}_{\sharp,m}(K_j)\\
        &\leq \mathrm{C}_{\sharp,m}(E).
    \end{align*}
    This finishes the proof.
\end{proof}

The next result generalizes \cite[Theorem 5.2]{FO87}.
\begin{lemma}\label{Lemma 5.14}
    Let $\Omega\subset\C^n$ be a bounded strictly $m$-pseudoconvex domain, $K\subset\Omega$ be a compact subset, and $\delta>0$ be a constant. Let $\theta := \left(dd^c (u^*_{K,\Omega,m}+\delta|z|^2)\right)^{m-1} \wedge \omega^{n-m}$, and $M^\theta$ be a conformal diffusion associated with the Dirichlet form $\mathcal{E}^\theta$. Then
    \begin{align*}
        \mathrm{C}_{\sharp, m}(K)  \leq \int_{\Omega}P_z^{\theta}(\sigma_K < \infty)dV(z) + (m-1)\delta  V(\Omega)\sup_\Omega |z|^2.
    \end{align*}
\end{lemma}
\begin{proof}
    Set $h:=u^*_{K,\Omega,m}$, $v:=h+\delta|z|^2$, and $\gamma:=\sup_\Omega |z|^2$. Compute that
    \begin{align}\label{computation}
        &dd^c h \wedge (dd^c v)^{m-1}\wedge \omega^{n-m}-(dd^c h)^m\wedge\omega^{n-m}\\
        &=\delta dd^c h\wedge dd^c|z|^2 \wedge\left(\sum_{j=0}^{m-2}(dd^c v)^{m-2-j}\wedge (dd^c h )^j\right) \wedge\omega^{n-m}\notag\\
        &\leq (m-1)\delta dd^c |z|^2 \wedge (dd^c v)^{m-1}\wedge\omega^{n-m}.\notag
    \end{align}
    It follows from Theorem \ref{Theorem:Lu} that $(dd^c h)^m\wedge\omega^{n-m}=0$ on $\Omega\setminus K$ in the sense of currents.  Hence (\ref{computation}) implies that $dd^c \left(h-(m-1)\delta|z|^2\right)\wedge\theta\leq 0$ on $\Omega\setminus K$. Then $w:=h-(m-1)\delta|z|^2\in \mathcal{F}^\theta_{loc}\cap L^\infty(\Omega)$ is $\mathcal{E}^\theta$-q.c., and $\mathcal{E}^\theta$-superharmonic on $\Omega\setminus K$. 

    Choose a function $\rho$ satisfying the condition in Definition \ref{Definition: bounded strictly m-psc domain} and $\rho|_K<-1$. 
    Let $U_j:=\{\rho<-\frac{1}{j}\}$ for every $j\in\N$. The Dirichlet form $\mathcal{E}^\theta$ satisfies (\ref{poincare inequality for lebesgue measure}), by \cite[Theorem 9.4 (ii)]{FO87}, there exists a subset $A\subset\Omega$ of zero measure $V$ such that
    \begin{align*}
         w(z) \geq E_z^{\theta}\left(w(Z^{\theta}_{\tau_{U_j} \wedge \sigma_K})\right),\quad \forall j\in\N,\quad \forall z\in \Omega\setminus A.
    \end{align*}
    Therefore,
    \begin{align}\label{good}
         -h(z) \leq -E_z^{\theta}\left(h(Z^{\theta}_{\tau_{U_j} \wedge \sigma_K})\right)+(m-1)\delta\gamma,\quad \forall j\in\N,\quad \forall z\in \Omega\setminus A.
    \end{align}
    
    We write
    \begin{align}\label{good5}
    E_z^{\theta}\left(h(Z^{\theta}_{\tau_{U_j} \wedge \sigma_K})\right)=\int_{\{\tau_{U_j}\leq \sigma_K\}}h(Z^{\theta}_{\tau_{U_j} }) dP^\theta_z + \int_{\{\tau_{U_j}> \sigma_K\}}h(Z^{\theta}_{\sigma_K}) dP^\theta_z.
    \end{align}

    Since $h\leq 0$, we have
    \begin{align}\label{good1}
        0\geq \int_{\{\tau_{U_j}\leq \sigma_K\}}h(Z^{\theta}_{\tau_{U_j} }) dP^\theta_z\geq \int_{\mathsf{W}}h(Z^{\theta}_{\tau_{U_j} }) dP^\theta_z.
    \end{align}
    Since $V(\Omega)<\infty$, by \cite[Theorem 9.4 (i)]{FO87}, there exists a subset $B\subset\Omega$ of zero measure $V$ such that
    \begin{align*}
        P^\theta_z(\zeta^\theta<\infty)=1,\quad \forall z\in \Omega\setminus B.
    \end{align*}
    Combining this with (\ref{equation: no killing inside corollary}), we obtain
    \begin{align}\label{good2}
        \tau_{U_j}<\zeta^\theta<\infty,\quad 
        \forall j\in\N, \quad P_z^\theta\text{-a.s.,}\quad \forall z\in\Omega\setminus B. 
    \end{align}
    For every $z\in\Omega\setminus B$, there exists $j(z)\in\N$ such that $z\in U_{j(z)}$. Hence, we have
    \begin{align}\label{good3}
        \tau_{U_j}>0,\quad \forall j\geq j(z),\quad P^\theta_z\text{-a.s.},\quad \forall z\in \Omega\setminus B.
    \end{align}
    Combining (\ref{good2}) with (\ref{good3}), we obtain
    \begin{align}\label{good4}
        0<\tau_{U_j}<\zeta^\theta<\infty,\quad \forall j\geq j(z),\quad P^\theta_z\text{-a.s.},\quad \forall z\in \Omega\setminus B.
    \end{align}
    Combining (\ref{good4}) with the continuous property of paths on $(0,\zeta^\theta)$ (\ref{equaton:diffusion}), we obtain
    \begin{align*}
        Z^\theta_{\tau_{U_j}}\in \partial U_j\quad \forall j\geq j(z),\quad P^\theta_z\text{-a.s.},\quad \forall z\in \Omega\setminus B.
    \end{align*}
    Since $0\geq h\geq \rho$, we then have
    \begin{align*}
        0\geq \liminf_{j\to\infty}h(Z^\theta_{\tau_{U_j}})\geq\lim_{j\to\infty}\rho(Z^\theta_{\tau_{U_j}}) =\lim_{j\to\infty}-\frac{1}{j} =0,\quad P^\theta_z\text{-a.s.},\quad \forall z\in \Omega\setminus B.
    \end{align*}
    Therefore, by the dominated convergence theorem,
    \begin{align*}
        \lim_{j\to\infty}\int_{\mathsf{W}}h(Z^{\theta}_{\tau_{U_j} }) dP^\theta_z=\int_{\mathsf{W}}\lim_{j\to\infty}h(Z^{\theta}_{\tau_{U_j} }) dP^\theta_z=0,\quad \forall z\in\Omega\setminus B.
    \end{align*}
    Combining this with (\ref{good1}), we obtain
    \begin{align}\label{good6}
        \lim_{j\to\infty}\int_{\{\tau_{U_j}\leq \sigma_K\}}h(Z^{\theta}_{\tau_{U_j} }) dP^\theta_z=0,\quad \forall z\in\Omega\setminus B.
    \end{align}

    On $\{\sigma_K<\infty\}$ we have $\sigma_K<\zeta^\theta$. Then by (\ref{equation: diffusion corollary}),
    \begin{align*}
        \{\sigma_K<\infty\}=\bigcup_{j=1}^\infty\{\tau_{U_j}>\sigma_K\},\quad P^\theta_z\text{-a.s.},\quad \forall z\in\Omega.
    \end{align*}
    The dominated convergence theorem gives
    \begin{align}\label{good7}
        \lim_{j\to\infty}\int_{\{\tau_{U_j}> \sigma_K\}}h(Z^{\theta}_{\sigma_K}) dP^\theta_z=\int_{\{ \sigma_K<\infty\}}h(Z^{\theta}_{\sigma_K}) dP^\theta_z,\quad \forall z\in\Omega.
    \end{align}
    Combining (\ref{good}), (\ref{good5}), (\ref{good6}), and (\ref{good7}), we obtain
    \begin{align*}
         -h(z) &\leq -\int_{\{ \sigma_K<\infty\}}h(Z^{\theta}_{\sigma_K}) dP^\theta_z+(m-1)\delta\gamma\\
         &\leq P^\theta_z(\sigma_K<\infty)+(m-1)\delta\gamma,\quad \forall z\in \Omega\setminus (A\cup B).
    \end{align*}
    Integrating both sides in $z\in\Omega$, we obtain the desired assertion.
\end{proof}

We now establish the stochastic characterization of the Choquet capacity $\mathrm{C}_{\sharp, m}$, generalizing \cite[Theorem 6.5]{FO87}.
\begin{theorem}\label{formula for Csharp}
    Let $\Omega \subset \C^n$ be a bounded strictly $m$-pseudoconvex domain, and $E \subset \Omega$ be a Borel subset.  For every $\theta\in \mathcal{M}_m^+(\Omega)$, let $M^{\theta}$ be a conformal diffusion associated with the Dirichlet form $\mathcal{E}^{\theta}$. Then
    \begin{align*}
        \mathrm{C}_{\sharp, m}(E) = \sup_{\theta \in \mathcal{M}_m^+(\Omega)}\int_{\Omega}P_z^{\theta}(\sigma_E < \infty)dV(z).
    \end{align*}
\end{theorem}
\begin{proof}
    One inequality follows from Lemma \ref{Lemma 5.13}.
    The reverse inequality follows from Lemma \ref{Lemma 5.14} when $E$ is compact, its general case follows from the compact case and the fact that $\mathrm C_{\sharp,m}$ is a Choquet capacity (Proposition \ref{proposition: Csharpm is Choquet}).
\end{proof}

We finally establish the following characterization of $m$-polar sets.
\begin{theorem}\label{Theorem:main result}
    Let $\Omega \subset \C^n$ be a bounded strictly $m$-pseudoconvex domain, and $E \subset \Omega$ be a Borel subset. For every $\theta\in \mathcal{M}_m^+(\Omega)$, let $M^{\theta}$ be a conformal diffusion associated with the Dirichlet form $\mathcal{E}^{\theta}$. Then the following statements are equivalent.
    \begin{itemize}
        \item [(i)] $u_{E, \Omega, m}^* \equiv 0$.
        \item [(ii)] $\mathrm{C}_{\sharp, m}(E) = 0$.
        \item [(iii)] $\mathrm{Cap}_{BT,m}(E, \Omega) = 0$.
        \item [(iv)] $\mathrm{Cap}_\theta(E) = 0$ for every $\theta \in \mathcal{M}_m^+(\Omega)$.
        \item [(v)] $P_z^{\theta}(\sigma_{E} < \infty) = 0$ $V$-a.e. on $\Omega$ for every $\theta \in \mathcal{M}_m^+(\Omega)$.       
    \end{itemize}
\end{theorem}
\begin{proof}
    (i) $\Leftrightarrow$ (iii) is from Theorem \ref{theorem: Lu's characterization}. (i) $\Leftrightarrow$ (ii) is from the definition of $\mathrm{C}_{\sharp, m}$. (iii) $\Leftrightarrow$ (iv) is from Corollary \ref{Corollary: zero BT capacity iff zero capacity CapTheta}. (ii) $\Leftrightarrow$ (v) is from Theorem \ref{formula for Csharp}. 
\end{proof}
\begin{remark}
    For every $\theta \in \mathcal{M}_m^+(\Omega)$, define $\mathrm{Cap}^0_{\theta}$ by replacing $\mathcal{E}_1^{\theta}$ with $\mathcal{E}^{\theta}$ in  Definition \ref{definition:CapTheta}. By (\ref{equation: Etheta and Etheta1}), for every subset $E \subset \Omega$,
    \begin{align*}
        \mathrm{Cap}_{\theta}(E) = 0 \Leftrightarrow \mathrm{Cap}^0_{\theta}(E) = 0.
    \end{align*}
    Hence, statement (iv) in the above theorem is equivalent to the following statement.
    \begin{itemize}
        \item[(vi)] $\mathrm{Cap}_\theta^0(E) = 0$ for every $\theta \in \mathcal{M}_m^+(\Omega)$.
    \end{itemize}
    
   By \cite[(9.11)]{FO87}, $\mathrm{Cap}_\theta^0$ can be expressed in terms of the hitting probability of the conformal diffusion associated with the Dirichlet form $\mathcal{E}^{\theta}$. More precisely, for every compact subset $K \subset \Omega$,
    \begin{align*}
        P_z^{\theta}(\sigma_K < \infty)=:e_K(z)\in \mathcal{F}^{\theta},\quad \mathrm{Cap}^0_{\theta}(K) = \mathcal{E}^{\theta}(e_K, e_K).
    \end{align*}
\end{remark}

\section*{Declaration of AI use}
After we completed this paper on our own, we used ChatGPT to check typos. We reviewed all suggested corrections and accepted some
of them.

\vspace{0.6cm}

\bibliographystyle{amsalpha}
\bibliography{ref.bib}

\vspace{1cm}

{\footnotesize

\begin{itemize}
    \item[] \textsc{Bowoo Kang\\Division of Mathematics, Department of Mathematics and Computer Science,\\University of Cologne\\ Weyertal 86-90, 50931 Cologne, Germany}\\
    \textit{E-mail address}: \texttt{bkang@uni-koeln.de}
\end{itemize}

\vspace{0.6cm}

\begin{itemize}
    \item[] \textsc{Jaechul Lee\\Department of Mathematical Sciences, KAIST\\291 Daehak-ro, Yuseong-gu, Daejeon 34141, Republic of Korea}\\
    \textit{E-mail address}: \texttt{jaechul1@kaist.ac.kr}
\end{itemize}

\vspace{0.6cm}

\begin{itemize}
    \item[] \textsc{Zhenghao Li\\Center for Complex Geometry, Institute for Basic Science\\ 55 Expo-ro, Yuseong-gu, Daejeon 34126, Republic of Korea}\\
    \textit{E-mail address}: \texttt{zhenghaoli@ibs.re.kr}
\end{itemize}
}
\end{document}